\documentclass[a4paper, 11pt]{amsart}
\usepackage[utf8]{inputenc}
\usepackage[english]{babel}
\usepackage{graphicx}
\usepackage{mathrsfs}
\usepackage{a4wide, comment}
\usepackage{xcolor}
\usepackage{amsfonts, amsmath, amsthm, amssymb}
\usepackage{enumitem, mathtools}
\usepackage{xcolor}
\usepackage{float}
\usepackage[numbers]{natbib}
\usepackage{xspace}
\usepackage[colorlinks=true,linkcolor=blue,urlcolor=blue,hypertexnames=false]{hyperref}

\allowdisplaybreaks[1]

\newtheorem{theorem}{Theorem}[section]

\newtheorem{assumption}[theorem]{Assumption}
\newtheorem{prop}[theorem]{Proposition}
\newtheorem{corollary}[theorem]{Corollary}
\newtheorem{lemma}[theorem]{Lemma}
\newtheorem{remark}[theorem]{Remark}

\newcommand{\I}{\mathcal{I}}

\def \E{\mathsf{E}}
\def \P{\mathsf{P}}
\def \R{\mathbb{R}}

\def \cF{\mathcal{F}}
\def \cG{\mathcal{G}}
\def \cI{\mathcal{I}}
\def \cB{\mathcal{B}}
\def \cL{\mathscr{L}}
\def \cT{\mathcal{T}}
\def \cO{\mathcal{O}}

\def \ve{\varepsilon}

\def \vol{b}

\def \lbracket{(}

\DeclareMathOperator*{\argmax}{arg\,max}
\newcommand\ind[1]{\mathbf{1}_{\{#1\}}}

\newcommand{\NJ}{N\!J}

\makeatletter 
\def\punctuationSpace{\@ifnextchar.{}{\@ifnextchar,{}{\@ifnextchar;{}{ }}}}
\makeatother

\newcommand{\cadlag}{c\`adl\`ag\punctuationSpace}
\newcommand{\caglad}{c\`agl\`ad\punctuationSpace}
\newcommand\as{\mbox{-a.s.}\punctuationSpace}

\makeatletter
\let\save@mathaccent\mathaccent
\newcommand*\if@single[3]{%
  \setbox0\hbox{${\mathaccent"0362{#1}}^H$}%
  \setbox2\hbox{${\mathaccent"0362{\kern0pt#1}}^H$}%
  \ifdim\ht0=\ht2 #3\else #2\fi
  }
\newcommand*\rel@kern[1]{\kern#1\dimexpr\macc@kerna}
\newcommand*\widebar[1]{\@ifnextchar^{{\wide@bar{#1}{0}}}{\wide@bar{#1}{1}}}
\newcommand*\wide@bar[2]{\if@single{#1}{\wide@bar@{#1}{#2}{1}}{\wide@bar@{#1}{#2}{2}}}
\newcommand*\wide@bar@[3]{%
  \begingroup
  \def\mathaccent##1##2{%
    \let\mathaccent\save@mathaccent
    \if#32 \let\macc@nucleus\first@char \fi
    \setbox\z@\hbox{$\macc@style{\macc@nucleus}_{}$}%
    \setbox\tw@\hbox{$\macc@style{\macc@nucleus}{}_{}$}%
    \dimen@\wd\tw@
    \advance\dimen@-\wd\z@
    \divide\dimen@ 3
    \@tempdima\wd\tw@
    \advance\@tempdima-\scriptspace
    \divide\@tempdima 10
    \advance\dimen@-\@tempdima
    \ifdim\dimen@>\z@ \dimen@0pt\fi
    \rel@kern{0.6}\kern-\dimen@
    \if#31
      \overline{\rel@kern{-0.6}\kern\dimen@\macc@nucleus\rel@kern{0.4}\kern\dimen@}%
      \advance\dimen@0.4\dimexpr\macc@kerna
      \let\final@kern#2%
      \ifdim\dimen@<\z@ \let\final@kern1\fi
      \if\final@kern1 \kern-\dimen@\fi
    \else
      \overline{\rel@kern{-0.6}\kern\dimen@#1}%
    \fi
  }%
  \macc@depth\@ne
  \let\math@bgroup\@empty \let\math@egroup\macc@set@skewchar
  \mathsurround\z@ \frozen@everymath{\mathgroup\macc@group\relax}%
  \macc@set@skewchar\relax
  \let\mathaccentV\macc@nested@a
  \if#31
    \macc@nested@a\relax111{#1}%
  \else
    \def\gobble@till@marker##1\endmarker{}%
    \futurelet\first@char\gobble@till@marker#1\endmarker
    \ifcat\noexpand\first@char A\else
      \def\first@char{}%
    \fi
    \macc@nested@a\relax111{\first@char}%
  \fi
  \endgroup
}
\makeatother

\DeclareFontFamily{U}{mathx}{\hyphenchar\font45}
\DeclareFontShape{U}{mathx}{m}{n}{<-> mathx10}{}
\DeclareSymbolFont{mathx}{U}{mathx}{m}{n}

\def\cadlag{c\`{a}dl\`{a}g\xspace}
\def\caglad{c\`{a}gl\`{a}d\xspace}

\title{Nonzero-Sum Stochastic Exit Games with Singular Belief Dynamics}
\author[Y.~Cao]{Yukun Cao}
\author[J.~Palczewski]{Jan Palczewski}
\address{Y.~Cao: School of Mathematics, University of Leeds, Woodhouse Lane, LS2 9JT Leeds, UK.}
\email{\href{mailto:mmyc@leeds.ac.uk}{mmyc@leeds.ac.uk}}
\address{J.~Palczewski: Faculty of Mathematics, Wroc\l{}aw University of Science and Technology, Wybrze\.{z}e Wyspia\'{n}skiego
27, 50-370, Wroc\l{}aw, Poland, and School of Mathematics, University of Leeds, Woodhouse Lane, LS2 9JT Leeds, UK.}\thanks{J.~Palczewski was partially supported by EPSRC grant no.~UKRI175.}
\email{\href{mailto:jan.palczewski@pwr.edu.pl}{jan.palczewski@pwr.edu.pl}}
\date{\today}

\numberwithin{equation}{section}

\begin{document}

\begin{abstract}
We study a two-player nonzero-sum stochastic exit game with asymmetric information. Each firm chooses the time to exit a duopoly and, conditional on remaining active, when to exit a monopoly. The asymmetry stems from privately observed exit values, while profit flows are driven by a general one-dimensional diffusion. Allowing for exit from monopoly induces singular dynamics of the associated belief process, so we construct it as a limit of regular approximations. Therefore, the underlying dynamics are not described by a system of SDEs, which requires a new mathematical approach. We construct a symmetric Nash equilibrium in which the duopoly exit is a hitting time of the belief process while the monopoly exit is described through the underlying diffusion. Our equilibrium exhibits novel phenomena: the market may be abandoned by both firms in finite time (market collapse), and firms with heterogeneous exit values may optimally exit simultaneously (clustering). These features are absent when monopoly is always preferred to exit.
\end{abstract}

\maketitle

\section{Introduction}

Technological progress and structural changes in the global economy often lead to the decline and eventual disappearance of once-prosperous industries. Firms operating in such environments face the fundamental strategic problem of determining whether and when to abandon a given line of business. In competitive settings, these decisions are inherently interactive: the timing of exit not only reflects underlying economic conditions but also depends on the anticipated behaviour of competing firms.

A game-theoretic approach to exit timing has therefore received considerable attention in the literature, beginning with the seminal works of Ghemawat and Nalebuff \cite{Ghemawat1985} and Fudenberg and Tirole \cite{fudenberg1986theory}, and extending to more recent contributions such as Georgiadis et al. \cite{georgiadis2022absence} and Kwon and Palczewski \cite{kwon}. In these models, firms operate in a duopoly and share profits that may deteriorate over time. Each firm chooses a stopping time at which it exits the market. The first firm to exit leaves the opponent as a monopolist. The transition to monopoly typically induces an increase in the remaining firm’s profit flow, creating a strategic trade-off between exiting early and remaining to capture monopoly rents.

Early contributions assume deterministic dynamics of profit flows, leading to equilibrium exit times determined through deterministic optimisation. Subsequent works introduce stochastic dynamics: in \cite{fine_equilibrium_1989} the model is in discrete time and driven by independent but not identically distributed random demand, in \cite{lambrecht2001,murto2004exit} the profit flow is given by a geometric Brownian motion, while in \cite{georgiadis2022absence,kwon} it is modelled as a continuous function of a one-dimensional diffusion. In these settings, equilibrium strategies are characterised using solutions to auxiliary optimal stopping problems. 

A further important extension, due to \cite{fudenberg1986theory}, introduces private information. When exiting the market, the firm may obtain a (possibly negative) exit value representing profits/costs from liquidating the line of business as well as expected profits from re-purposing freed facilities and human capital. While \cite{Ghemawat1985} assumes zero exit value, \cite{fudenberg1986theory} gives it a fundamental role in the model -- the exit value is random and only observed by the firm itself. This represents the uncertainty about the exact costs and opportunities of competing firms. Only the distributions of firms' exit values are common knowledge. This approach is further extended to a stochastic profit flow by \cite{kwon}.

A key assumption maintained in \cite{kwon} is that a firm strictly prefers operating as a monopolist rather than exiting the market. As a consequence, once a firm becomes a monopolist, it optimally remains active indefinitely. While analytically convenient, this assumption rules out the possibility that the market may eventually become unviable even for a monopolist, a phenomenon frequently observed in practice, e.g., the VHS industry worldwide or coal mining in many countries.

This paper contributes to the literature by explicitly constructing an equilibrium in a non-zero-sum exit game combining the three features presented above:
\begin{enumerate}[label=(\roman*)]
\item The dynamics of economic factors driving profit flows are described by a general one-dimensional diffusion.
\item The exit values (types) of firms are unknown to competitors, thereby creating information asymmetry.
\item The monopolist no longer remains active indefinitely but optimally chooses an exit time.
\end{enumerate}
Referring to the extant literature, the present paper extends \cite{fudenberg1986theory} by allowing for stochastic dynamics. It also generalises \cite{lambrecht2001,murto2004exit} by adding stochastic and private exit values and general dynamics of profit flows. It goes beyond \cite{kwon} where an optimally acting monopolist never exits the market. 

Compared to \cite{kwon}, allowing the monopolist to exit constitutes a single but structural change that leads to a qualitatively different class of stopping games. In \cite{kwon}, the strict preference for monopoly ensures that each firm faces a single stopping decision and yields a well-behaved belief process characterised as a Carathéodory solution to an ODE with a bounded right-hand side. This process reflects the inference of the distribution of opponents’ private types (exit values) and governs the exit from the duopoly. Once the monopolist's exit is allowed, the structure of strategies changes fundamentally: each firm must solve a double stopping problem, deciding when to exit both the duopoly and the monopoly. The monopolist’s value becomes endogenous, creating a feedback between the two phases of the game. As a consequence, the equilibrium behaviour exhibits features absent in \cite{kwon}: the game may terminate in finite time (industry collapse), exit decisions become highly sensitive to the proximity of the state to the monopoly exit region, and firms of different types may optimally exit simultaneously (clustering). These effects are reflected in the belief dynamics which lose the regularity and may exhibit singular features, including discontinuities and unbounded drift. The handling of these new belief dynamics requires a new approach to their definition as a limit of approximating processes. Consequently, the dynamics is not representable in an explicit integral form as in \cite{kwon}, thus requiring new and non-standard techniques for studying the equilibrium best response optimal stopping problems. In this sense, the paper does not merely extend the analysis of \cite{kwon}, but identifies and resolves a qualitatively different class of stochastic exit games. More generally, the paper develops analytical tools for studying asymmetric information games in which the belief dynamics are singular and do not admit explicit representations.

From a broader methodological perspective, this paper contributes to the growing literature on stopping games with asymmetric information. Much of this literature focuses on zero-sum settings, beginning with \cite{grun2013dynkin}, where asymmetry is driven by a privately observed finite-valued parameter and the state follows an SDE; in such models, equilibria typically require mixed (randomised) stopping times. Subsequent work considers alternative sources of asymmetry, including privately observed signals affecting the drift of a diffusion \cite{DeAngelis2020} and models based on privately observed finite-state Markov chains \cite{Gensbittel2019}. General existence results for zero-sum games with asymmetric information are provided in \cite{DAMP2021}, while \cite{deangelis2025} develops a martingale framework for characterising players' values and equilibrium strategies.
In contrast, non-zero-sum stopping games with asymmetric information are less developed and typically rely on guess-and-verify arguments, as in \cite{DeAngelis2020ghosts}. A distinct strand is formed by exit games with private exit values, including \cite{fudenberg1986theory}, \cite{kwon}, and the present paper. In these models, equilibria are constructed in pure strategies, although private information induces behaviour that, from the opponent’s perspective, resembles randomised stopping. A related mechanism also appears in the investment timing game of \cite{decamps2004investment}, where firms with privately known investment costs learn about a project of uncertain quality both from public news and from the rival's behaviour. As in the present paper, the equilibrium play is driven by the dynamics of beliefs about the opponent's private type.

Finally, decisions on managing involvement in declining industries have also been studied from different angles, leading to mathematical problems of a different nature. Relaxing the all-or-nothing nature of exit turns the war of attrition into one with variable levels of concession: in \cite{ghemawat1990devolution}, firms reduce capacity continuously and irreversibly, so that strategies are decreasing capacity paths rather than stopping times, while in \cite{reynolds1988plant,whinston1988exit} capacity is withdrawn plant by plant, which leads to multiple stopping problems. These models assume deterministic decline and complete information. Equilibria are constructed by backward induction from the date at which even a monopolist ceases to be profitable \cite{reynolds1988plant,whinston1988exit}, with continuous-time results obtained either directly \cite{reynolds1988plant} or as limits of discrete-time games when the period length tends to zero \cite{ghemawat1990devolution,whinston1988exit}. Without strategic interaction, exit under uncertainty is a classical optimal stopping problem \cite{dixit1989entry,kwon2010invest}.

The remainder of the paper develops the above ideas in a systematic way. Section \ref{sec:setup} introduces the model and reformulates the original double stopping non-zero sum game as a single stopping game. Section \ref{sec:single_player} analyses the associated optimal stopping problems for the duopoly and the monopoly, which serve as key building blocks for the equilibrium. Section \ref{sec:sketch} outlines the main ideas underlying the construction of the equilibrium and provides intuition for the role of the belief process. Section \ref{sec:belief} constructs the belief process and establishes its properties, including its behaviour under singular dynamics. Sections \ref{sec:value} and \ref{sec:proof-thm} characterise the optimal stopping regions for best response problems and verify the equilibrium condition. Finally, Section \ref{sec:example} presents numerical illustrations and highlights the qualitative differences relative to \cite{kwon}.

As a matter of convention throughout the paper, we write increasing/decreasing for non-strict monotonicity, and strictly increasing/decreasing for strict monotonicity.


\section{Problem setup}\label{sec:setup}
 Consider a market which initially supports two competing firms. They extract identical profit flows at the rate $D(X_t)$, where $(X_t)$ is a stochastic process observed by both firms and $D$ is a known function. Both firms face a decision whether to continue in the duopoly or exit the market. If firm $i$ exits the market, it receives a one-off payoff $\theta_i$. The other firm $j$ becomes a monopolist with the profit flow given by $M(X_t)$, where $M$ is a function known to both players and, naturally, $M > D$. When the monopolist market becomes too thin to support firm $j$, it exits and receives a one-off payoff $\theta_j$. The exit payoffs $\theta_i$ and $\theta_j$ are the same irrespective of whether the firm exits the duopoly or the monopoly market; they model the opportunities and costs of winding down the line of business.

Although both firms observe the same process $(X_t)$, the value of $\theta_i$, $i=1,2$, is not known to the other firm. We refer to it as the \emph{private type} or \emph{exit value} of firm $i$. The private type is random and its distribution, supported on $[\theta_L, \theta_U]$, is known to both firms. The realisation of $\theta_i$ is announced to firm $i$ at time $0$ while the competitor can only infer the value of $\theta_i$ from actions of firm $i$ as the private type is independent from $(X_t)$ and from the private type of the competitor.

Formally, there is a probability space $(\Omega, \cF,\P)$ with a filtration $(\cF_{t})$ satisfying the usual conditions supporting a one-dimensional $(\cF_t)$-Brownian motion $(W_t)$. The process $(X_{t})$ follows time-homogeneous dynamics
\begin{equation}\label{eqn:SDE}
    dX_{t}=\mu(X_{t})dt+\vol(X_{t})dW_{t}, \qquad X_0 = x.
\end{equation}
Denoting by $\cI = (x_L, x_U)$ the interval in which the process $(X_t)$ takes values, we assume that the functions $\mu, \vol:\cI \to \R$ are Lipschitz continuous and $\vol(\cdot) > 0$. So, for any initial value $x \in \cI$, there is a unique strong solution which is a strong Markov process. We further require that the boundaries $x_L, x_U$ of the interval  $\cI$ are not attainable for the process $(X_t)$. We allow one or both boundaries to be infinite.\footnote{We note that the Lipschitz continuity of $\mu, \vol$ can be relaxed as long as \eqref{eqn:SDE} has a unique strong solution for each $x \in \cI$, it defines a strong Markov process with non-attainable boundaries $x_L, x_U$ and the fundamental solutions are of class $C^2$. We adopt a more restrictive assumption of global Lipschitz continuity to avoid additional technicalities in the paper.}

The probability space $(\Omega, \cF, \P)$ supports additionally two random variables $\theta_{i}$, $i=1,2$, such that $\theta_1$, $\theta_2$ and $\cF_\infty$ are independent, so $(\theta_1, \theta_2)$ is also independent of the process $(X_t)$. These variables lead to two filtrations obtained as initial enlargements of filtration $(\cF_t)$:
\[
\cF^i_t = \cF_t \vee \sigma(\theta_i), \quad t \ge 0,\ i = 1,2.
\]
The filtration $(\cF^i_t)$ models the information flow of firm $i$. Note that it is right-continuous as $(\cF_t)$ is right-continuous and $\theta_i$ is independent.

A \emph{strategy} for firm $i$ is a pair of times $(\tau_i, \sigma_i)$ such that $\tau_i$ is an $(\cF^i_t)$-stopping time and $\sigma_i$ is an $(\hat\cF^i_t)$-stopping time, where 
\[
\hat\cF^i_t = \bigcap_{s > t} \cF^i_s \vee \sigma(\tau_1 \wedge \tau_2 \wedge s)
\]
is the smallest right-continuous complete filtration such that $\{\tau_1 \wedge \tau_2 \le u\} \in \hat\cF^i_t$, for all $u \le t$, and $\cF^i_t \subset \hat\cF^i_t$. At time $\tau_i$, firm $i$ exits the duopoly market, so at time $\tau_1 \wedge \tau_2$, which we assume to be observable without delay by both firms, the market shifts from duopoly to monopoly. The filtration $(\hat\cF^i_t)$ encompasses the additional knowledge whether the market is in the monopoly state. The role of $\sigma_i$ is to indicate the time when firm $i$, active in the monopoly market, exits. We therefore require that $\sigma_i \ge \tau_1 \wedge \tau_2$. 

Given a strategy profile $((\tau_1, \sigma_1), (\tau_2, \sigma_2))$, the expected payoff of firm 1 reads
\begin{multline}\label{eqn:bulky_payoffs}
J_1(x, (\tau_1, \sigma_1), (\tau_2, \sigma_2))\\ =  \E_{x}\bigg[\int^{\tau_{1}\wedge\tau_{2}}_{0}e^{-rt}D(X_{t})dt+\ind{\tau_1\leq\tau_2}e^{-r\tau_{1}}\theta_{1}+\ind{\tau_1>\tau_2}\Big(\int_{\tau_2}^{\sigma_1} e^{-rt} M(X_t) dt + e^{-r\sigma_{1}} \theta_1\Big)\bigg],
\end{multline}
where the subscript $x$ in $\E_x$ denotes the initial value of the process $(X_t)$. The expected payoff of firm 2 is analogous:
\begin{multline}\label{eqn:bulky_payoff2}
J_2(x, (\tau_1, \sigma_1), (\tau_2, \sigma_2))\\ =  \E_{x}\bigg[\int^{\tau_{1}\wedge\tau_{2}}_{0}e^{-rt}D(X_{t})dt+\ind{\tau_2\leq\tau_1}e^{-r\tau_{2}}\theta_{2}+\ind{\tau_2>\tau_1}\Big(\int_{\tau_1}^{\sigma_2} e^{-rt} M(X_t) dt + e^{-r\sigma_{2}} \theta_2\Big)\bigg].
\end{multline}

Our goal is to find a \emph{Nash equilibrium} in this non-zero-sum game with asymmetric information, i.e., for any $x \in \cI$, we want to find a strategy profile $((\tau^*_1, \sigma^*_1), (\tau^*_2, \sigma^*_2))$ such that for any other strategy profile $((\tau_1, \sigma_1), (\tau_2, \sigma_2))$ we have
\begin{equation}
\begin{aligned}
J_1(x, (\tau^*_1, \sigma^*_1), (\tau^*_2, \sigma^*_2)) &\ge J_1(x, (\tau_1, \sigma_1), (\tau^*_2, \sigma^*_2)),\label{eqn:nash1}\\
J_2(x, (\tau^*_1, \sigma^*_1), (\tau^*_2, \sigma^*_2)) &\ge J_2(x, (\tau^*_1, \sigma^*_1), (\tau_2, \sigma_2)).
\end{aligned}
\end{equation}

\subsection{Standing assumptions}

The analysis of this paper will be performed under the following standing assumptions which we collect here for the convenience of the reader.
\begin{assumption}$\ $\label{ass:main}
\begin{enumerate}[label=(\roman*)]
 \item Random variables $\theta_{i}$ have support $[\theta_{L},\theta_{U}]$ and a continuous cumulative distribution function $F$ which is strictly increasing in the interior of the support. The inverse of $F$ is Lipschitz continuous.
 \item Functions $D, M: \cI \to [0,\infty)$ are continuous, increasing and bounded, and furthermore $M>D$.
 \item The discount rate $r>0$.
\end{enumerate}
\end{assumption}

\begin{remark}
Condition (i) of the above assumption is satisfied if $F$ has a density that is separated from $0$ on $(\theta_L, \theta_U)$ by some constant $\delta > 0$.
\end{remark}

To formulate assumptions which involve the interplay between the dynamics of $(X_t)$ and the profit flow functions $D, M$, we introduce resolvents of $D$ and $M$:
\begin{equation}\label{eqn:resolvent}
\hat D(x) = \E_x \Big[ \int_0^\infty e^{-rt} D(X_t) dt \Big], 
\qquad
\hat M(x) = \E_x \Big[ \int_0^\infty e^{-rt} M(X_t) dt \Big].
\end{equation}

\begin{assumption}$\ $\label{ass:dm}
\begin{enumerate}[label=(\roman*)]
\item For each $\theta\in[\theta_{L},\theta_{U}]$, there exists a critical point $c_D(\theta) \in \overline{\cI}$ such that $D(x) < r\theta$ for $x<c_D(\theta)$ and $D(x)>r\theta$ for $x>c_D(\theta)$, where $x$ is always taken in $\cI$.
\item For each $\theta\in[\theta_{L},\theta_{U}]$, there exists a critical point $c_{M}(\theta)\in \overline{\cI}$ such that $M(x) \le r\theta$ for $x<c_{M}(\theta)$ and $M(x)>r\theta$ for $x>c_{M}(\theta)$, where $x$ is always taken in $\cI$.

\item Let $\phi$ be the decreasing fundamental solution to the ordinary differential equation $(\cL-r)u=0$, where
    \begin{equation*}
        \cL:=\frac{1}{2}\vol^{2}(x)\frac{\partial^{2}}{\partial x^{2}}+\mu(x)\frac{\partial}{\partial x}.
    \end{equation*}
    For each $\theta \in [\theta_L, \theta_U]$, define functions
    \begin{equation*}
        a_{\theta}(x):=\frac{\theta-\hat D(x)}{\phi(x)}, \qquad  b_{\theta}(x):=\frac{\theta-\hat M(x)}{\phi(x)}, \qquad x \in \cI.
    \end{equation*}
    Assume that they can be extended to $a_\theta, b_\theta:\overline{\cI} \to \overline\R$ in a continuous way, where $\overline\R = \R \cup \{-\infty, +\infty\}$, and that they attain unique global maxima at $\alpha(\theta), \beta(\theta) \in \overline{\cI}$, respectively. 
\end{enumerate}
\end{assumption}

In conditions (i) and (ii) above, the critical points \( c_D(\theta) \) and \( c_M(\theta) \) are allowed to take the boundary values $x_L$ or $x_U$. The phrase ``where $x$ is always taken in $\cI$'' means that if either \( c_D(\theta) \) or \( c_M(\theta) \) equals \( x_L \), then the region \( x < x_L \) is empty and the condition in that region is automatically satisfied; similarly, if either equals \( x_U \), then the condition in the region \( x > x_U \) is satisfied. 
The strict inequality $D(x) < r \theta$ for $x < c_D(\theta)$ is required for arguments involving the belief process $(Z_t)$ defining the equilibrium (see Lemma \ref{lem:stop-Z-transform}); it is equivalent to strict monotonicity at points $x$ for which $D(x) \in [r\theta_L, r\theta_U]$. Such strict monotonicity is not required for the function $M$. Condition (ii) merely defines $c_M(\theta)$ as the largest $x$ for which $M(x) \le r\theta$. Sufficient conditions for the continuous extension of $a_\theta, b_\theta$ to $\overline\cI$ in (iii) are provided in Remark \ref{rem:a_theta}. In particular, it is shown that $a_\theta(x_L) = b_\theta(x_L) = 0$.

Assumption \ref{ass:dm} ensures that the following two optimal stopping problems have solutions of threshold type and allows us to study those thresholds in more detail. For $(x, \theta) \in \cI \times [\theta_L, \theta_U]$, let
\begin{equation}\label{eqn:stop_dm}
 \begin{aligned}
    m(x;\theta)&=\sup_{\sigma \in \cT(\cF_t)}\E_{x}\Big[\int^{\sigma}_{0}e^{-rs}M(X_{s})ds+e^{-r\sigma}\theta\Big],\\
    d(x;\theta)&=\sup_{\sigma \in \cT(\cF_t)}\E_{x}\Big[\int^{\sigma}_{0}e^{-rs}D(X_{s})ds+e^{-r\sigma}\theta\Big],
 \end{aligned}
\end{equation}
where $\cT(\cF_t)$ denotes the set of $(\cF_t)$-stopping times. The value functions and optimal boundaries for these stopping problems provide the main building blocks for an explicit construction of a Nash equilibrium for our game.

We list the main properties of the above value functions below and defer their proof and more detailed analysis to Section \ref{sec:single_player}.
\begin{prop}\label{prop:dm}
Recall functions $\alpha, \beta:[\theta_L, \theta_U] \to \overline \cI$ from Assumption \ref{ass:dm}.
\begin{enumerate}[label=(\roman*)]
  \item Functions $m, d: \cI \times [\theta_L, \theta_U] \to \R$ are (jointly) continuous.
   \item For each $\theta \in [\theta_L, \theta_U]$, the stopping set for $m(x; \theta)$ is $\{x \in \cI: x \le \beta(\theta) \}$, and the function $\beta$ is continuous, increasing and satisfies $\beta(\theta) \le c_M(\theta)$.
  \item For each $\theta \in [\theta_L, \theta_U]$, the stopping set for $d(x; \theta)$ is $\{x \in \cI: x \le \alpha(\theta) \}$, and the function $\alpha$ is continuous, increasing and satisfies $\alpha(\theta) \le c_D(\theta)$.
  \item Functions $m, d$ are continuously differentiable in $x$ (i.e., the smooth fit holds), of class $C^2$ in $x$ away from their stopping boundaries $\beta(\cdot), \alpha(\cdot)$ and, for almost every $x \in \cI$, 
\begin{equation}\label{eqn:a1}
(\cL - r) m(x; \theta) + M(x) \le 0, \qquad (\cL - r) d(x; \theta) + D(x) \le 0,
\end{equation}
with the equality for $x > \beta(\theta)$ and $x > \alpha(\theta)$, respectively.
\end{enumerate}
\end{prop}

\subsection{Derivation of the stopping game}\label{subsec:stopping-game}
The strategy profile (and the Nash equilibrium) has a complex structure comprising a pair of stopping times with respect to different filtrations. In this subsection, we will show that $(\tau^*_1, \tau^*_2)$ in the Nash equilibrium is also an equilibrium of a non-zero-sum stopping game and $(\sigma^*_1, \sigma^*_2)$ are solutions to optimal stopping problems. 

Observe that the payoff of firm 1 in \eqref{eqn:bulky_payoffs} does not depend on the second element $\sigma_2$ of firm 2's strategy and, likewise, the payoff of firm 2 in \eqref{eqn:bulky_payoff2} does not depend on $\sigma_1$. Conditions \eqref{eqn:nash1} lead us to the following statement about $\sigma_1$ and $\sigma_2$ in equilibrium. An easy proof is omitted.
\begin{lemma}\label{lem:Nash_sigma}
In a Nash equilibrium $((\tau^*_1, \sigma^*_1), (\tau^*_2, \sigma^*_2))$, on the set $\{\tau^*_1 > \tau^*_2\}$, the stopping time $\sigma^*_1$ coincides with a solution to the optimal stopping problem
\begin{equation}\label{eqn:sigma1}
\sup_{\sigma \ge \tau^*_1 \wedge \tau^*_2,\ \sigma \in \cT(\hat\cF^1_t)} \E_x \Big[ \ind{\tau^*_1 > \tau^*_2} \Big(\int_{\tau^*_1 \wedge \tau^*_2}^{\sigma} e^{-rt} M(X_t) dt + e^{-r\sigma} \theta_1 \Big)\Big],
\end{equation}
and, on the set $\{\tau^*_2 > \tau^*_1\}$, the stopping time $\sigma^*_2$ coincides with a solution to the optimal stopping problem
\begin{equation}\label{eqn:sigma2}
\sup_{\sigma \ge \tau^*_1 \wedge \tau^*_2,\ \sigma \in \cT(\hat\cF^2_t)} \E_x \Big[ \ind{\tau^*_2 > \tau^*_1} \Big( \int_{\tau^*_1 \wedge \tau^*_2}^{\sigma} e^{-rt} M(X_t) dt + e^{-r\sigma} \theta_2 \Big) \Big].
\end{equation}
\end{lemma}

In other words, once the players select $\tau^*_1, \tau^*_2$, the stopping times $\sigma^*_1, \sigma^*_2$ are automatically determined. This observation suggests that the game can be reformulated as a game with a single decision variable per player: $\tau_1$ and $\tau_2$. The next few results will help establish this formulation.

Recalling that $(X_t)$ is a Markov process, the stopping problems \eqref{eqn:sigma1}-\eqref{eqn:sigma2} seem to be closely related to the optimal stopping problems $m(X_{\tau^*_1\wedge\tau^*_2}; \theta_1)$ and $m(X_{\tau^*_1\wedge\tau^*_2}; \theta_2)$, respectively, where $m$ is defined in \eqref{eqn:stop_dm}. In order to state formally the relationship between \eqref{eqn:stop_dm} and \eqref{eqn:sigma1}-\eqref{eqn:sigma2}, and their respective solutions, we need to describe in more detail the structure of $(\cF^i_t)$-stopping times (Lemma \ref{lem:structure_Fi_times}) and the relationship between filtrations $(\cF_t)$, $(\cF^i_t)$ and $(\hat\cF^i_t)$ (Lemma \ref{lem:immersion}).

\begin{lemma}\label{lem:structure_Fi_times}
\leavevmode
\begin{enumerate}[label=(\roman*)]
\item If the mapping $\hat{\tau}:(\Omega,\cF)\otimes([\theta_{L},\theta_{U}],\cB([\theta_{L},\theta_{U}]))\to([0,\infty],\cB([0,\infty]))$ is measurable and 
\begin{equation}\label{eqn:t}
\{(\omega, z) \in \Omega \times [\theta_L, \theta_U]:\ \hat{\tau}(\omega, z) < t \} \in \cF_t \otimes \cB([\theta_{L},\theta_{U}])
\end{equation}
for all $t > 0$, then $\tau=\hat{\tau}(\cdot,\theta_{i})$ is an $(\cF_{t}^{i})$-stopping time. Furthermore, $\hat\tau(\cdot, \theta)$ is an $(\cF_t)$-stopping time for every $\theta \in [\theta_L, \theta_U]$.
\item If $\tau$ is an $(\cF_{t}^{i})$-stopping time, then there exists $\hat \tau$ as in (i) such that $\tau = \hat{\tau}(\cdot,\theta_{i})$.
\end{enumerate}
\end{lemma}
\begin{proof}
Fix $i \in \{1,2\}$. Define $\Psi(\omega) = (\omega, \theta_i(\omega))$. Then $\cF^i_t = \Psi^{-1}\big(\cF_t \otimes \cB([\theta_{L},\theta_{U}])\big)$, cf. \cite[Lemma~1.4]{Kallenberg2002}. We proceed with the proof of (i). Let $\tau=\hat{\tau}(\cdot,\theta_{i})$ and take any $t \ge 0$. We have 
\[
\{\tau < t\} = \{\hat{\tau}(\cdot,\theta_{i}) < t\} = \Psi^{-1} \left( \{(\omega, z) \in \Omega \times [\theta_L, \theta_U]:\ \hat{\tau}(\omega, z) < t \} \right) \in \cF^i_t
\]
by \eqref{eqn:t}. Due to the right-continuity of $(\cF^i_t)$, this implies that $\tau$ is an $(\cF^i_t)$-stopping time. Noting that $\hat\tau(\cdot, \theta)$ is a section of the set in \eqref{eqn:t}, we obtain $\{ \hat\tau(\cdot, \theta) < t \} \in \cF_t$, which, by the right-continuity of $(\cF_t)$ implies that $\hat\tau(\cdot, \theta)$ is an $(\cF_t)$-stopping time.

Assume now that $\tau$ is an $(\cF^i_t)$-stopping time. For any rational $q > 0$, the set $\{ \tau < q \} \in \cF^i_q$, so by the comment at the beginning of the proof, $\{ \tau < q \} = \{ \omega \colon (\omega, \theta_i(\omega)) \in B_q \}$ for some set $B_q \in \cF_q \otimes \cB([\theta_{L},\theta_{U}])$. Define $\hat\tau(\omega, \theta) = \inf\{q \in \mathbb{Q}\cap (0, \infty)\colon (\omega, \theta) \in B_q \}$. Then we have
\[
\hat\tau(\omega, \theta_i(\omega)) = \inf\{q \in \mathbb{Q}\cap (0, \infty)\colon (\omega, \theta_i(\omega)) \in B_q \}
= \inf\{q \in \mathbb{Q}\cap (0, \infty)\colon \tau(\omega) < q \} = \tau(\omega),
\]
where the second equality is by the definition of $B_q$ and with the convention that $\inf \emptyset = \infty$. To show \eqref{eqn:t}, fix $t > 0$ and notice that 
\begin{align*}
\{(\omega,z) \colon \hat{\tau}(\omega, z) < t\} 
&= \{ (\omega,z)\colon (\omega, z) \in B_q \text{ for some $q < t$}\}\\
&= \bigcup_{q \in \mathbb{Q} \cap (0, t)} B_q \in \cF_t \otimes \cB([\theta_{L},\theta_{U}]).
\end{align*} 
\end{proof}

The above lemma generalises Proposition~1 in \cite{kwon} which relies on the product form of the probability space. The analogue of condition \eqref{eqn:t} there employs a non-strict inequality and the method of proof is different.

A filtration $(\cG^1_t)$ is said to be \emph{immersed} into a filtration $(\cG^2_t)$ if any $(\cG^1_t)$-martingale is a $(\cG^2_t)$-martingale (cf., \cite[Def.~3.1]{aksamit2017}).

\begin{lemma}\label{lem:immersion}
Filtration $(\cF_t)$ is immersed into $(\cF^i_t)$ and filtration $(\cF^i_t)$ is immersed into $(\hat\cF^i_t)$, $i=1,2$.
\end{lemma}
\begin{proof}
We start the proof with an auxiliary result: for $\sigma$-algebras $\cG', \cG'', \cG^\bot$ with $\cG' \subset \cG''$ and $\cG''$ independent from $\cG^\bot$, we have
\begin{equation}\label{eqn:cond_indep}
\E[Z | \cG'] = \E[Z | \cG'']
\end{equation}
for any $\cG' \vee \cG^\bot$-measurable and integrable random variable $Z$. The statement is immediate for $Z = Z' Z^\bot$, where $Z'$ is $\cG'$-measurable and $Z^\bot$ is $\cG^\bot$-measurable. Indeed, $\E[Z | \cG'] = \E[Z' Z^\bot| \cG'] = Z' \E[Z^\bot] = \E[Z' Z^\bot|\cG'']$ because $Z^\bot$ is independent from $\cG''$ and $Z'$ is $\cG''$-measurable. The equality is extended to all bounded $\cG' \vee \cG^\bot$-measurable random variables $Z$ by the Monotone Class Theorem \cite[Thm.~II.3.1]{RogersWilliams1}. A further extension to integrable $Z$ is achieved via the Monotone Convergence Theorem applied to truncations of positive and negative parts of $Z$.

To prove immersion of $(\cF_t)$ into $(\cF^i_t)$, we use an equivalent condition shown in \cite[Thm.~3.2]{aksamit2017}: the immersion holds iff $\E[Z | \cF_t] = \E[Z|\cF_\infty]$ for any integrable $(\cF^i_t)$-measurable random variable $Z$. This follows from \eqref{eqn:cond_indep} with $\cG' = \cF_t$, $\cG'' = \cF_\infty$ and $\cG^\bot = \sigma(\theta_i)$. We note that the immersion of $(\cF_t)$ into $(\cF^i_t)$ can also be deduced from \cite[Prop.~1.21]{aksamit2017} for which there is only a sketch of the proof.

For the immersion property of $(\cF^i_t)$ into $(\hat \cF^i_t)$, we first notice that $\hat \cF^i_t \subset \tilde \cF_t := \cF^i_t \vee \sigma(\theta_j)$, where $j = 3-i$ is the index of the remaining firm, because $\tau_j$ is an $(\cF^j_t)$-stopping time and $\cF^j_t = \cF_t \vee \sigma(\theta_j) \subset \tilde \cF_t$. Due to the independence of $\theta_j$ from $\cF_\infty$ and from $\theta_i$, we can use the same arguments as in the above paragraph to show that the filtration $(\cF^i_t)$ is immersed into $(\tilde \cF_t)$. From $(\cF^i_t) \subset (\hat \cF^i_t) \subset (\tilde\cF_t)$, \cite[Prop.~1.24\lbracket a)]{aksamit2017} yields that $(\cF^i_t)$ is immersed into $(\hat \cF^i_t)$.
\end{proof}

The above two lemmas pave the way for linking the optimal stopping problems \eqref{eqn:sigma1}-\eqref{eqn:sigma2} with $m(x; \theta_1)$ and $m(x; \theta_2)$.

\begin{lemma}\label{lem:equiv_to_m}
For a strategy profile $((\tau_1, \sigma_1), (\tau_2, \sigma_2))$, we have
\begin{align*}
J_1(x, (\tau_1, \sigma_1^*), (\tau_2, \sigma_2)) &\ge J_1(x, (\tau_1, \sigma_1), (\tau_2, \sigma_2)),\\
J_2(x, (\tau_1, \sigma_1), (\tau_2, \sigma^*_2)) &\ge J_2(x, (\tau_1, \sigma_1), (\tau_2, \sigma_2)),
\end{align*}
where
\begin{equation}\label{eqn:sigma_star}
\sigma_i^* = \inf\{t \ge \tau_1 \wedge \tau_2: X_t \le \beta(\theta_i) \} , \quad i=1,2,
\end{equation}
$\beta(\theta)$ is the stopping boundary from Proposition \ref{prop:dm}(ii) corresponding to the optimal stopping problem $m(x; \theta)$.
Furthermore,
\[
J_1(x, (\tau_1, \sigma_1^*), (\tau_2, \sigma_2)) =  \E_{x}\bigg[\int^{\tau_{1}\wedge\tau_{2}}_{0}e^{-rt}D(X_{t})dt+1_{\tau_1\leq\tau_2}e^{-r\tau_{1}}\theta_{1}+1_{\tau_1>\tau_2}e^{-r\tau_2} m(X_{\tau_2}; \theta_1)\bigg],
\]
and
\[
J_2(x, (\tau_1, \sigma_1), (\tau_2, \sigma_2^*)) =  \E_{x}\bigg[\int^{\tau_{1}\wedge\tau_{2}}_{0}e^{-rt}D(X_{t})dt+1_{\tau_2\leq\tau_1}e^{-r\tau_{2}}\theta_{2}+1_{\tau_2>\tau_1} e^{-r\tau_1} m(X_{\tau_1}; \theta_2) \bigg].
\]
\end{lemma}
\begin{proof}
We only consider $J_1$, since the arguments for $J_2$ are analogous.
We first show that the process 
\[
U_t = \int_0^t e^{-ru} M(X_u) du + e^{-rt} m(X_t; \theta_1), \qquad t \ge 0,
\]
is an $(\hat\cF^1_t)$-supermartingale for any $\P_x$. We start by proving that $(U_t)$ is an $(\cF^1_t)$-supermartingale for any $\P_x$ and then use the immersion property of $(\cF^1_t)$ into $(\hat\cF^1_t)$. Fix $x \in \cI$ and consider the increment for $t_2 > t_1$. By the Monotone Class Theorem, the supermartingale condition $\E_x [ U_{t_2} - U_{t_1} | \cF^1_{t_1}] \le 0$, $\P_x$\as, is equivalent to showing that $\E_x [ \ind{A} \ind{B} (U_{t_2} - U_{t_1}) ] \le 0$ for any $A \in \cF_{t_1}$ and $B \in \sigma(\theta_1)$. This is obviously implied by
\begin{equation}\label{eqn:before_freeze}
\E_x [ \ind{A} (U_{t_2} - U_{t_1}) | \theta_1 ] \le 0, \qquad \P_x\as,
\end{equation}
for any $A \in \cF_{t_1}$, the inequality that we will show now. 

By the Freezing Lemma \cite[Lemma~4.1]{baldi2017} the conditional expectation on the left-hand side of \eqref{eqn:before_freeze} equals $\varphi(\theta_1)$, where
\[
\varphi(\theta) 
= 
\E_x \bigg[\ind{A} \bigg(\int_{t_1}^{t_2} e^{-ru} M(X_u) du + e^{-rt_2} m(X_{t_2}; \theta) - e^{-rt_1} m(X_{t_1}; \theta) \bigg)\bigg].
\]
Thanks to the smoothness of $m$ asserted in Proposition \ref{prop:dm}, we can apply generalised It\^o's formula (\cite[Ch.~3, Problem~7.3]{karatzas2012brownian}) for any fixed value $\theta$ and get
\begin{multline}\label{eqn:a2}
e^{-rt_2} m(X_{t_2}; \theta) - e^{-rt_1} m(X_{t_1}; \theta)\\
= \int_{t_1}^{t_2} e^{-ru} (\cL - r)m(X_u; \theta) du  + \int_{t_1}^{t_2} e^{-ru} \frac{\partial}{\partial x} m(X_u; \theta) \vol(X_u) dW_u.
\end{multline}
Indeed, $m$ is continuously differentiable with the second derivative potentially discontinuous at the stopping boundary $\beta(\theta)$ but locally bounded. 

Let $(\eta_n)$ be the localising sequence for the local martingale $s \mapsto \int_{t_1}^{s} e^{-ru} \frac{\partial}{\partial x} m(X_u; \theta) \vol(X_u) dW_u$, $s \ge t_1$. By the dominated convergence theorem (as $m$ is continuous and bounded, and $M$ is bounded) and then using \eqref{eqn:a2}, we have
\begin{align}
\varphi(\theta)
&=
\lim_{n \to \infty} \E_x \bigg[\ind{A} \bigg(\int_{t_1}^{t_2\wedge \eta_n} e^{-ru} M(X_u) du + e^{-r(t_2 \wedge \eta_n)} m(X_{t_2\wedge \eta_n}; \theta) - e^{-rt_1} m(X_{t_1}; \theta) \bigg)\bigg]\notag\\
&=
\lim_{n \to \infty} \E_x \bigg[\ind{A} \int_{t_1}^{t_2 \wedge \eta_n} e^{-ru} \big( M(X_u) + (\cL - r)m(X_u; \theta) \big) du\notag\\
&\hspace{55pt}+ \ind{A} \int_{t_1}^{t_2 \wedge \eta_n} e^{-ru} \frac{\partial}{\partial x} m(X_u; \theta) \vol(X_u) dW_u \bigg]\notag\\
&=
\lim_{n \to \infty} \E_x \bigg[\ind{A} \int_{t_1}^{t_2 \wedge \eta_n} e^{-ru} \big( M(X_u) + (\cL - r)m(X_u; \theta) \big) du \bigg],\label{eqn:e1}
\end{align}
where in the last equality we used the martingale property of the stochastic integral and the fact that $A \in \cF_{t_1}$. Inequality \eqref{eqn:a1} in Proposition \ref{prop:dm} shows that the expression under the limit in \eqref{eqn:e1} is non-positive for each $n$, so $\varphi(\theta) \le 0$ for any $\theta \in [\theta_L, \theta_U]$. This implies \eqref{eqn:before_freeze} and, consequently, that $(U_t)$ is an $(\cF^1_t)$-supermartingale for any $\P_x$.

We extend the supermartingale property to the filtration $(\hat\cF^1_t)$ using the immersion property. Notice that the supermartingale $(U_t)$ is uniformly integrable; the integral is bounded (in absolute value) by $\sup_x |M(x)| / r$ and the function $|m|$ is bounded by $\sup_x |M(x)| / r + \max(|\theta_U|, |\theta_L|)$. It, therefore, has a Doob-Meyer decomposition as the sum of an $(\cF^1_t)$-martingale and a decreasing $(\cF^1_t)$-adapted process. By the immersion property, Lemma \ref{lem:immersion}, the martingale is also a $(\hat \cF^1_t)$-martingale, hence $(U_t)$ is also a $(\hat\cF^1_t)$-supermartingale. 

Take any $\sigma \in \cT(\hat\cF^1_t)$ such that $\sigma \ge \tau_1 \wedge \tau_2$. By the optional sampling theorem, recalling that $(U_t)$ is a uniformly integrable $(\hat\cF^1_t)$-supermartingale, we have
\begin{align*}
&\E_x \Big[ \int_{\tau_1 \wedge \tau_2}^{\sigma} e^{-rt} M(X_t) dt + e^{-r\sigma} \theta_1 \Big| \hat \cF^1_{\tau_1 \wedge \tau_2}\Big]\\
&\le
\E_x \Big[ \int_{\tau_1 \wedge \tau_2}^{\sigma} e^{-rt} M(X_t) dt + e^{-r\sigma} m(X_\sigma; \theta_1) \Big| \hat \cF^1_{\tau_1 \wedge \tau_2} \Big]\\
&=
\E_x [ U_\sigma | \hat \cF^1_{\tau_1 \wedge \tau_2} ] - \int_0^{\tau_1 \wedge\tau_2} e^{-ru} M(X_u) du\\
&\le 
U_{\tau_1 \wedge\tau_2} - \int_0^{\tau_1 \wedge\tau_2} e^{-ru} M(X_u) du
=
e^{-r (\tau_1 \wedge\tau_2)} m(X_{\tau_1 \wedge\tau_2}; \theta_1),
\end{align*}
where $\theta \le m(x; \theta)$ is used in the first inequality and the last equality is by the definition of $(U_t)$. Using Proposition \ref{prop:dm} and the classical optimal stopping theory, it is immediate to see that the above inequalities become equalities for $\sigma = \sigma^*_1$ from the statement of the lemma. These findings and the tower property justify the statement of the lemma:
\begin{align*}
J_1(x, (\tau_1, \sigma_1), (\tau_2, \sigma_2)) 
&\le  \E_{x}\bigg[\int^{\tau_{1}\wedge\tau_{2}}_{0}e^{-rt}D(X_{t})dt+1_{\tau_1\leq\tau_2}e^{-r\tau_{1}}\theta_{1}+1_{\tau_1>\tau_2}e^{-r\tau_2} m(X_{\tau_2}; \theta_1)\bigg]\\
&= J_1(x, (\tau_1, \sigma_1^*), (\tau_2, \sigma_2)).
\end{align*}
\end{proof}

\begin{corollary}\label{cor:equiv_to_m}
The optimal stopping times for \eqref{eqn:sigma1} and \eqref{eqn:sigma2} are given by 
\begin{equation*}
\sigma_i^* = \inf\{t \ge \tau^*_1 \wedge \tau^*_2: X_t \le \beta(\theta_i) \} , \quad i=1,2,
\end{equation*}
where $\beta(\theta)$ is the stopping boundary from Proposition \ref{prop:dm}(ii) corresponding to the optimal stopping problem $m(x; \theta)$.
The optimal values are 
\begin{equation}\label{eqn:cond_m}
\E_x \big[ \ind{\tau^*_i > \tau^*_{j}} e^{-r \tau^*_j} m(X_{\tau^*_j}; \theta_i)\big], \qquad i = 1,2,\ j = 3-i. 
\end{equation}
\end{corollary}

Stopping times $\sigma^*_1, \sigma^*_2$ in a Nash equilibrium are determined by Lemma \ref{lem:Nash_sigma} in the sense that for any Nash equilibrium $((\tau^*_1, \sigma^*_1), (\tau^*_2, \sigma^*_2))$ the stopping times $\sigma^*_1, \sigma^*_2$ can be replaced by those in \eqref{eqn:sigma_star} and the equilibrium condition \eqref{eqn:nash1} still holds. Hence, the only remaining decision variables for the firms are stopping times $\tau_1, \tau_2$. In view of \eqref{eqn:cond_m}, the corresponding payoffs take the form
\begin{equation}\label{eqn:payoff}
\begin{aligned}
J_1(x,\tau_{1},\tau_{2})&=\E_{x}\big[\int^{\tau_{1}\wedge\tau_{2}}_{0}e^{-rt}D(X_{t})dt+1_{\tau_1\leq\tau_2}e^{-r\tau_{1}}\theta_{1}+1_{\tau_1>\tau_2}e^{-r\tau_{2}}m(X_{\tau_{2}};\theta_1)\big],\\
J_2(x,\tau_{1},\tau_{2})&=\E_{x}\big[\int^{\tau_{1}\wedge\tau_{2}}_{0}e^{-rt}D(X_{t})dt+1_{\tau_2\leq\tau_1}e^{-r\tau_{2}}\theta_{2}+1_{\tau_2>\tau_1}e^{-r\tau_{1}}m(X_{\tau_{1}};\theta_2)\big].
\end{aligned}
\end{equation}
The Nash equilibrium condition becomes
\begin{equation}\label{eqn:nash2}
\begin{aligned}
J_1(x, \tau^*_1, \tau^*_2) &\ge J_1(x, \tau_1, \tau^*_2),\\
J_2(x, \tau^*_1, \tau^*_2) &\ge J_2(x, \tau^*_1, \tau_2),
\end{aligned}
\end{equation}
for any strategy profile $(\tau_1, \tau_2) \in \cT(\cF^1_t) \times \cT(\cF^2_t)$. By Lemma \ref{lem:equiv_to_m}, any equilibrium $(\tau^*_1, \tau^*_2)$ yields an equilibrium for the original game when $\sigma^*_1, \sigma^*_2$ are defined by \eqref{eqn:sigma_star}.

To indicate that we deal with the above reformulation of the original game, we will refer to the participants of this game as \emph{players} with the term firm reserved for the original game. The rest of the paper will consider only the above reformulation.

\begin{remark}\label{rem:kwon1}
In \cite{kwon}, the authors consider a game of the form \eqref{eqn:payoff} but with $m(X_{\tau_i}; \theta_j)$, $i \ne j$, replaced by $\hat M(X_{\tau_i})$ which is assumed to satisfy $\hat M(x) > \theta_U$ for all $x \in \cI$. When $\hat M(x)> \theta_U$, $x \in \cI$, it is easy to see that $\beta(\theta) = x_L$, and $m(x; \theta) = \hat M(x)$ for any $\theta \in [\theta_L, \theta_U]$ and $x \in \cI$, so our setup generalises \cite{kwon}. As we argue in the introduction, the condition $\hat M(\cdot) > \theta_U$ is very restrictive as it implies that it is optimal to operate in the monopoly market for every exit value $\theta$ and every state $X_t$ of the underlying factors. This contradicts the observed business reality where some lines of business are permanently abandoned by all participants, for example, due to technology advancements. From a mathematical perspective, the relaxation of this condition and allowing the monopolist to exit the market has a profound effect on the solution. The monopolist payoff $m(x; \theta)$ not only depends on the exit value of the player but it is also given in a non-explicit way and is not separated from the exit value. As a result, players' behaviour in a Nash equilibrium is structurally different leading to the game being potentially terminated simultaneously by players of different types (clustering) and exit decisions being highly sensitive to the proximity of the monopoly stopping boundary, something that cannot happen in the setting of \cite{kwon}. 
\end{remark}

We conclude this section with a sufficient condition for a Nash equilibrium using the structure of $(\cF^i_t)$-stopping times from Lemma \ref{lem:structure_Fi_times}. To this end, we define a payoff functional for a deterministic exit value: for $\sigma\in\cT(\cF_{t})$, $\theta\in[\theta_{L},\theta_{U}]$ and a random time $\tau$ on $(\Omega,\cF,\P)$ we set
\begin{equation}
J(x,\sigma,\tau;\theta)=\E_{x}\bigg[\int_{0}^{\sigma\wedge\tau}e^{-rt}D(X_{t})dt+\ind{\sigma\le\tau}e^{-r\sigma}\theta+\ind{\sigma>\tau}e^{-r\tau}m(X_{\tau}; \theta)\bigg].\label{eqn:J}
\end{equation}
This functional corresponds to the expected payoff to player $i$ of type $\theta$ whose strategy is to stop at $\sigma$ while his competitor's strategy is to stop at $\tau$. The following proposition is proved similarly to \cite[Cor.~1]{kwon} using the structure of $(\cF^i_t)$-stopping times established in Lemma \ref{lem:structure_Fi_times}.

\begin{prop}\label{prop:equil} Assume that $\hat{\tau}_{1},\hat{\tau}_{2}:\Omega\times[\theta_{L},\theta_{U}]\to[0,\infty]$
are as in Lemma \ref{lem:structure_Fi_times}. Define $\tau_{1}=\hat{\tau}_{1}(\cdot,\theta_{1})$
and $\tau_{2}=\hat{\tau}_{2}(\cdot,\theta_{2})$. If, for $F$-almost every $\theta\in[\theta_{L},\theta_{U}]$ (i.e., for the set of probability $1$ with respect to the Borel probability measure induced by the CDF $F$),
we have 
\[
\hat{\tau}_{1}(\cdot,\theta)\in\argmax_{\sigma\in\cT(\cF_{t})}J(x,\sigma,\tau_{2};\theta)
\qquad\text{and}\qquad
\hat{\tau}_{2}(\cdot,\theta)\in\argmax_{\sigma\in\cT(\cF_{t})}J(x,\sigma,\tau_{1};\theta)
\]
then $(\tau_{1},\tau_{2})$ is a Nash equilibrium. 
\end{prop} 

This proposition indicates the direction we will take in constructing the equilibrium. In Section \ref{sec:sketch}, we show that this approach naturally leads to a particular form of equilibrium strategy.

\section{Auxiliary optimal stopping problems}\label{sec:single_player}

In this section, we prove Proposition \ref{prop:dm} and derive further properties of the stopping boundaries $\alpha, \beta$ for optimal stopping problems \eqref{eqn:stop_dm}. These problems have already been used in the derivation of the game \eqref{eqn:payoff} and they will also play an important role in the derivation of a Nash equilibrium. 

We start by establishing smoothness properties of the fundamental solutions $\phi$ (decreasing) and $\psi$ (increasing) of $(\cL - r)u = 0$ and of the resolvents $\hat M, \hat D$ defined in \eqref{eqn:resolvent}.
\begin{lemma}\label{lem:smoothness}
We have $\phi, \psi, \hat M, \hat D \in C^2(\cI)$.
\end{lemma}
\begin{proof}
Due to the assumption that the boundaries of $\cI$ are not attainable, they are singular in the terminology of \cite[Sec.~II.1]{Borodin1996}. Recalling that $b (\cdot) > 0$ and $\mu, b$ are Lipschitz continuous, by \cite[Ch.~II, \S10]{Borodin1996}, the fundamental solutions $\phi$ and $\psi$ are classical solutions of the second-order ODE $(\cL -r)u = 0$, hence they are twice continuously differentiable. Using standard arguments, we can show that $\hat M$ is a classical solution to the ODE $(\cL - r ) \hat M = -M$: for any bounded interval $[a,b] \subset \cI$, the ODE $(\cL - r ) u = -M$ with the boundary conditions $u(a) = \hat M(a)$ and $u(b) = \hat M(b)$ has a unique classical solution which, using It\^o's formula, can be shown to coincide with $\hat M$. Analogous arguments apply to $\hat D$.
\end{proof}

\begin{remark}\label{rem:a_theta}
It is assumed in Assumption \ref{ass:dm} that functions $a_\theta$ and $b_\theta$ can be extended in a continuous way to $x_L, x_U$. Since the boundaries $x_L, x_U$ are singular, we have $\lim_{x \downarrow x_L} \phi(x) = \infty$ and $\lim_{x \uparrow x_U} \phi(x) \ge 0$ (\cite[Ch.~II, \S10]{Borodin1996}). By the boundedness of $\hat D$ and $\hat M$, we have $\lim_{x \downarrow x_L} a_\theta(x) = \lim_{x \downarrow x_L} b_\theta(x) = 0$.

The existence of the limit of $a_\theta$ at $x_U$ is implied by either $\lim_{x \uparrow x_U} \phi(x) > 0$ or $\lim_{x \uparrow x_U} \hat D(x) \ne \theta$.  Consider first the case $\lim_{x \uparrow x_U} \phi(x) > 0$ (entrance-not-exit boundary). The continuity of $\hat D$ established in Lemma \ref{lem:smoothness} and its monotonicity (by Assumption \ref{ass:main}(ii)) ensure that the limit $\lim_{x \uparrow x_U} {\hat D}(x)$ exists and is finite. Consequently, $\lim_{x \uparrow x_U} a_\theta(x)$ exists and is finite. However, when $\lim_{x \uparrow x_U} \phi(x) = 0$ (natural boundary) and $\hat D(x_U-) \coloneqq \lim_{x \uparrow x_U} \hat D(x) \ne \theta$, the limit $\lim_{x \uparrow x_U} a_\theta(x) = \pm \infty$ depending on the sign of $\lim_{x \uparrow x_U} (\theta - \hat D(x))$. Analogous arguments apply to $b_\theta$.
\end{remark}

\begin{lemma}\label{lem:hatDx_U}
Assume that $x_U$ is a natural boundary. Then $\hat D(x_U-) \coloneqq \lim_{x \uparrow x_U} \hat D(x) = D(x_U-) / r$ and $\hat M(x_U-) \coloneqq \lim_{x \uparrow x_U} \hat M(x) = M(x_U-) / r$.
\end{lemma}
\begin{proof}
By the monotonicity of $D$, we have
\[
\hat D(x) = \E_x \Big[ \int_0^\infty e^{-rs} D(X_s) ds \Big] \le \E_x \Big[ \int_0^\infty e^{-rs} D(x_U-) ds \Big] = D(x_U-) / r,
\]
so $\hat D(x_U-) \le D(x_U-) / r$. To show the opposite inequality, we fix $y\in\cI$, take $x>y$, and let $\tau_y=\inf\{t:X_t\le y\}$. Since $X_s\ge y$ on $[0,\tau_y)$ and $D$ is increasing, we have $D(X_s) \ge D(y)$ for $s \le \tau_y$ $\P_x$\as and so
\[
\hat D(x)\ \ge\ \E_x\Big[\int_0^{\tau_y}e^{-rs}D(X_s)\,ds\Big]\ \ge\ \frac{D(y)}{r}\big(1-\E_x[e^{-r\tau_y}]\big)=\frac{D(y)}{r}\Big(1-\frac{\phi(x)}{\phi(y)}\Big).
\]
Let $x\uparrow x_U$ (so $\phi(x)\to0$) and then $y\uparrow x_U$. This gives $\hat D(x_U-) \ge D(x_U-)/r$. The arguments for $\hat M$ are analogous.
\end{proof}
\begin{lemma}\label{lem:a_theta diff}
Let $x_L < y < x < x_U$ and $\theta \in [\theta_L, \theta_U]$. Then
\[
a_{\theta}(y) - a_{\theta}(x) = \frac1{\phi(x)} \E_x \Big[ \int_0^{\tau_y} e^{-rs} \big(D(X_s) - r\theta\big) ds \Big],
\]
where $\tau_y = \inf\{ t \ge 0\colon X_t = y\}$.
\end{lemma}
\begin{proof}
By the strong Markov property of $(X_t)$, we have
\begin{align*}
\theta - \hat D(x) 
&= \E_x \Big[ \int_0^{\tau_y} e^{-rs} \big(r\theta - D(X_s)\big) ds + e^{-r\tau_y} \big(\theta - \hat D(y) \big)\Big]\\
&= \E_x \Big[ \int_0^{\tau_y} e^{-rs} \big(r\theta - D(X_s)\big) ds\Big] + \frac{\phi(x)}{\phi(y)} \big(\theta - \hat D(y) \big).
\end{align*}
Dividing by $\phi(x)$ and rearranging yields the identity in the statement of the lemma.
\end{proof}

We turn our attention to the characterisation of the value functions $d, m$ defined in \eqref{eqn:stop_dm} and the corresponding optimal stopping times. The statement of the lemma uses the notation introduced in Assumption \ref{ass:dm}.

\begin{lemma}\label{lem:alpha}
For each $\theta \in [\theta_L, \theta_U]$, the optimal stopping time for $d(x; \theta)$ is given by
\begin{equation}\label{eqn:tau_theta}
    \tau_{\theta}^*=\inf\{t\geq 0:X_{t}\leq\alpha(\theta)\},
\end{equation}
where $\alpha$ is defined in Assumption \ref{ass:dm}.  The function $\alpha$ is continuous, increasing and satisfies $\alpha(\theta) \le c_D(\theta)$. It is, furthermore, strictly increasing for those $\theta$ for which it takes values in $\cI$. 

The value function $d(x;\theta) $ is given by
    \begin{equation}
        d(x;\theta)=
\begin{cases}
 a_{\theta}(\alpha(\theta))\phi(x)+\hat D(x),   & \mathrm{for} \: x>\alpha(\theta),\\ \theta, & \mathrm{for}\: x\leq \alpha(\theta),
\end{cases}
\label{uvalue}
    \end{equation}
$d:\cI\times [\theta_L, \theta_U]\to \R$ is jointly continuous, $x \mapsto d(x; \theta)$ is differentiable for any fixed $\theta$, and $d(x;\theta)>\theta$ for all $x>\alpha(\theta)$.
\end{lemma}
\begin{proof}
Fix $\theta \in [\theta_L, \theta_U]$. Due to the continuity, boundedness and monotonicity of $D$, the value function $d(x; \theta)$ is continuous and increasing in $x$. As the stopping payoff is constant, this implies that the stopping set is of the form $(x_L, x^*]$ for some $x^* \in \overline{\cI}$. Denote by $\hat d(x; y, \theta)$ the value corresponding to the stopping strategy given by the first hitting time of the interval $(x_L, y]$. It has an explicit form: for $y \in \overline{\cI}$,
\[        
\hat d(x;y, \theta)
=
\begin{cases}
 a_{\theta}(y)\phi(x)+\hat D(x),   & \mathrm{for} \: x>y,\\ \theta, & \mathrm{for}\: x\le y.
\end{cases}
\]
This formula is clearly true for $y \in \cI$ by \cite[Ch.~II, \S10]{Borodin1996}. When $y = x_L$, $\hat d(x; x_L, \theta) = \hat D(x)$ which agrees with the above formula because $a_\theta(x_L) = 0$ as established in Remark \ref{rem:a_theta}. The case of $y = x_U$ trivially holds as $\hat d(x; x_U, \theta) = \theta$. We also note that $a_{\theta}(x)\phi(x)+\hat D(x) = \theta$ for $x \in \cI$ which is used frequently in the remaining part of the proof. 

Due to the already asserted form of the optimal stopping set, we have
\begin{equation}\label{eqn:d_sup}
d(x; \theta) = \sup_{y \in \overline{\cI}} \hat d(x; y, \theta).
\end{equation}

After this preliminary discussion, we turn our attention to proving that $\alpha(\theta)$ is the optimal stopping boundary. We consider first the case when $\alpha(\theta) \in \cI$. Take any $x > \alpha(\theta)$. From $a_{\theta}(\alpha(\theta)) > a_\theta(x)$, we have
\begin{equation}\label{eqn:hatd}
\hat d(x; \alpha(\theta), \theta) = a_{\theta}(\alpha(\theta))\phi(x)+\hat D(x) > a_{\theta}(x)\phi(x)+\hat D(x) = \theta.
\end{equation}
Hence, recalling \eqref{eqn:d_sup}, $d(x; \theta) > \theta$, so $x$ is in the continuation set. From the arbitrariness of $x > \alpha(\theta)$, we conclude that the interval $(\alpha(\theta), x_U)$ is in the continuation set for $d(x; \theta)$. To show that $(x_L, \alpha(\theta)]$ is the stopping set, we only need to show that $\alpha(\theta)$ is in the stopping set: $d(\alpha(\theta); \theta) = \theta$.  Since $a_\theta$ attains its unique global maximum at $\alpha(\theta)$, we have 
\[
\hat d(\alpha(\theta); \alpha(\theta), \theta) =  \theta = a_{\theta}(\alpha(\theta))\phi(\alpha(\theta))+\hat D(\alpha(\theta))
> a_{\theta}(y)\phi(\alpha(\theta))+\hat D(\alpha(\theta)) = \hat d(\alpha(\theta); y, \theta)
\]
for any $y \in (x_L, \alpha(\theta))$. Hence, $d(\alpha(\theta); \theta) = \sup_{y \in \overline{\cI}} \hat d(\alpha(\theta); y, \theta) =  \hat d(\alpha(\theta); \alpha(\theta), \theta) = \theta$, so $\alpha(\theta)$ is in the stopping set. This concludes the proof that $\tau^*_\theta$ given by \eqref{eqn:tau_theta} is the optimal stopping time and the value function is given by \eqref{uvalue}.

It remains to study the case when $\alpha(\theta) \notin \cI$. Assume first that $\alpha(\theta) = x_L$, so that the candidate stopping set is empty. For any $x \in \cI$, we have $\hat D(x) = a_{\theta}(x_L)\phi(x)+\hat D(x) > a_{\theta}(x)\phi(x)+\hat D(x) = \theta$, where we used $a_\theta(x_L) = 0$ as demonstrated in Remark \ref{rem:a_theta} and the fact that $a_\theta$ attains a strict maximum at $x_L$. This yields $d(x; \theta) \ge \hat d(x; x_L, \theta) > \theta$, so $x$ is in the continuation set. The arbitrariness of $x \in \cI$ implies that $\cI$ is the continuation set and it is optimal never to stop, which is realised by the stopping time \eqref{eqn:tau_theta} with $\alpha(\theta) = x_L$ as the boundaries are not attainable.

The last case to study is $\alpha(\theta) = x_U$. 
Arguing by contradiction, assume that the stopping set is $(x_L, x^*]$ for some $x^* < x_U$. Then $\hat d(y; x^*, \theta) = d(y; \theta) \ge \hat d(y; y, \theta)$ for any $y \in (x^*, x_U)$, which implies $a_\theta(x^*) \ge a_\theta(y)$ for all $y \in (x^*, x_U)$. By the continuity of $a_\theta$ on $\overline{\cI}$ (see Assumption \ref{ass:dm}(iii)), we conclude that $a_\theta(x^*) \ge a_\theta(x_U)$. This contradicts that the unique global maximum of $a_\theta$ is attained at $x_U$. In turn, all points in $\cI$ are in the stopping set, $d(x; \theta) = \theta$ for $x \in \cI$, and the optimal stopping time is given by \eqref{eqn:tau_theta}. 

Direct verification using \eqref{uvalue} immediately shows the continuity of $d(x; \theta)$ in $x$ for a fixed $\theta$. The continuity in $\theta$ for a fixed $x$ comes from the original definition of the optimal stopping problem: for $\theta, \theta' \in [\theta_L, \theta_U]$ we obtain the bound $|d(x; \theta) - d(x; \theta')| \le |\theta - \theta'|$. The monotonicity of $d(x; \theta)$ in $\theta$ allows us to use the results of \cite{kruse1969joint} to conclude the joint continuity in $(x, \theta)$. Alternatively, we can argue directly referring to the continuity in $\theta$ which is uniform in $x$.

The continuous differentiability of $d$ in $x$ away from $x = \alpha(\theta)$ is clear as $\phi$ and $\hat D$ are continuously differentiable (see Lemma \ref{lem:smoothness}). To prove the smooth fit at $\alpha(\theta) \in \cI$, we note that the left derivative is $0$ and the right derivative is equal to
\begin{align*}
\lim_{x \downarrow \alpha(\theta)} \big[ a_{\theta}(\alpha(\theta)) \phi'(x) + \hat D'(x)\big]
&=
\lim_{x \downarrow \alpha(\theta)} \frac{\big(\theta - \hat D(\alpha(\theta))\big) \phi'(x) + \phi(\alpha(\theta)) \hat D'(x)}{\phi(\alpha(\theta))}\\
&=
\frac{\big(\theta - \hat D(\alpha(\theta))\big) \phi'(\alpha(\theta)) + \phi(\alpha(\theta)) \hat D'(\alpha(\theta))}{\phi(\alpha(\theta))},
\end{align*}
where the last equality follows from the continuity of $\phi'$ and $\hat D'$. The first order necessary condition for the global maximum of $a_\theta$ at $\alpha(\theta)$ gives
\[
0 = a'_{\theta}(\alpha(\theta)) = \frac{-\hat D'(\alpha(\theta)) \phi(\alpha(\theta)) - \big(\theta - \hat D(\alpha(\theta))\big) \phi'(\alpha(\theta))}{\phi^2(\alpha(\theta))}.
\]
Therefore, we can conclude that the right derivative of $d(x; \theta)$ at $x = \alpha(\theta)$ is $0$ and the smooth fit is proved. The continuous differentiability of $d$ in $x \in \cI$ when $\alpha(\theta) \in \{x_L, x_U\}$ is immediate.

The fact that $\alpha(\theta) \le c_D(\theta)$ comes from the identity
\[
d(x; \theta) 
= \sup_{\sigma \in \cT(\cF_t)} \E_x \Big[\int_0^\sigma e^{-rs} D(X_s) ds + e^{-r\sigma} \theta \Big]
= \theta + \sup_{\sigma \in \cT(\cF_t)} \E_x \Big[\int_0^\sigma e^{-rs} \big( D(X_s) - r\theta) ds \Big]
\]
from which it is immediate that it is suboptimal to stop whenever $D(X_s) > r\theta$. This identity also shows that the stopping set for $d$ is the same as the stopping set $\{x \in \cI: \tilde d(x; \theta) = 0 \}$ for
\begin{equation}\label{eqn:tilde d}
\tilde d(x; \theta) = \sup_{\sigma \in \cT(\cF_t)} \E_x \Big[\int_0^\sigma e^{-rs} \big( D(X_s) - r\theta) ds \Big].
\end{equation}
The function $\tilde d$ is decreasing in $\theta$. This implies that the stopping set grows with $\theta$ and therefore $\alpha$ is increasing.

To see that $\alpha$ is strictly increasing when it takes values in $\cI$, fix $y_1 < y_2$ such that $\alpha(y_1), \alpha(y_2) \in \cI$ and $\alpha(y_1) = \alpha(y_2)$. From the first order optimality condition $a_{y_i}'(\alpha(y_i)) = 0$, $i=1,2$, we have
\[
y_i \phi'(\alpha(y_i)) 
= (\hat D \phi' - \hat D' \phi)(\alpha(y_i)), \qquad i=1,2.
\]
The right-hand side does not depend on $i$, so $(y_2-y_1) \phi'(\alpha(y_i))=0$ and, consequently, we must have $\phi'(\alpha(y_i))=0$. However, recalling that $\phi$ is a fundamental solution, it satisfies $(\cL - r) \phi(\alpha(y_i)) = 0$, so $\frac{b^2(\alpha(y_i))}{2} \phi''(\alpha(y_i)) = r \phi(\alpha(y_i)) > 0$. This implies that $\phi$ has a strict local minimum at $\alpha(y_i)$, which contradicts that $\phi$ is strictly decreasing.
This combined with the monotonicity of $\alpha$ completes the proof of the strict monotonicity of $\alpha$.

We turn our attention to the continuity of $\alpha$. Assume first that $x_U$ is not a natural boundary, i.e., $\phi(x_U) > 0$ and so $|a_\theta(x_U)| < \infty$. We proceed by contradiction. Assume that there is a sequence $(\theta_{n})\subset[\theta_{L},\theta_{U}]$ converging to $\theta$ and $\alpha(\theta_{n})\to\hat{\alpha}\ne\alpha(\theta)$. The mapping $\overline\cI \times [\theta_L, \theta_U] \ni (x,\theta)\mapsto a_{\theta}(x)$ is jointly continuous due to the continuity in $x$ and in $\theta$ and the monotonicity in $\theta$ (see \cite{kruse1969joint}). Hence, $a_{\theta_{n}}(\alpha(\theta))\to a_{\theta}(\alpha(\theta))$ and $a_{\theta_{n}}(\alpha(\theta_{n}))\to a_{\theta}(\hat{\alpha})$. Since $\alpha(\theta_{n})$ is a global maximum of $a_{\theta_{n}}$, we have $a_{\theta_{n}}(\alpha(\theta_{n}))\ge a_{\theta_{n}}(\alpha(\theta))$, so $a_{\theta}(\hat{\alpha})\ge a_{\theta}(\alpha(\theta))$. This contradicts that $\alpha(\theta)$ is the unique global maximum of $a_{\theta}(\cdot)$.

When $x_U$ is a natural boundary, we have $\phi(x_U-) = 0$, so $a_\theta(x_U) = \infty$ for $\theta > \theta_0 \coloneqq D(x_U-)/r$, see Remark \ref{rem:a_theta} and Lemma \ref{lem:hatDx_U}. Hence, $\alpha(\theta) = x_U$ for $\theta > \theta_0$. Analogously, $\alpha(\theta) < x_U$ for $\theta < \theta_0$ as $a_{\theta}(x_U) = -\infty$. The continuity of $\alpha$ on $[\theta_L,\theta_0)$ is proved identically as in the above paragraph when $\hat\alpha < x_U$ as the joint continuity of $[x_L, x_U) \times [\theta_L, \theta_0) \ni (x,\theta) \mapsto a_\theta(x)$ can be asserted due to the monotonicity and coordinate-wise continuity. This is, however, not true when the function takes infinite values, so these arguments cannot be extended to $\hat\alpha = x_U$ and for asserting the continuity on $[\theta_L, \theta_U]$. We first show that $\hat\alpha$ cannot be equal to $x_U$. Take $\theta'\in (\theta,\theta_0)$. For all sufficiently large $n$, we have $a_{\theta_n}(\alpha(\theta_n))\le a_{\theta'}(\alpha(\theta_n))\to a_{\theta'}(x_U-)=-\infty$. As $\alpha(\theta_n)$ is the maximum of $a_{\theta_n}$, we have $a_{\theta_n}(\alpha(\theta_n))\ge a_{\theta_n}(\alpha(\theta))\to a_\theta(\alpha(\theta)) > -\infty$, hence a contradiction.

The continuity on $(\theta_0, \theta_U]$ is obvious, so we only need to consider $\theta = \theta_0$. We will show that $\alpha(\theta_0) = x_U$ and that $\lim_{\theta \uparrow \theta_0} \alpha(\theta) = x_U$. Notice that $D(x) - r\theta_0 = D(x) - D(x_U-) \le 0$ for any $x \in \cI$, so $x \mapsto a_{\theta_0}(x)$ is increasing by Lemma \ref{lem:a_theta diff}. Hence it attains the maximum at $x_U$, so $\alpha(\theta_0) = x_U$.  Consider now the left-hand limit of $\alpha$ at $\theta_0$: $\hat\alpha = \lim_{\theta \uparrow \theta_0} \alpha(\theta)$. Assume by contradiction that $\hat\alpha < x_U$. Choose $x_0 \in (\hat\alpha, x_U)$ such that $a_{\theta_0}(x_0) > a_{\theta_0}(\hat\alpha) + \delta$ for some $\delta > 0$; such $x_0$ exists because $x_U$ is the unique global maximum of $a_{\theta_0}(\cdot)$. By definition, $\lim_{\theta_n \uparrow \theta_0} a_{\theta_n}(x) = a_{\theta_0}(x)$ for any $x \in \cI$, so there exists $n$ such that $a_{\theta_0}(x_0) - a_{\theta_n}(x_0) \le \delta/2$. Then
\[
a_{\theta_n}(x_0) \ge a_{\theta_0}(x_0) - \delta/2 \ge a_{\theta_0}(\hat\alpha) + \delta/2 \ge a_{\theta_0}(\alpha(\theta_n)) + \delta/2 \ge a_{\theta_n}(\alpha(\theta_n)) + \delta/2,
\]
where the second inequality is by the monotonicity of $x \mapsto a_{\theta_0}(x)$ and $\alpha(\theta_n) \le \hat\alpha$, and the third one is by the monotonicity of $\theta \mapsto a_{\theta}(x)$ for any $x \in \cI$. This is a contradiction with $\alpha(\theta_n)$ being the global maximum of $a_{\theta_n}(\cdot)$. This completes the proof that $\hat\alpha = x_U$.
\end{proof}

We remark that whenever $\alpha(\theta)$ is in $\cI$, the form of the value function and the optimal stopping strategy could be obtained from \cite[Thm.~3]{Alvarez2001} under an additional assumption about the monotonicity of $x \mapsto a_\theta(x)$. This theorem, however, does not apply to the cases when $\alpha(\theta) \in \{x_L, x_U\}$ and does not cover the remaining claims in the statement of the lemma. Our proof is also more direct thanks to the monotonicity of the value function which we can assert a priori. 

An analogous statement holds for the optimal stopping problem $m(x; \theta)$.
\begin{lemma}\label{lem:beta}
For each $\theta \in [\theta_L, \theta_U]$, the optimal stopping time for $m(x; \theta)$ is given by
\begin{equation*}
   \tau_{\theta}^*=\inf\{t\geq 0:X_{t}\leq\beta(\theta)\},
\end{equation*}
where $\beta$ is defined in Assumption \ref{ass:dm}. The function $\beta$ is continuous, increasing and satisfies $\beta(\theta) \le c_M(\theta)$. It is, furthermore, strictly increasing for those $\theta$ for which it takes values in $\cI$.

The value function $m(x;\theta)$ is given by
    \begin{equation}
        m(x;\theta)=
\begin{cases}
 b_{\theta}(\beta(\theta))\phi(x)+\hat M(x), & \mathrm{for} \: x>\beta(\theta)\\ \theta, & \mathrm{for}\: x\leq \beta(\theta),
\end{cases}
\label{vvalue}
    \end{equation}
$m:\cI\times [\theta_L, \theta_U]\to \R$ is jointly continuous, $x \mapsto m(x; \theta)$ is differentiable for any fixed $\theta$, and $m(x;\theta)>\theta$ for all $x>\beta(\theta)$.
\end{lemma}

In the following corollary, we highlight the monotonicity of $d$ and $m$ in $x$ which was used in the proofs above and follows directly from the monotonicity of $D$ and $M$ and the form of the stopping functional.
\begin{corollary} \label{cor:monot_dm}
Mappings $x \mapsto d(x;\theta)$ and $x \mapsto m(x;\theta)$ are increasing for each $\theta \in [\theta_L, \theta_U]$.
\end{corollary}

\begin{proof}[Proof of Proposition \ref{prop:dm}]
All claims apart from part (iv) follow immediately from Lemmas \ref{lem:alpha} and \ref{lem:beta}. We will present arguments for $d(x; \theta)$; the proof for $m(x; \theta)$ is analogous. 

We first show that the value function $d(\cdot; \theta)$ is $C^2$ everywhere apart from $\alpha(\theta)$. Since $d(x; \theta) = \theta$ for $x \le \alpha(\theta)$, it is trivially of class $C^2$. For $x > \alpha(\theta)$ (on the continuation set), the functions $\phi$ and $\hat D$ are of class $C^2$ by Lemma \ref{lem:smoothness} and so is $d(x; \theta)$ by \eqref{uvalue}. The fact that $d(\cdot; \theta) \in C^1(\cI)$ has been established in Lemma \ref{lem:alpha}.

The inequalities \eqref{eqn:a1} are understood in the classical sense everywhere apart from one point $\alpha(\theta)$. For $x < \alpha(\theta)$, we have $(\cL -r)d(x; \theta) + D(x) = -r\theta + D(x) \le 0$ since $\alpha(\theta) \le c_D(\theta)$. For $x > \alpha(\theta)$, using the explicit form \eqref{uvalue} of $d(x; \theta)$ and that $\phi$ is a fundamental solution of $(\cL -r)\phi = 0$, we have
\[
(\cL -r)d(x; \theta) + D(x) = (\cL - r) \hat D(x) + D(x) = 0,
\]
where in the last equality we recall that $\hat D$ is the resolvent and, due to the smoothness asserted in Lemma \ref{lem:smoothness}, it is a classical solution of $(\cL -r)\hat D + D = 0$.
\end{proof}

Crucial for our equilibrium is the separation of stopping boundaries for the duopoly stopping problem $d(x; \theta)$ and the monopoly stopping problem $m(x; \theta)$ resulting from the strict inequality $M > D$.
\begin{lemma}\label{lem:alpha_beta_ineq}
We have $\alpha (\theta)>\beta(\theta)$ for $\theta \in [\theta_L, \theta_U]$ whenever $\alpha(\theta) < x_U$ and $\beta(\theta) > x_L$.
\end{lemma}
\begin{proof}
We first prove $\alpha \ge \beta$ and then show that the equality is impossible. 

Fix $\theta \in [\theta_L, \theta_U]$. Since $D \le M$ (Assumption \ref{ass:main}(ii)), we have $d(x; \theta) \le m(x; \theta)$. We also have $d(x; \theta) \ge \theta$ as $\theta$ is the stopping payoff. The stopping regions are of the form $(x_L, \alpha(\theta)] = \{x \in \cI: d(x; \theta) = \theta\}$ and $(x_L, \beta(\theta)] = \{x \in \cI: m(x; \theta) = \theta\}$. The above inequality between $d$ and $m$ implies that $\{x \in \cI: m(x; \theta) = \theta\} \subseteq \{x \in \cI: d(x; \theta) = \theta\}$, so $\beta(\theta) \le \alpha(\theta)$.

Take $\theta \in [\theta_L, \theta_U]$ such that $\alpha(\theta) < x_U$ and $\beta(\theta) > x_L$. The inequality $\alpha(\theta) > \beta(\theta)$ is easy to show whenever at least one of the boundaries is not in $\cI$. Indeed, the only possibilities that the boundaries are not in $\cI$ are $\alpha(\theta) = x_L$ and/or $\beta(\theta) = x_U$. When $\alpha(\theta) = x_L$, we must have $\beta(\theta) \le \alpha(\theta)$ but $\beta(\theta) \ne x_L$. When $\beta(\theta) = x_U$, there must be $\alpha(\theta) = x_U$ which violates the conditions of the lemma.

It remains to prove the strict inequality when $\alpha(\theta), \beta(\theta) \in \cI$. Assume by contradiction that $\alpha(\theta)=\beta(\theta)$ for some $\theta$. By Proposition \ref{prop:dm}(iv), the value functions $d(\cdot; \theta)$ and $m(\cdot; \theta)$ satisfy the following ODEs for $x \in (\alpha(\theta), x_U)$
\begin{equation*}
    \begin{aligned}
        (\cL-r)d(x;\theta)&=-D(x),\\ 
        d(\alpha(\theta); \theta)&=\theta, \\
        d'(\alpha(\theta); \theta)&=0,
    \end{aligned}
\qquad \text{and} \qquad
    \begin{aligned}
        (\cL-r)m(x; \theta)&=-M(x),\\ 
        m(\alpha(\theta); \theta)&=\theta, \\
        m'(\alpha(\theta); \theta)&=0,
    \end{aligned}
\end{equation*}
where we recall that the first derivative is zero due to the smooth fit. One should note that $m(x; \theta) > d(x; \theta)$ for $x>\alpha(\theta)$ as $M > D$. Consider a non-negative function $f(x) = m(x; \theta) - d(x; \theta)$. It satisfies
\begin{equation*}
    \begin{aligned}
        (\cL-r)f&=-(M-D), \quad x \in (\alpha(\theta), x_U)\\ 
        f(\alpha(\theta))&=0,\\
        f'(\alpha(\theta))&=0.
    \end{aligned}
\end{equation*}
By the continuity of $f$ and $f'$ and the boundary conditions at $\alpha(\theta)$, we have 
\[
\limsup_{x\downarrow\alpha(\theta)}f''(x)=\limsup_{x\downarrow\alpha(\theta)}\frac{2}{b^2(x)}\Big[-M(x)+D(x)+rf(x)-\mu(x)f'(x)\Big]<0,
\]
where we also used that $D$ and $M$ are continuous and $M > D$. Hence, there exists $\epsilon > 0$ such that for any $x\in(\alpha(\theta),\alpha(\theta)+\epsilon)$ we have $f''(x)<0$. This implies that $f'$ is strictly decreasing for $x\in(\alpha(\theta),\alpha(\theta)+\epsilon)$ and, given $f'(\alpha(\theta)) = 0$, it is negative. Likewise, $f$ is strictly decreasing in this interval and negative. This contradicts $f(x) \ge 0$. 
\end{proof}

We conclude this section with a technical result that is instrumental for the construction of our equilibrium.

\begin{lemma}\label{decreasem}
    For any $x \in \cI$, the mapping $\theta \mapsto m(x;\theta)-\theta$ is strictly decreasing on the interval $[\theta_L, \theta^*(x)) \coloneqq \{ \theta \in [\theta_L, \theta_U)\colon \beta(\theta) < x\}$ for some increasing function $\theta^*(x)$ with the convention that the interval $[\theta_L, \theta_L)$ is interpreted as the empty set.
\end{lemma}
\begin{proof}
Fix $x \in \cI$. Since $\beta(\cdot)$ is continuous and increasing, the set $\{ \theta \in [\theta_L, \theta_U)\colon \beta(\theta) < x\}$ has the form $[\theta_L, \theta^*)$ for some $\theta^* \eqqcolon \theta^*(x)$. The fact that $\theta^*(\cdot)$ is increasing follows immediately from the monotonicity of $\beta(\cdot)$.

Fix $x \in \cI$ and denote $\theta^* = \theta^*(x)$. If $\theta^* = \theta_L$, there is nothing to prove. Otherwise, let $f(\theta) = m(x; \theta) - \theta$.  Take $\theta \in [\theta_L, \theta^*)$ and $\ve > 0$ such that $\theta + \ve < \theta^*$. Using the explicit formula \eqref{vvalue} for $m(x; \theta)$ when $x > \beta(\theta)$, we have
\begin{align*}
f(\theta+\ve)-f(\theta)
&=
\phi(x)\big[b_{\theta+\ve}(\beta(\theta+\ve)) -b_{\theta}(\beta(\theta)) \big]-\ve\\
&=
\phi(x)\big[b_{\theta+\ve}(\beta(\theta+\ve)) - b_{\theta}(\beta(\theta+\ve)) + b_{\theta}(\beta(\theta+\ve)) -b_{\theta}(\beta(\theta)) \big]-\ve\\
&\le
\phi(x)\big[b_{\theta+\ve}(\beta(\theta+\ve)) - b_{\theta}(\beta(\theta+\ve)) \big]-\ve\\
&=
\phi(x)\bigg[\frac{\theta+\ve-\hat M(\beta(\theta+\ve))}{\phi(\beta(\theta+\ve))}-\frac{\theta-\hat M(\beta(\theta+\ve))}{\phi(\beta(\theta+\ve))}\bigg]-\ve\\
&=\ve\Big(\frac{\phi(x)}{\phi(\beta(\theta+\ve))}-1\Big),
\end{align*}
where the inequality is justified by $b_{\theta}(x)$ attaining its global maximum at $\beta(\theta)$, see Assumption \ref{ass:dm}. Since $\phi(x)$ is strictly decreasing and $x > \beta(\theta+\ve) \ge \beta(\theta)$, we have that $\frac{\phi(x)}{\phi(\beta(\theta+\ve))}-1<0$. Hence $f(\theta+\ve) - f(\theta) < 0$. From the arbitrariness of $\theta$ and $\ve$, we conclude that $\theta \mapsto m(x; \theta) - \theta$ is strictly decreasing on $[\theta_L, \theta^*(x))$.
\end{proof}

\begin{remark}
By continuity, the strict monotonicity from the above lemma can be extended to the closed set $[\theta_L, \theta^*(x)]$.
\end{remark}

\section{Outline of the equilibrium construction}\label{sec:sketch}

We start from a heuristic motivation for the structure of a Nash equilibrium. We recall Proposition \ref{prop:equil} which states that players' strategies are given in terms of $\hat \tau_i(\cdot, \theta)$, $i=1,2$, which are solutions to the optimal stopping problems $\sup_{\sigma \in \cT(\cF_t)} J(x, \sigma, \tau_j; \theta)$, where $\theta \in [\theta_L, \theta_U]$ is deterministic, $j = 3-i$ is the competitor's index. The functional $J$ is defined in \eqref{eqn:J} and rewritten here for the convenience of the argument:
\[
J(x,\sigma,\tau;\theta)=\E_{x}\bigg[\int_{0}^{\sigma\wedge\tau}e^{-rt}D(X_{t})dt+ e^{-r(\sigma \wedge \tau)} \Big(\ind{\sigma\le\tau}\theta+\ind{\sigma>\tau}m(X_{\tau}; \theta)\Big)\bigg].
\]
The stopping payoff in brackets is increasing in $\theta$ while the running payoff is independent of $\theta$. It is therefore natural to expect that the larger the value of $\theta$, the smaller the optimal stopping time, or, formally, $\hat \tau_i(\cdot, y_1) \ge \hat \tau_i(\cdot, y_2)$ whenever $y_1 < y_2$. This motivates the introduction of a belief process $Y^i_t = \inf \{ y \in [\theta_L, \theta_U]: \hat \tau_i(\cdot, y) < t \}$. Its left-continuous inverse would then determine $\hat \tau_i$.

Instead of formalising the above heuristic arguments, we will look for an equilibrium defined by a belief process. This approach is further motivated by a related problem in \cite{kwon} where an equilibrium of this form is constructed. Due to the \emph{symmetry} of the game (the payoffs for each player are identical), we will look for a \emph{symmetric equilibrium} in the sense that the belief processes for both players are identical. Formally, given an $(\cF_t)$-adapted decreasing process $(Y_t)$ with values in $[\theta_L, \theta_U]$, we define
\begin{equation}\label{eqn:def_hat_tau}
\hat \tau_i(\omega, \theta) = \inf \{ t \ge 0:\ Y_t(\omega) < \theta \}.
\end{equation}
In order to remove the dependence of $(Y_t)$ on the cumulative distribution function $F$ of $\theta_i$ (see Assumption \ref{ass:main}), we will work with $Z_t = F(Y_t)$ which retains the continuity and monotonicity properties of $(Y_t)$ due to the continuity and strict monotonicity of $F$. To facilitate understanding, we introduce functions
\begin{equation}\label{definitionofnotation}
Z(y)=F(y), \qquad Y(z)=F^{-1}(z),
\end{equation}
which map the process $(Y_t)$ into $(Z_t)$ and $(Z_t)$ into $(Y_t)$, respectively. Using the process $(Z_t)$, the definition \eqref{eqn:def_hat_tau} takes the form
\begin{equation}\label{eqn:def_hat_tau_Z}
\hat \tau_i(\omega, \theta) = \inf \{ t \ge 0: Z_t(\omega) < Z(\theta) \}.
\end{equation}
We will refer to $(Y_t)$ (and also, with a slight abuse of terminology, to $(Z_t)$) as a \emph{belief process} due to the discussion at the start of this section.

Informed by the equilibrium found in \cite{kwon} in a highly simplified version of the current setting, we postulate that the belief process satisfies the following relation
\begin{equation}\label{eqn:Z}
    dZ_{t}=-Z_t \lambda(X_{t},Y(Z_{t}))dt, \qquad Z_{0} = 1,
\end{equation}
with
\begin{equation}\label{eqn:lambda}
    \lambda(x,y)=\frac{ry-D(x)}{m(x;y)-y}\ind{x\le\alpha(y)},
\end{equation}
where $\alpha(y)$ is the stopping boundary corresponding to the value function $d(x; y)$, see Section \ref{sec:single_player}. Notice that the initial condition $Z_0 = 1$ corresponds to $Y_0 = Y(Z_0) = \theta_U$. By Assumption \ref{ass:dm}, since $\alpha(y) \le c_D(y)$, we have $ry - D(x) \ge 0$ for $x \le \alpha(y)$. Obviously, $m(x; y) - y \ge 0$, so $\lambda (x,y) \ge 0$; it may explode when $m(x;y)-y=0$, see the discussion in the next paragraph. This means, at least formally, that the process $(Z_t)$ is decreasing and bounded by $0$ from below. 

It is tempting to interpret the equation \eqref{eqn:Z} omega by omega as an ordinary differential equation. This is not possible since the right-hand side may explode (the denominator becomes $0$ while the numerator is strictly positive when $X_t \le \beta(Y(Z_t))$) and it may have potentially infinitely many discontinuities due to jumps whenever $X_t$ crosses the boundary $\alpha(Y(Z_t))$  (the infinite variation of $X_t$ means that there may be infinitely many such jumps in finite time). The latter problem was resolved in \cite{kwon} by constructing the maximal Carathéodory solution to an equation closely related to \eqref{eqn:Z}. This construction does not work for unbounded functions $\lambda$ as in our framework. We will, therefore, work with an appropriate approximation scheme which will allow us to make sense of \eqref{eqn:Z} as a limit. The process $(Z_t)$ that we will construct will not be a Carathéodory solution as it may have a jump. It will, however, have a key property that the pair $(X_t, Z_t)$ is a strong Markov process albeit without explicit dynamics. This will be instrumental for the application of the theory of optimal stopping of Markov processes to verify best response conditions of Proposition \ref{prop:equil}.

\begin{figure}[tb]
    \centering
    \includegraphics[width=0.7\textwidth]{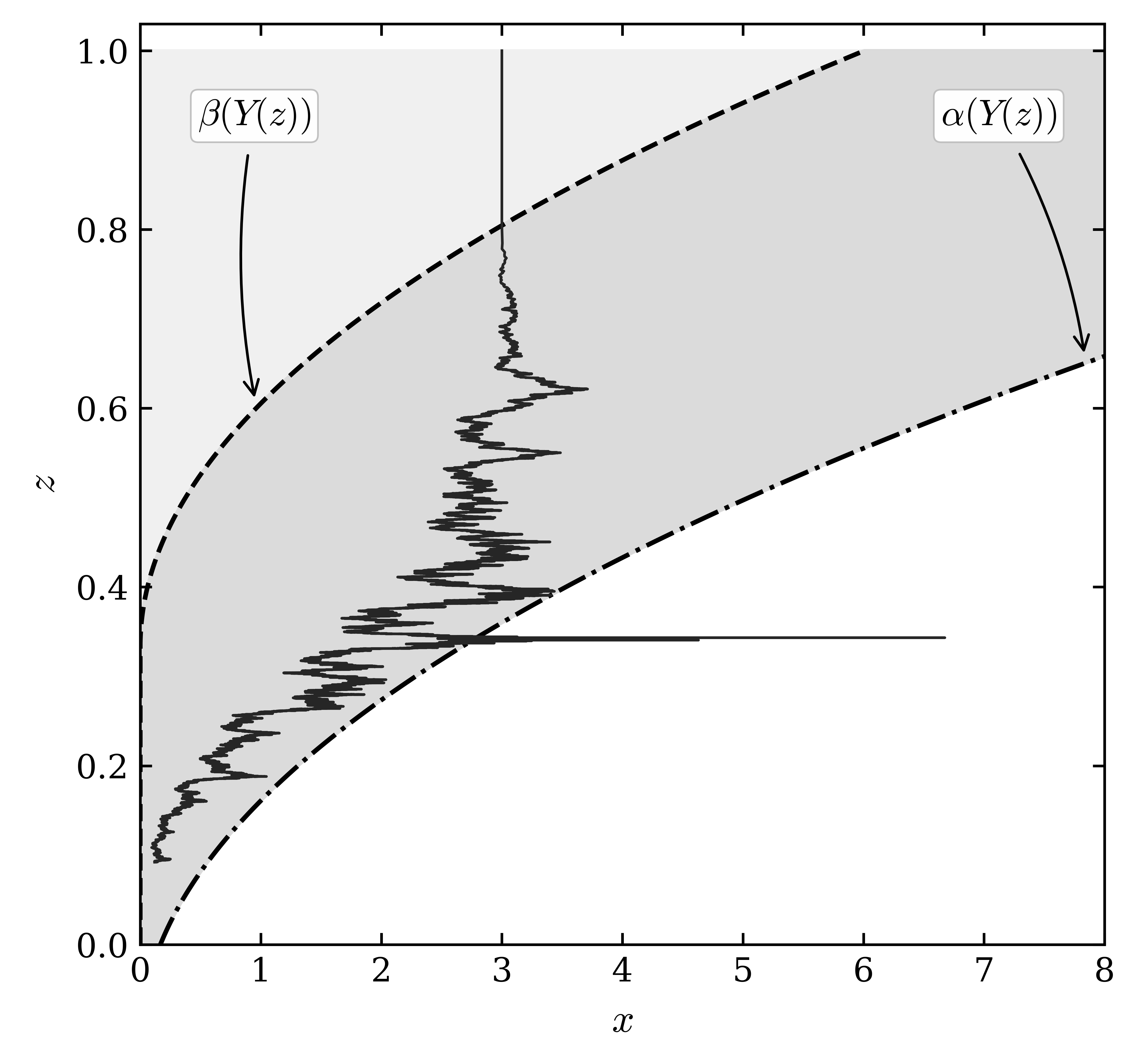}
    \caption{A plot of $(X_t,Z_t)$ (solid thick line).}
    \label{fig:example_plot}
\end{figure}

Before turning to the formal arguments, we describe the structure of the equilibrium and sketch a roadmap for the remainder of the paper. The state space, depicted in Figure \ref{fig:example_plot}, consists of points $(x, z) \in \cI \times [0,1]$. The variable $x$ corresponds to the underlying process $(X_t)$ while $z$ is the value of the belief $(Z_t)$. The solid curve shows a sample trajectory of the process $(X_t, Z_t)$. The belief process starts from $1$ and decreases to $0$. The state space is divided into 3 regions: light grey, grey and white. In the light grey region, the process is immediately shifted downwards to the dashed line -- the curve $z \mapsto \beta(Y(z))$ defining the stopping boundary of the monopoly exit problem $m(x; Y(z))$. This is because when $X_t \le \beta(Y(Z_t))$, even the monopolist cannot profitably remain in the market. Hence, all types $\theta$ of both players that are larger than $\beta^{-1}(X_t)$ (with the inverse appropriately defined if $X_t$ is not in the range of $\beta$) exit immediately driving the belief process $Z_t$ down until it reaches the dashed boundary or $0$. In the grey region, the coordinate $z$ decreases in a continuous way. This corresponds to the situation in which $X_t$ lies below the duopoly stopping boundary $\alpha(Y(Z_t))$ but still in the monopoly continuation region. Here, the players delay exit in equilibrium. When the process enters the white region, the process $(Z_t)$ remains constant resulting in horizontal segments traced out by the movement of $(X_t)$. We will show that the process $(X_t, Z_t)$ does not evolve into the light grey region after time $0$ and before reaching the bottom of the graph $z=0$ when the game finishes for all types of players.

We emphasise that, conditional on their type, each player follows a pure stopping time given by \eqref{eqn:def_hat_tau_Z}. However, because the type of the competitor is unknown, the analysis naturally adopts an averaging perspective: the strategies of the remaining types of the competitor are aggregated according to the conditional distribution of types implied by the equilibrium strategy.

The equilibrium just described is fundamentally different from that obtained in \cite{kwon}, see Remark \ref{rem:kwon1} and the introduction. In \cite{kwon}'s equilibrium, the belief process has a bounded rate of change and does not exhibit jumps. The game never finishes as the monopolist stays active forever. This corresponds to Figure \ref{fig:agreecvae} in Section \ref{sec:example} where there is no light-grey region and the grey region spans the whole space to the left of the $\alpha$-boundary. Our equilibrium allows for a richer spectrum of behaviours and requires new mathematical techniques for the construction of the belief process and for the analysis of the best response optimal stopping problems.

For the convenience of the reader we present a roadmap for the construction of our equilibrium:
\begin{enumerate}
\item In Section \ref{sec:belief}, we construct a process $(Z_t)$ that formally solves \eqref{eqn:Z}. It is defined as a limit of an appropriate approximating sequence $(Z^{\delta_n}_t)$, see \eqref{eqn:Z_def}. Each $(Z^{\delta_n}_t)$ is a unique Carathéodory solution to an ODE (i.e., is a solution in an integral form). Theorem \ref{thm:markov_Z} proves the strong Markov property of $(X_t, Z_t)$ and Theorem \ref{thm:prop_Z} shows the properties of $(Z_t)$ such as the continuity of trajectories after time $0$ and that $(X_t, Z_t)$ does not enter the light grey region illustrated in Figure \ref{fig:example_plot}. The process $(Z_t)$ defines the equilibrium via \eqref{eqn:def_hat_tau_Z}.
\item The best response optimal stopping problem from Proposition \ref{prop:equil} is studied in Section \ref{sec:value}. It is shown to have a Markovian structure on the state space consisting of $(X_t, Z_t)$ with the value function defined in \eqref{eqn:v_def}. Although it is in a non-standard form, we adapt results from \cite{shiryaev} to show that the hitting time of the stopping set is optimal (Theorem \ref{optimalstop}). 
\item The stopping set is derived in Section \ref{sec:proof-thm}. We divide the state space into disjoint regions, and, using different methods, show whether they belong to the stopping or the continuation sets. The inclusion in the stopping region is achieved by showing that the value function is dominated either by an appropriately scaled $m(x; \theta)$ (Lemma \ref{optimaliszero}) or by an appropriately scaled $d(x; \theta)$ (Lemma \ref{alphainclude}). The proof of the latter is particularly interesting, exploiting the characterisation of the value function as the smallest majorant of the payoff satisfying an appropriate supermartingale condition. It is then shown that $(x,z)\mapsto zd(x; \theta)$ majorises the payoff and satisfies the same supermartingale condition, so it must also dominate the value function. The continuation region is identified by constructing stopping times which lead to a strictly better payoff than immediate stopping. These results are also deeply technical and use the separation of boundaries $\alpha(\theta) > \beta(\theta)$ when at least one of them takes values in $\cI$. The stopping set in its entirety is shown in Theorem \ref{thm:stopping-set}. 
\item The stopping set is closed and stated for the two-dimensional process $(X_t, Z_t)$. However, the postulated strategy \eqref{eqn:def_hat_tau_Z} is a hitting time of an open set by the process $(Z_t)$ only. The equivalence between those two stopping times is established in Lemma \ref{lem:stop-Z-transform}.
\end{enumerate}

\section{Dynamics of the belief process}\label{sec:belief}

In this section, we will construct a sequence of processes $(Z^\delta_t)$ which, in the limit, defines a solution of \eqref{eqn:Z}. We will further study the properties of the pair $(X_t, Z_t)$, in particular, their continuity and strong Markov property. In the construction of $(Z_t)$ we will provide frequent references to \cite{kwon} to indicate similarities, often at the level of ideas rather than mathematical details. We emphasise that the present setting is substantially more difficult.

\begin{remark}
Define formally $A_t = -\log(Z_t)$ with the convention that $\log(0) = -\infty$. Again, formally, the process $(A_t)$ satisfies
\[
dA_{t}= \lambda(X_{t},\hat Y(A_{t}))dt, \qquad A_0 = -\log(Z_0),
\]
where $\hat Y(a) = Y(e^{-a})$. Notice that $Z_0 = 1$ corresponds to $A_0 = 0$ and the process $(A_t)$ is increasing as $\lambda \ge 0$. Such a process is studied in \cite{kwon}, where $\lambda$ is bounded. Here, however, $\lambda$ may explode so it is more convenient to work with the process $(Z_t)$ taking values in the interval $[0,1]$ for which the explosion of $\lambda$ may drive the process to $0$ while $(A_t)$ would explode to $\infty$.
\end{remark}

\subsection{Construction of $(X_t, Z_t)$}

For $\delta>0$, consider a process $(Z^\delta_t)$ that, formally, satisfies
\begin{equation}\label{eqn:Z_delta}
    dZ^\delta_{t}=-Z^\delta_{t}\lambda^\delta(X_{t},Y(Z^\delta_{t}))dt,
\end{equation}
with
\begin{equation*}
    \lambda^\delta(x,y)=\frac{ry-D(x)}{(m(x;y)-y)\vee\delta}\ind{x\leq\alpha(y)}.
\end{equation*}
The function $\lambda^\delta$ is bounded since $y \in [\theta_L, \theta_U]$ and $D$ is a bounded function. We will first make sense of the solution of \eqref{eqn:Z_delta} following arguments of \cite{kwon} with additional complications arising from the denominator of $\lambda^\delta$ depending on $y$ in a more complex way than in the aforementioned paper. In particular, the following lemma requires a non-trivial property of $m$ which is stated in Lemma \ref{decreasem}. 

\begin{lemma}\label{lambdadeltamontone}
    The mapping $(x,y)\mapsto \lambda^\delta(x,y)$ is decreasing in $x$ and increasing in $y$.
\end{lemma}
\begin{proof}
Fix $x \in \cI$. Lemma \ref{decreasem} shows that $y \mapsto m(x;y)-y$ is strictly decreasing in $y$ whenever $\beta(y) < x$ and equal to $0$ when $\beta(y) \ge x$. Hence, $y \mapsto (m(x;y)-y)\vee \delta$ is decreasing and also bounded from below by $\delta$. The numerator of $\lambda^\delta$ is strictly increasing in $y$ and non-negative when $x \le \alpha(y)$. Furthermore, $y \mapsto \ind{x \le \alpha(y)}$ is increasing, see Proposition \ref{prop:dm}(iii). Hence $\lambda^\delta(x, y)$ is increasing in $y$.

Fix $y \in [\theta_L, \theta_U]$. The mapping $x\mapsto ry-D(x)$ is decreasing since the function $D$ is increasing. Clearly, $x \mapsto \ind{x \le \alpha(y)}$ is decreasing. The mapping $x\mapsto m(x; y)-y$ is increasing, see Corollary \ref{cor:monot_dm}.  Combining these observations with the fact that the numerator and denominator are non-negative when $x \le \alpha(y)$ shows that $\lambda^\delta(x, y)$ is decreasing in $x$.
\end{proof}

\begin{lemma}\label{lem:odesolutionconstruction}
For each $(x, z) \in \cI \times [0,1]$, denoting by $(X^x_t)$ the diffusion $(X_t)$ with $X_0 = x$, there is a unique continuous process $(Z^{x,z, \delta}_t)$ that satisfies
\begin{equation}\label{eqn:A_delta}
    Z^{x,z,\delta}_{t}=z-\int^{t}_{0}Z^{x,z,\delta}_s\, \lambda^\delta(X^x_{s},Y(Z^{x,z,\delta}_{s}))ds.
\end{equation}
\end{lemma}
\begin{proof}
The proof borrows ideas from the proof of \cite[Proposition 4.5]{kwon} but there are significant additional complications that merit its complete presentation. The uniqueness claim is new. In the proof we suppress superscripts $x,z$ in the notation for $(X_t)$ and $(Z^\delta_t)$.

Let 
\begin{equation}\label{eqn:ldelta}
    l^\delta(x,y)=\frac{ry-D(x)}{(m(x;y)-y)\vee\delta},
\end{equation}
so that $\lambda^\delta(x,y) = l^{\delta}(x,y) \ind{x \le \alpha(y)}$. We want to rewrite the indicator function in terms of an explicit condition on $y$ given $x$. Recalling from Lemma \ref{lem:alpha} that $\alpha(y)$ is strictly increasing for those $y$ that satisfy $\alpha(y) \in \cI$, its inverse is well defined for $x \in (\alpha(\theta_L), \alpha(\theta_U))$. This allows us to define
\[
\hat\alpha^{-1}(x) = 
\begin{cases}
\theta_L, & x \le \alpha(\theta_L),\\
\alpha^{-1}(x), & x \in (\alpha(\theta_L), \alpha(\theta_U)),\\
\theta_U + (x - \alpha(\theta_U)), & x \ge \alpha(\theta_U).
\end{cases}
\]
The odd value in the above definition for $x \ge \alpha(\theta_U)$ is required to obtain the equality
\[
\{ (x, y) \in \cI \times [\theta_L, \theta_U]: x \le \alpha(y) \} = \{ (x, y) \in \cI \times [\theta_L, \theta_U]: y \ge \hat\alpha^{-1}(x)\}.
\]
This restatement of the condition between $x$ and $y$ is used in the following definition of a continuous approximation of $\lambda^\delta$:
\begin{equation}
\begin{aligned}\label{eqn:lambda_e_d}
&\lambda^{\epsilon,\delta}(x,y)=l^\delta(x,y)\ind{y\ge \hat\alpha^{-1}(x)}+l^\delta(x,\hat\alpha^{-1}(x) \wedge \theta_U) \frac{1}{\epsilon}(y-\hat\alpha^{-1}(x)+\epsilon)^+ \ind{y<\hat\alpha^{-1}(x)}.
\end{aligned}
\end{equation}
Note that the second term of $\lambda^{\epsilon, \delta}(x,y)$ converges to $0$ with $\epsilon \to 0$ for $x > \alpha(\theta_U)$ because $(y-\hat\alpha^{-1}(x)+\epsilon)^+ \le \big(\epsilon - (x - \alpha(\theta_U))\big)^+$ when $x > \alpha(\theta_U)$. Similarly, it converges to $0$ for $x > \alpha(y)$, $y < \theta_U$. Hence, $\lambda^{\epsilon, \delta}$ converges pointwise to $\lambda^\delta$ as $\epsilon \to 0$. 

Take $y_2 > y_1 \ge \hat\alpha^{-1}(x)$ and notice that
\[
\begin{aligned}
|l^\delta(x, y_2) - l^\delta(x, y_1)|
&\le \left|\frac{ry_2-D(x)}{(m(x;y_2)-y_2)\vee\delta} - \frac{ry_1-D(x)}{(m(x;y_2)-y_2)\vee\delta}\right|\\
&\hspace{12pt}+
|ry_1 - D(x)| \left|\frac{1}{(m(x;y_2)-y_2)\vee\delta} - \frac{1}{(m(x;y_1)-y_1)\vee\delta}\right|\\
&\le 
\frac{r}{\delta} |y_2-y_1| + C \left| \frac{ (m(x;y_1)-y_1) \vee \delta - (m(x;y_2)-y_2)\vee\delta }{\big((m(x;y_2)-y_2)\vee\delta\big) \big((m(x;y_1)-y_1)\vee\delta\big)}\right|\\
&\le
\frac{r}{\delta} |y_2-y_1| + C \frac{ |m(x; y_1) - m(x; y_2)| + |y_2 - y_1|}{\delta^2}\\
&\le
\frac{r}{\delta} |y_2-y_1| + C \frac{2}{\delta^2} |y_2-y_1|
\end{aligned}
\]
where $C = r(|\theta_U|\vee|\theta_L|) + \sup_{x' \in \cI} |D(x')|$ and in the last inequality we used that $|m(x; y_2) - m(x; y_1)| \le |y_2 - y_1|$. The latter bound is obtained from the definition of the stopping problem $m(x; y)$ as follows. Take $y_2 > y_1$ and let $\tau^*$ be the optimal stopping time for $m(x; y_2)$ which exists by Lemma \ref{lem:beta}. Then
\[
0 \le m(x; y_2) - m(x; y_1) \le \E_x \big[e^{-r\tau^*} (y_2 - y_1) \big] \le y_2 - y_1,
\]
where the first inequality is because $y \mapsto m(x; y)$ is increasing and $y_2 > y_1$, and the second inequality arises from a suboptimal payoff for $m(x; y_1)$ and the optimal payoff for $m(x; y_2)$. 

Having demonstrated that $l^\delta(x, y)$ is Lipschitz continuous in $y \ge \hat\alpha^{-1}(x)$ with the constant independent of $x$, it suffices to notice that $\hat \alpha^{-1}$ is continuous in order to see that $\lambda^{\epsilon, \delta}$ is continuous as a function of $(x,y)$ and Lipschitz continuous in $y$ uniformly in $x$; indeed, the Lipschitz continuity of $y \mapsto \lambda^{\epsilon, \delta}$ for $y < \hat\alpha^{-1}(x)$ with the Lipschitz coefficient $C/(\delta\epsilon)$ is immediate from the definition.

Consider an ODE for each fixed $\omega$ (suppressed in the notation for the sake of conciseness)
\begin{equation}\label{epsilondeltaode}
    dZ^{\epsilon,\delta}_{t}=-Z^{\epsilon,\delta}_{t}\lambda^{\epsilon,\delta}(X_{t},Y(Z^{\epsilon,\delta}_{t}))dt, \qquad Z^{\epsilon,\delta}_0 = z.
\end{equation}
Since $Y(\cdot)$ is Lipschitz continuous by Assumption \ref{ass:main}, the right-hand side of the above ODE is Lipschitz continuous in $z$ for any trajectory $t \mapsto X_t(\omega)$. Therefore, there exists a unique solution. Notice that the solution $(Z^{\epsilon, \delta}_t)$ is decreasing in $\epsilon$ (by the comparison principle for ODEs, \cite[Theorem III-5-1]{ode}) and it is bounded from above by $1$ (because the right-hand side is non-negative and $Z_0^{\epsilon,\delta} \in [0,1]$). Hence $Z^{0,\delta}_{t}\coloneqq\lim_{\epsilon\downarrow 0}Z^{\epsilon,\delta}_{t}$ exists and is decreasing in $t$. It is also lower semicontinuous (as an increasing limit of continuous functions), which together with the monotonicity implies that it is right-continuous (and adapted), which will be used to apply results from stochastic analysis later in the proof.

Consider $x \in \cI$ and a sequence $[\theta_L, \theta_U] \ni y^\epsilon \uparrow y \ne \hat\alpha^{-1}(x)$. Then $\lambda^{\epsilon, \delta}(x, y^\epsilon)$ converges to $\lambda^{\delta}(x,y)$. This convergence cannot be asserted when $y = \hat\alpha^{-1}(x)$ as the limit depends on the speed of convergence of $y^\epsilon$ to $y$. Since $\lambda^{\epsilon,\delta}(x,y)\ge 0$, we can apply Fatou's lemma to obtain
\begin{multline*}
\liminf_{\epsilon\downarrow 0}\int^{t}_{0} Z^{\epsilon,\delta}_{s} \lambda^{\epsilon,\delta}(X_{s},Y(Z^{\epsilon,\delta}_{s})) \ind{\hat \alpha^{-1}(X_{s}) < Y(Z^{0,\delta}_{s})}ds\\
\geq\int^{t}_{0}Z^{0,\delta}_{s} \lambda^\delta(X_{s},Y(Z^{0, \delta}_{s})) \ind{\hat \alpha^{-1}(X_{s}) < Y(Z^{0,\delta}_{s})} ds.
\end{multline*}
For the upper bound, we use the inequality
 \begin{equation*}
     \lambda^{\epsilon,\delta}(x,y)\leq l^\delta(x,y)\ind{\hat\alpha^{-1}(x)<y}+l^\delta(x,\hat\alpha^{-1}(x) \wedge \theta_U)\ind{\hat\alpha^{-1}(x)-\epsilon<y\leq\hat\alpha^{-1}(x)}
\end{equation*}
and the observation that 
\[
\limsup_{\epsilon\downarrow 0}\, \{\hat \alpha^{-1}(X_{s})-\epsilon<Y(Z^{\epsilon,\delta}_{s})\leq \hat\alpha^{-1}(X_{s})\} 
\subset 
\{Y(Z^{0, \delta}_{s})=\hat\alpha^{-1}(X_{s})\},
\]
which follows from the fact that the mapping $\epsilon \mapsto Y(Z^{\epsilon, \delta}_s)$ is decreasing. Reverse Fatou's lemma yields
\begin{align*}
&\limsup_{\epsilon\downarrow 0}\int^{t}_{0}Z^{\epsilon,\delta}_{s} \lambda^{\epsilon,\delta}(X_{s},Y(Z^{\epsilon,\delta}_{s}))ds\\
&\leq
\int^{t}_{0} Z^{0,\delta}_{s} \left(l^\delta\big(X_{s},Y(Z^{0,\delta}_{s})\big)\ind{\hat\alpha^{-1}(X_{s})<Y(Z^{0,\delta}_{s})}+l^\delta\big(X_{s},\hat\alpha^{-1}(X_{s}) \wedge \theta_U\big)\ind{Y(Z^{0,\delta}_{s}) = \hat \alpha^{-1}(X_{s})} \right) ds.
\end{align*}
We note that the map $s \mapsto Z^{0, \delta}_{s}$ is decreasing and hence of finite variation. It is also right-continuous, so we can apply \cite[Lemma~4.3]{kwon} to conclude that 
\begin{equation}\label{eqn:Z equal alpha}
\int^{t}_{0} Z^{0,\delta}_{s}\ind{Y(Z^{0,\delta}_{s}) = \hat \alpha^{-1}(X_{s})}  ds = 0, \qquad \P\as
\end{equation}

Combining the lower and upper bounds derived above with the identity \eqref{eqn:Z equal alpha}, we conclude that  the following limit exists $\P$\as with the set of measure zero independent of $t$:
\begin{equation*}
\lim_{\epsilon\downarrow 0}\int^{t}_{0}Z^{\epsilon, \delta}_{s} \lambda^{\epsilon, \delta}(X_{s},Y(Z^{\epsilon, \delta}_{s}))ds=\int^{t}_{0}Z^{0, \delta}_{s}\lambda^{\delta}(X_{s},Y(Z^{0, \delta}_{s}))ds.
\end{equation*}
Hence $Z^\delta := Z^{0, \delta}$ satisfies \eqref{eqn:A_delta}, from which we further conclude that it is continuous.

It remains to prove uniqueness.\footnote{The authors are grateful to Alexander Veretennikov for discussions which led to this claim and its proof.}
Take $z_1, z_2 \in (0,1]$, $x \in \cI$ and recall that $z \mapsto \lambda^\delta(x, Y(z))$ is increasing (Lemma \ref{lambdadeltamontone} and the fact that $z \mapsto Y(z)$ is increasing). Using this monotonicity and $\lambda^\delta \ge 0$, we will show that
\begin{equation}\label{eqn:z1}
(z_1-z_2)\big(z_1 \lambda^\delta(x, Y(z_1)) - z_2 \lambda^\delta(x, Y(z_2)) \big) \ge 0.
\end{equation}
When $z_1 = z_2$ the inequality is trivially satisfied. Assume $z_1 > z_2$. The first factor is positive and we study the second factor:
\begin{equation}\label{eqn:zlambda}
z_1 \lambda^\delta(x, Y(z_1)) - z_2 \lambda^\delta(x, Y(z_2)) 
\ge z_2 \big(\lambda^\delta(x, Y(z_1)) - \lambda^\delta(x, Y(z_2))\big) \ge 0,
\end{equation}
because the expression in the brackets on the right-hand side is non-negative by the monotonicity of $z \mapsto \lambda^\delta(x, Y(z))$. The proof of \eqref{eqn:z1} when $z_1 < z_2$ is analogous.  

Let $(\tilde Z^\delta_t)$ be another continuous process satisfying \eqref{eqn:A_delta} with the initial condition $z$. By the chain rule, we obtain
\begin{align*}
&(Z^{\delta}_t - \tilde Z^\delta_t)^2=\int^t_0 -2(Z^{\delta}_s - \tilde Z^\delta_s) \Big(Z^{\delta}_s \lambda^{\delta}(X_s,Y(Z^{\delta}_s)) - \tilde Z^\delta_s \lambda^\delta(X_s,Y(\tilde Z^\delta_s))\Big)ds \le 0,
\end{align*}
with the inequality following from \eqref{eqn:z1}. As the left-hand side is non-negative, this shows that $Z^\delta_t = \tilde Z^\delta_t$.
\end{proof}

\begin{remark}
Processes that satisfy \eqref{eqn:A_delta} are called Carathéodory solutions. Carathéodory solutions of ODEs lacking classical solutions may not be unique. However, we prove the uniqueness of our solution exploiting the monotonicity properties of $\lambda^\delta$. An analogous proof shows that the Carathéodory solution of \cite[Eq.~(4.2)]{kwon} is unique; this improves on the results of that paper, where it is only shown that it is a maximal solution.
\end{remark}

Having constructed the process $(Z^{x,z,\delta}_t)$ we explore its dependence on the initial condition $(x,z)$. These results are needed to show the continuous dependence of $(X^x_t, Z^{x,z,\delta}_t)$ on $(x, z)$ and that this process satisfies a strong Markov property. 

\begin{lemma}\label{lem:comparison_Z}
The process $(Z^{x,z,\delta}_t)$ is increasing in $x,z$ in the sense that $\P (Z^{x, z, \delta}_t \ge Z^{x', z', \delta}_t \ \forall t \ge 0) = 1$ for $x \ge x'$ and $z \ge z'$.
\end{lemma}
\begin{proof}
Take $z > z'$ and $t > 0$. By Tanaka's formula \cite[Ch.~VI, Thm.~1.2]{yor} and the inequality \eqref{eqn:zlambda}, we have 
\begin{align*}
&(Z^{x, z', \delta}_{t}-Z^{x, z, \delta}_{t})^+\\
&=
(Z^{x, z', \delta}_{0}-Z^{x, z, \delta}_{0})^+\\
&\quad- \int^t_0 \ind{Z^{x, z', \delta}_s > Z^{x, z, \delta}_{s}}\Big(Z^{x, z', \delta}_s\lambda^\delta(X^x_{s},Y(Z^{x, z', \delta}_s))-Z^{x, z, \delta}_s\lambda^\delta(X^x_{s},Y(Z^{x, z, \delta}_{s}))\Big)ds\\
&\le (Z^{x, z', \delta}_{0}-Z^{x, z, \delta}_{0})^+ = (z' - z)^+ = 0,
\end{align*}
so $Z^{x, z', \delta}_{t} \le Z^{x, z, \delta}_{t}$. 

A slightly longer argument is needed to prove the monotonicity in $x$. Let $x > x'$ and apply again Tanaka's formula:
\begin{equation}\label{eqn:tza}
\begin{aligned}
&(Z^{x', z, \delta}_{t}-Z^{x, z, \delta}_{t})^+\\
&=
0 - \int^t_0 \ind{Z^{x', z, \delta}_s > Z^{x, z, \delta}_{s}}\Big(Z^{x', z, \delta}_s\lambda^\delta(X^{x'}_{s},Y(Z^{x', z, \delta}_s))-Z^{x, z, \delta}_s\lambda^\delta(X^{x}_{s},Y(Z^{x, z, \delta}_{s}))\Big)ds.
\end{aligned}
\end{equation}
For the expression under the integral, we have
\begin{equation}\label{eqn:ZZZ}
\begin{aligned}
&Z^{x', z, \delta}_s\lambda^\delta(X^{x'}_{s},Y(Z^{x', z, \delta}_s))-Z^{x, z, \delta}_s\lambda^\delta(X^x_{s},Y(Z^{x, z, \delta}_{s}))\\
&=
Z^{x', z, \delta}_s \Big(\lambda^\delta(X^{x'}_{s},Y(Z^{x', z, \delta}_s))- \lambda^\delta(X^{x}_{s},Y(Z^{x', z, \delta}_{s}))\Big)\\
&\quad+
Z^{x', z, \delta}_s\lambda^\delta(X^{x}_{s},Y(Z^{x', z, \delta}_s))-Z^{x, z, \delta}_s\lambda^\delta(X^x_{s},Y(Z^{x, z, \delta}_{s})).
\end{aligned}
\end{equation}
The comparison principle for SDEs \cite[Ch.~IX, Thm.~3.7]{yor} implies that $X^{x}_s \ge X^{x'}_s$ for all $s \ge 0$. This and Lemma \ref{lambdadeltamontone} show $\lambda^\delta(X^x_s,y)\le\lambda^\delta(X^{x'}_s,y)$ for all $y$ so that the first difference in \eqref{eqn:ZZZ} is non-negative. The second difference is non-negative by \eqref{eqn:zlambda} on the event $Z^{x', z, \delta}_s > Z^{x, z, \delta}_{s}$. Hence, we conclude from \eqref{eqn:tza} that $(Z^{x', z, \delta}_{t}-Z^{x, z, \delta}_{t})^+ \le 0$, so $Z^{x', z, \delta}_{t} \le Z^{x, z, \delta}_{t}$.

The above arguments show that $Z^{x', z', \delta}_{t} \le Z^{x, z, \delta}_{t}$ for all rational $t \ge 0$ $\P$\as whenever $z \ge z'$ and $x \ge x'$. It is then extended by the continuity to all $t \ge 0$.
\end{proof}

\begin{lemma}\label{lem:z_lipsch}
For any $x \in \cI$ and $\delta > 0$, we have $|Z^{x, z, \delta}_t - Z^{x, z', \delta}_t| \le |z-z'|$, $z, z' \in [0,1]$.
\end{lemma}
\begin{proof}
Take $z > z'$ and fix $x, t$. Then $Z^{x, z, \delta}_s \ge Z^{x, z', \delta}_s$ for all $s \ge 0$ by Lemma \ref{lem:comparison_Z}. Using \eqref{eqn:A_delta}, we write
\begin{align*}
Z^{x, z, \delta}_t - Z^{x, z', \delta}_t 
&= Z^{x, z, \delta}_0 - Z^{x, z', \delta}_0 
- \int^{t}_{0} \big[ Z^{x,z,\delta}_s \lambda^\delta(X^x_{s},Y(Z^{x, z, \delta}_{s})) - Z^{x, z', \delta}_s \lambda^\delta(X^x_{s},Y(Z^{x, z', \delta}_{s})) \big] ds\\
&\le
z - z',
\end{align*}
where the inequality follows from \eqref{eqn:zlambda}. This shows that $z \mapsto Z^{x,z,\delta}_t$ is Lipschitz continuous with the constant independent of $x$ and $\omega$.
\end{proof}

\begin{lemma}\label{continuitywrtstart}
For any sequence $(x_n, z_n)$ converging to $(x,z)$, we have $Z^{x_n,z_n,\delta}_t \to Z^{x,z,\delta}_t$ $\P$\as for any $t \ge 0$. 
\end{lemma}
\begin{proof}
By Lemma \ref{lem:z_lipsch}, $z \mapsto Z^{x,z,\delta}_t$ is Lipschitz continuous with the constant independent of $x$ and $\omega$.

To prove the continuity in $x$, we first take $x_n\uparrow x_0$ and fix $z$ and $t > 0$. We recall a standard fact that under our assumptions on the coefficients of $(X_t)$, for any $T > 0$, we have $X^{x_n}_s \to X^{x_0}_s$ for all $s \in [0, T]$ $\P$\as We will write $Z^{x_n}$ for $Z^{x_n, z, \delta}$. By the monotonicity of $x \mapsto Z^x$ established in Lemma \ref{lem:comparison_Z}, there exists a limit $\Bar{Z}_t = \lim_{n \to \infty} Z^{x_n}_t$ and $\Bar{Z}_t \le Z^{x_0}_t$. Recall that $Z^{x_n}$ satisfies the integral equation \eqref{eqn:A_delta}. This allows us to obtain the lower bound
\begin{align*}
Z^{x_n}_t &= z - \int_0^t Z^{x_n}_s \lambda^\delta(X^{x_n}_s, Y(Z^{x_n}_s)) ds
\ge z - \int_0^t Z^{x_n}_s \lambda^\delta(X^{x_n}_s, Y(Z^{x_0}_s)) ds\\
&\ge z - \int_0^t Z^{x_0}_s \lambda^\delta(X^{x_n}_s, Y(Z^{x_0}_s)) ds,
\end{align*}
where we used that $\lambda^\delta$ is increasing in $y$ for the first inequality. Reverse Fatou's lemma allows us to take the limit $n\to\infty$:
\begin{align*}
\Bar{Z}_t &\ge z - \limsup_{n \to \infty}\int_0^t Z^{x_0}_s \lambda^\delta(X^{x_n}_s, Y(Z^{x_0}_s)) ds
\ge z - \int_0^t Z^{x_0}_s \limsup_{n \to \infty} \lambda^\delta(X^{x_n}_s, Y(Z^{x_0}_s)) ds\\
&= z - \int_0^t Z^{x_0}_s \lambda^\delta(X^{x_0}_s, Y(Z^{x_0}_s)) ds = Z^{x_0}_t,
\end{align*}
with the penultimate equality justified by the left-continuity of $x \mapsto \lambda^\delta(x, y)$ for any $y$ and $X^{x_n}_s \uparrow X^{x_0}_s$ for all $s \in [0, T]$ $\P$\as Recalling that $\Bar{Z}_t \le Z^{x_0}_t$ proves that $Z^{x_n}_t \to Z^{x_0}_t$ when $n\to\infty$.

Consider now $x_n \downarrow x_0$. Then $\Bar{Z}_t \ge Z^{x_0}_t$. To obtain the opposite inequality we follow an analogous argument as above. We notice that
\[
Z^{x_n}_t \le z - \int_0^t Z^{x_0}_s \lambda^\delta(X^{x_n}_s, Y(Z^{x_0}_s)) ds.
\]
We now apply Fatou's lemma but it requires more attention as $\lambda^\delta$ is not right-continuous in $x$. Recalling the notation \eqref{eqn:ldelta}, we have
\begin{align*}
\Bar{Z}_t 
&\le z - \int_{0}^t Z^{x_0}_s \liminf_{n\to\infty} \lambda^\delta(X^{x_n}_s, Y(Z^{x_0}_s)) ds
\le z - \int_{0}^t Z^{x_0}_s l^\delta(X^{x_0}_s, Y(Z^{x_0}_s))\ind{X^{x_0}_s < \alpha(Y(Z^{x_0}_s))} ds\\
&= z - \int_{0}^t Z^{x_0}_s \lambda^\delta(X^{x_0}_s, Y(Z^{x_0}_s))ds  + \int_{0}^t Z^{x_0}_s l^\delta(X^{x_0}_s, Y(Z^{x_0}_s))\ind{X^{x_0}_s = \alpha(Y(Z^{x_0}_s))} ds = Z^{x_0}_t,
\end{align*}
where in the second inequality we used $\liminf_{n\to\infty}\ind{X^{x_n}_s \le \alpha(Y(Z^{x_0}_s))} \ge \ind{X^{x_0}_s < \alpha(Y(Z^{x_0}_s))}$ and the continuous dependence of $X^x_t$ on $x$: $X^{x_n}_s \downarrow X^{x_0}_s$ for all $s \in [0, T]$ $\P$\as; the last equality is a consequence of \cite[Lemma~4.3]{kwon} which shows that $\int_{0}^t \ind{X^{x_0}_s = \alpha(Y(Z^{x_0}_s))} ds = 0$ $\P$\as

The above assertions over monotone sequences $(x_n)$ imply that $Z^{x_n}_t \to Z^x_t$ $\P$\as for an arbitrary sequence $x_n \to x$. The convergence of $Z^{x_n, z_n, \delta}_t$ to $Z^{x, z, \delta}_t$ follows from the continuity in $z$ which is uniform in $x$ and the convergence over $x_n$ for a fixed $z$:
\[
\lim_{n \to \infty} Z^{x_n, z_n, \delta}_t = \lim_{n \to \infty} \big(Z^{x_n, z_n, \delta}_t - Z^{x_n, z, \delta}_t\big) + \lim_{n \to \infty} Z^{x_n, z, \delta}_t
= \lim_{n \to \infty} Z^{x_n, z, \delta}_t = Z^{x, z, \delta}_t.
\]
\end{proof}

\begin{theorem}\label{thm:odesolutionconstruction}
\leavevmode
\begin{enumerate}
\item The map $(x,z) \mapsto P^\delta_t f(x,z) \coloneqq \E[f(X^x_t, Z^{x,z,\delta}_t)]$ is continuous for any $t \ge 0$ and any bounded continuous function $f$.
\item The map $(x,z) \mapsto P^\delta_t f(x,z)$ is Borel measurable for any $t \ge 0$ and any bounded Borel function $f$.
\item The process $(X^x_{t},Z^{x,z,\delta}_{t})_{t\geq 0}$ is a strong Markov process, i.e., for any finite $(\cF_t)$-stopping time $\eta$, $t \ge 0$ and a bounded Borel function $f$, we have 
\[
\E\big[f(X^x_{\eta+t}, Z^{x,z,\delta}_{\eta + t}) \big| \cF_\eta\big] = P^{\delta}_t f(X^x_{\eta}, Z^{x,z,\delta}_{\eta}), \quad \P\as
\]
\end{enumerate}
\end{theorem}
\begin{proof}
For assertion (1), we use Lemma \ref{continuitywrtstart}, the aforementioned continuous dependence of $X^x_t$ on $x$ and the dominated convergence theorem. Item (2) is obtained from (1) via a standard monotone class argument. The proof of the strong Markov property of $(X_t, Z^\delta_t)$ is analogous to that in the proof of \cite[Prop.~4.6]{kwon} and, hence, skipped.
\end{proof}

We will now establish the monotonicity of $(Z^{x,z,\delta}_t)$ in $\delta$, which will allow us to take the limit $\delta \to 0$ and define the object which we will call the solution to \eqref{eqn:Z}.

\begin{lemma}\label{amono}
The process $(Z^{x,z,\delta}_t)$ is increasing in $\delta > 0$ in the sense that $\P \big(Z^{x,z,\delta_1}_t \le Z^{x,z,\delta_2}_t \ \forall t \ge 0\big) = 1$ for $0 < \delta_1 < \delta_2$.
\end{lemma}
\begin{proof}
We suppress $(x,z)$ in the notation. Fix $\delta_2 > \delta_1 > 0$. By the comparison theorem for ODEs \cite[Theorem III-5-1]{ode}, we have that $Z^{\epsilon,\delta_1}_t \le Z^{\epsilon,\delta_2}_t$ for each $\epsilon > 0$ and any $t \ge 0$. Then we can take the limit as $\epsilon \downarrow 0$ on both sides and conclude that $Z^{\delta_1}_t \le Z^{\delta_2}_t$ $\P$\as The inequality holds for rational $t$ outside of the same $\P$-measure zero set. Due to the continuity of both processes, it is extended to all real $t \ge 0$ outside of the common null set.
\end{proof}

The monotonicity in $\delta$ established in Lemma \ref{amono} allows us to take the limit $\delta \to 0$. We fix a sequence $\delta_n \downarrow 0$ over which we define the limit:
\begin{equation}\label{eqn:Z_def}
Z^{x,z}_t := \lim_{n \to \infty} Z^{x,z,\delta_n}_t = \inf_{n} Z^{x,z,\delta_n}_t.
\end{equation}
This is a decreasing limit of continuous processes, so $(Z^{x,z}_t)$ has upper semicontinuous trajectories. Those trajectories are decreasing in $t$ which implies that they are \caglad. The right continuity may fail, for example, at time $0$ as we shall see later.

\begin{theorem}\label{thm:markov_Z}
Define $P_t f(x,z) = \E[f(X^x_t, Z^{x,z}_t)]$, $t \ge 0$, for a Borel-measurable bounded function $f$. The mapping $(x,z) \mapsto P_t f(x,z)$ is Borel measurable. Furthermore, the process $(X_{t},Z_{t})$ satisfies the strong Markov property.
\end{theorem}
\begin{proof}
By the dominated convergence theorem, we have, for any bounded and continuous $f$,
\[
\lim_{n \to \infty} \E[f(X^x_t, Z^{x,z, \delta_n}_t)] = \E[f(X^x_t, Z^{x,z}_t)].
\]
Hence, $P_t f = \lim_{n\to \infty} P^{\delta_n}_t f$ and $P_t f$ is a Borel function as a pointwise limit of Borel functions (see Theorem \ref{thm:odesolutionconstruction}). Borel measurability of $P_t f$ for bounded Borel $f$ is achieved by a monotone class argument.

For the strong Markov property, we start from studying a unique Carathéodory solution to the ODE:
\begin{equation}\label{eqn:barZ_delta_u}
        d\Bar{Z}^{\delta, \phi}_{t}=-\Bar{Z}^{\delta, \phi}_{t}\lambda^{\delta}(X_{u+t}, Y(\Bar{Z}^{\delta, \phi}_{t}))dt,\qquad \Bar{Z}^{\delta,\phi}_{0}=\phi,
\end{equation}
constructed in the same way as in the proof of Lemma \ref{lem:odesolutionconstruction} with a random variable $\phi \in [0,1]$. Notice that the arguments used in the proof of Lemma \ref{lem:z_lipsch} are $\omega$ by $\omega$, so they apply to this more general setting and yield
\begin{equation}\label{eqn:z_phi}
|\Bar{Z}^{\delta, \phi}_t - \Bar{Z}^{\delta, \phi'}_t| \le |\phi - \phi'|,
\end{equation}
for any random variables $\phi, \phi'$ with values in $[0, 1]$. 

Let us now establish the Markov property of $(Z_t)$. Fix an arbitrary $u > 0$. Let $(\Bar{Z}^{\delta}_{t}) := (\Bar{Z}^{\delta, Z_u}_t)$ and denote $\Bar{Z}_{t}=\lim_{n \to \infty}\Bar{Z}^{\delta_n}_{t}$. The Markov property of $(Z_t)$ is equivalent to showing that for any $s\geq 0$, we have $Z_{u+s}=\Bar{Z}_{s}$. This is not immediate as it is possible to have $Z^\delta_u>Z_u$ for all $\delta > 0$ which implies that $\Bar{Z}^\delta_t<Z^\delta_{u+t}$ for all small enough $t > 0$. 

By the uniqueness of solutions to \eqref{eqn:barZ_delta_u}, we have $Z^{\delta_n}_{u+t} = \Bar{Z}^{\delta_n, Z^{\delta_n}_u}_t$. Using \eqref{eqn:barZ_delta_u} and \eqref{eqn:z_phi}, we have
\[
|Z^{\delta_n}_{u+t} - \Bar{Z}^{\delta_n}_t| = |\Bar{Z}^{\delta_n, Z^{\delta_n}_u}_t - \Bar{Z}^{\delta_n, Z_u}_t| \le |Z^{\delta_n}_u - Z_u|.
\]
Taking $n\to\infty$ gives $Z_{u+t} = \Bar{Z}_t$ and completes the proof of the Markov property of $(Z_t)$. 

In the above arguments, we argued $\omega$ by $\omega$. Hence, the strong Markov property is proved analogously by replacing $u$ with a stopping time $\sigma$ and noting that $(Z^{\delta_n}_t)$ is a strong Markov process by Theorem \ref{thm:odesolutionconstruction}.
\end{proof}

Since the process $(X_t, Z_t)$ is Markovian, we extend the subscript notation $\P_{xz}$ to indicate the initial condition of both coordinates. 

\subsection{Properties of $(X_t, Z_t)$}

It is clear from the construction of $(Z_t)$ that, outside a set of $\P$ zero measure, it is defined as a monotone limit of $(Z^{\delta_n}_t)$ which satisfy \eqref{eqn:A_delta}. By the standard theory, $(X_t)$ is continuous apart from the exceptional set of measure zero. Finally, we also show in Lemma \ref{lem:Z_cont} that $(Z_t)$ is continuous for $t > 0$ apart from an exceptional set of $\P$-measure $0$. In what follows, we will not mention the zero measure set and write ``for all $\omega$''.

We have already noted that the difficulty with the definition of $(Z_t)$ was related to the explosion of $\lambda$ defined in \eqref{eqn:lambda} when $X_t$ approaches $\beta(Y(Z_t))$. We will now study how this affects the dynamics of the process $(X_t, Z_t)$. To this end we define a \emph{jump set}, where $\lambda = \infty$,
\begin{equation*}
  J:=\{(x,z) \in \cI \times (0, 1] :\ x<\beta(Y(z))\},
\end{equation*}
and a \emph{no-jump set}, where $\lambda$ is finite,
\begin{equation*}
  \NJ:=\{(x,z )\in \cI \times (0, 1] :\ x>\beta(Y(z))\}.
\end{equation*}
The reason for this terminology will become clear soon. In short, the process $(Z_t)$ jumps out of the region $J$ immediately and evolves continuously in $\NJ$. Notice that due to the continuity of $\beta$ and $Y$ both sets are open in $\cI \times (0,1]$. Let $\widebar{J}$ and $\widebar{\NJ}$ denote the closure of $J$ and $\NJ$ in $\cI \times (0,1]$. 

The following theorem summarises the properties of $(X_t, Z_t)$.

\begin{theorem}\label{thm:prop_Z}
Fix $(x, z) \in \cI \times [0,1]$. All statements are $\P_{xz}$\as
\begin{enumerate}
 \item For $t > 0$, the process $(X_t, Z_t)$ takes values in $\widebar\NJ \cup (\cI \times \{0\})$ and has continuous trajectories. If $(x,z) \in \widebar\NJ \cup (\cI \times \{0\})$, then $(Z_t)$ is also continuous at $t=0$.
 \item If $(x,z)\in J \cap ((\beta(\theta_L),\beta(\theta_U)]\times (0,1])$, then the process $(Z_t)$ experiences a jump at time $0$ and $Z_{0+} = Z(\beta^{-1}(x))$.
 \item If $x \le \beta(\theta_L)$, then $Z_{0+}=0$.
\end{enumerate}
\end{theorem}

Assertions (2)-(3) indicate that whenever $z > 0$ and $x < \beta(Y(z))$, there is an initial jump at time $t=0$. If $x > \beta(\theta_L)$, this jump goes to $Z(\beta^{-1}(x))$. Otherwise, the process $(Z_t)$ jumps to $0$. The proof of this theorem is delivered by Lemmas \ref{notimeintervallemma}--\ref{lem:Z_cont}.

\begin{lemma} \label{notimeintervallemma}
$\P_{xz}(\exists s>0\  \exists \epsilon>0\ \forall t\in[s,s+\epsilon]\ (X_t,Z_t)\in J)=0$.
\end{lemma}
\begin{proof}
We will argue $\omega$ by $\omega$ and omit $\omega$ in the notation. We assume, by contradiction, that there is a non-trivial time interval $[s,s+\epsilon]$ such that $(X_t(\omega),Z_t(\omega))\in J$ for all $t\in[s,s+\epsilon]$. Recall
\begin{equation}\label{eqn:4.16}
Z^\delta_{t}=Z^\delta_s-\int^{t}_{s}\frac{Z^\delta_u (rY({Z}^{\delta}_u)-D(X_u))}{(m(X_u;Y(Z_u^{\delta}))-Y(Z_u^{\delta}))\vee\delta}\ind{X_u\leq\alpha(Y( Z_u^{\delta}))}du
\end{equation}
and note that $\lim_{n\to\infty} Z^{\delta_n}_t = Z_t$ for $t \ge s$ and the limit is monotone from above.
This means that
\begin{equation}\label{eqn:4.17}
X_t < \beta(Y(Z_t)) \le \beta(Y(Z^{\delta_n}_t)) \le \alpha(Y(Z^{\delta_n}_t))\qquad \text{for $t \in [s, s+\epsilon]$ and all $n$,} 
\end{equation}
where we used that $\beta$ and $Y$ are increasing. Recall Proposition \ref{prop:dm}(ii) and Assumption \ref{ass:dm}(ii); they imply that $M(x') - r\theta \le 0$ for all $x' < \beta(\theta)$, which immediately yields
\[
r Y(Z^{\delta_n}_t)-D(X_t)\geq r Y(Z^{\delta_n}_t)-M(X_t)+(M-D)(X_t)\geq (M-D)(X_t) \qquad \text{for $t \in [s, s + \epsilon]$.}
\]
Hence, using that $(Z_t)$ is decreasing, the numerator under the integral in \eqref{eqn:4.16} is bounded from below by 
\[
Z_{s+\epsilon} \inf_{t \in [s, s + \epsilon]} (M - D)(X_t) =: \eta > 0.
\]
We also have $m(X_t; Y(Z^{\delta_n}_t)) = Y(Z^{\delta_n}_t)$ since $X_t \le \beta(Y(Z^{\delta_n}_t))$, $t \in [s, s + \epsilon]$, so the denominator under the integral equals $\delta_n$. Combining these observations shows that
\[
Z^{\delta_n}_{s + \epsilon} \le Z^{\delta_n}_s - \int_s^{s + \epsilon} \frac{\eta}{\delta_n} dt = Z^{\delta_n}_s - \frac{\eta}{\delta_n} \epsilon.
\]
Taking the limit $n \to \infty$ leads to $Z_{s + \epsilon} \le -\infty$ which contradicts $Z_{s + \epsilon} \ge 0$. 
\end{proof}

\begin{lemma}\label{tgreaterthan0lemma}
    $(X_t,Z_t) \notin J$ for all $t > 0$, $\P_{xz}$\as
\end{lemma}
\begin{proof}
We will argue $\omega$ by $\omega$ and omit $\omega$ in the notation. We fix $\omega$  such that $(X_t,Z_t) \in J$ for some $t>0$. By the definition of $J$, there is $\epsilon > 0$ such that $[X_t-\epsilon,X_t+\epsilon]\times[Z_t,1]\subset J$. We note that this holds even for $Z_t=1$. Due to the continuity of $t\mapsto X_t$, there must exist $\eta < t$ such that for $s\in[t-\eta,t]$ we have $X_s\in [X_t-\epsilon,X_t+\epsilon]$. We observe that $Z_t$ is decreasing. This means $(X_s,Z_s)\in J$ for $s\in[t-\eta,t]$. By the arbitrariness of $\omega$, we have
\[
\{(X_t,Z_t) \in J\text{ for some $t > 0$} \} \subseteq \{\exists s>0\ \exists \epsilon>0 \ (X_t,Z_t)\in J \text{ for all }t\in[s,s+\epsilon]\}.
\]
The latter set has $\P_{xz}$-probability $0$ by Lemma \ref{notimeintervallemma}. This completes the proof.
\end{proof}

\begin{lemma}\label{lem:3.14}
If $x\in (\beta(\theta_L),\beta(\theta_U)]$, then $Z_{0+}\leq Z(\beta^{-1}(x))$, $\P_{xz}\as$ for any $z \in [0,1]$. If $x \le \beta(\theta_L)$ then $Z_{0+} = 0$ $\P_{xz}$\as
\end{lemma}
\begin{proof} 
Due to Lemma \ref{tgreaterthan0lemma}, $(X_t, Z_t) \notin J$ for all $t > 0$ $\P_{xz}$\as By the definition of $J$, this means that either $Z_t > 0$ and $X_t \ge \beta(Y(Z_t))$ or $Z_t = 0$. 

Consider first $x\in (\beta(\theta_L),\beta(\theta_U)]$. The inverse $\beta^{-1}$ is well defined. When $Z_{0+} = 0$, the inequality $Z_{0+}\leq Z(\beta^{-1}(x))$ is trivially satisfied. When $Z_{0+} > 0$, we have $Z_t > 0$ for all sufficiently small $t$. As $(X_t, Z_t) \notin J$ for $t > 0$ by Lemma \ref{tgreaterthan0lemma}, we have $X_t \ge \beta(Y(Z_t))$.  Taking the limit $t \downarrow 0$ and using the continuity of $(X_t)$, we obtain $x \ge \beta(Y(Z_{0+}))$ which can be equivalently written as $Z(\beta^{-1}(x)) \ge Z_{0+}$. 

When $x \le \beta(\theta_L)$, we proceed by contradiction. Assume that $Z_{0+} > 0$, so $Z_t > 0$ for all sufficiently small $t$ and then arguing as in the previous paragraph we have $x \ge \beta(Y(Z_{0+}))$. As $Z_{0+} > 0$, we have $Y(Z_{0+}) > \theta_L$, which combined with the strict monotonicity of $\beta$ contradicts that $x \le \beta(\theta_L)$. 
\end{proof}

\begin{remark}
We note that $Z(\theta_U)=1$ and $Z_t\leq 1$. Hence, if we make the convention that $Z(\beta^{-1}(x))=1$ for all $x>\beta(\theta_U)$, the above lemma holds for all $x\in\cI$.
\end{remark}

\begin{lemma}
If $(x,z)\in \widebar{J}\cap ((\beta(\theta_L),\beta(\theta_U)]\times (0,1])$, then $(X_{0+},Z_{0+})=(x,Z(\beta^{-1}(x)))$ $\P_{xz}$\as
    \label{0lemma}
\end{lemma}
\begin{proof}
 We will argue $\omega$ by $\omega$ omitting $\omega$ in the notation. We have $X_{0+}=x$ due to the continuity of $(X_t)$. By Lemma \ref{lem:3.14}, $Z_{0+}\leq Z(\beta^{-1}(x))$. The proof will now follow by contradiction. Assume that $Z_{0+}<Z(\beta^{-1}(x))$ so that $(x, Z_{0+}) \in \NJ \cup (\cI \times \{0\})$. The idea of the argument is that the function $\lambda$ is finite on $\NJ$, so the process $(Z_t)$ cannot be instantaneously shifted within $\NJ$. Formal details follow.  

We claim that there exists $\epsilon \in (0, Z(\beta^{-1}(x)) - Z_{0+})$ such that $[x-\epsilon,x+\epsilon]\times(0,Z_{0+}+\epsilon]\subset \NJ$. When $Z_{0+} > 0$, we have $(x, Z_{0+}) \in \NJ$, so this follows from the openness and the form of $\NJ$. When $Z_{0+} = 0$, the assumption that $x > \beta(\theta_L)$ gives $Z(\beta^{-1}(x)) > 0$. Noticing that $z \mapsto \beta(Y(z))$ is increasing and continuous, the set $[x-\epsilon, x+\epsilon] \times \big(0, Z(\beta^{-1}(x-\epsilon))\big)$ is non-empty and contained in $\NJ$ for a sufficiently small $\epsilon$.

By the continuity of $(X_t)$, there exists $t^*>0$ such that $X_t\in [x-\epsilon,x+\epsilon]$ for $t\le t^*$. We note that $Z_0=Z_0^\delta\ge Z(\beta^{-1}(x)) >Z_{0+}+\epsilon$ by the definition of $\widebar{J}$ and the choice of $\epsilon$. Fix any $t < t^*$. By the definition of $Z_t$, there is $N$ such that $Z^{\delta_n}_t < Z_{t}+\epsilon$ for all $n \ge N$. Recalling that $Z_t \le Z_{0+}$, we have $Z^{\delta_n}_t < Z_{0+}+\epsilon$. Since $Z^{\delta_n}_0 > Z_{0+} + \epsilon$ and $(Z^{\delta_n}_t)$ is continuous, there is $u_n \in (0, t)$ such that $Z^{\delta_n}_{u_n} = Z_{0+}+\epsilon$. We have
\begin{align*}
Z^{\delta_n}_t&=Z^{\delta_n}_{u_n}-\int^{t}_{u_n}Z^{\delta_n}_s \frac{r Y(Z^{\delta_n}_s)-D(X_s)}{(m(X_s; Y(Z^{\delta_n}_s))- Y(Z^{\delta_n}_s))\vee\delta_n}\ind{X_s\leq\alpha( Y(Z^{\delta_n}_s))}ds\\
&\ge Z_{0+}+\epsilon-\frac{(t-u_n)c_2}{\rho \vee \delta_n}\\
&\ge Z_{0+}+\epsilon-\frac{t c_2}{\rho},
\end{align*}
where 
\begin{align*}
c_2 &= \sup_{(x',z')\in[x-\epsilon,x+\epsilon]\times(0,Z_{0+}+\epsilon]} z'\big(r Y(z') - D(x')\big) \ind{x'\leq\alpha(Y(z'))} < \infty, \\
\rho &= \inf_{(x',z')\in[x-\epsilon,x+\epsilon]\times(0,Z_{0+}+\epsilon]} m(x'; Y(z'))- Y(z') > 0.
\end{align*}
The fact that $c_2 < \infty$ is immediate by the continuity as the set over which the supremum is taken can be extended to a compact $[x-\epsilon,x+\epsilon]\times[0,Z_{0+}+\epsilon]$ as the expression is $0$ for $z' = 0$. To show that $\rho > 0$, we recall Lemma \ref{decreasem} which implies that $\rho = \inf_{x'\in[x-\epsilon,x+\epsilon]} m(x'; Y(Z_{0+}+\epsilon))- Y(Z_{0+}+\epsilon)$. The infimum is now taken over a compact set from a continuous function that is strictly positive on $\NJ$, so it is strictly positive itself.

Taking the limit $n \to \infty$, the above estimate yields $Z_t \ge Z_{0+}+\epsilon-\frac{t c_2}{\rho}$. By the arbitrariness of $t > 0$, we obtain
$Z_{0+} = \lim_{t \downarrow 0} Z_t \ge Z_{0+} + \epsilon$, a contradiction. This means that $(x, Z_{0+}) \notin \NJ$, which together with $Z_{0+}\leq Z(\beta^{-1}(x))$ completes the proof. 
\end{proof}
\begin{remark}
We note that for $x\le\beta(\theta_L)$, we have that $(X_{0+},Z_{0+})=(x,0)$ due to Lemma \ref{lem:3.14}. For $x>\beta(\theta_U)$, it is impossible that $(x,z)\in \widebar{J}$ due to the definition of $J$.
\end{remark}

\begin{lemma}\label{lem:Z_cont}
The process $(Z_t)$ is continuous $\P_{xz}$\as apart from time $0$. If $(x,z) \in \overline\NJ \cup (\cI \times \{0\})$, then $(Z_t)$ is also $\P_{xz}$\as continuous at $t=0$.
\end{lemma}
\begin{proof}
The process $(Z_t)$ is decreasing, so there can only be countably many discontinuities. It is also \caglad, so it suffices to show that it is impossible to have $Z_{t+} < Z_t$ for any $t > 0$. The proof will go by contradiction and $\omega$ by $\omega$. As before, we omit $\omega$ in the notation. We draw attention to the similarity of arguments to those in the proof of Lemma \ref{0lemma}. 

Assume, by contradiction, that $Z_{t+} < Z_t$. Surely, $X_t \le \beta(\theta_U)$ as otherwise the dynamics of $(Z_t)$ must be continuous by construction. We first discuss the case $X_t>\beta(\theta_L)$. By Lemma \ref{tgreaterthan0lemma}, we have $Z_t\leq Z(\beta^{-1}(X_t))$. Hence, $Z_{t+}<Z(\beta^{-1}(X_t))$ and so $(X_t, Z_{t+}) \in \NJ \cup (\cI \times \{0\})$. Arguing as in the proof of Lemma \ref{0lemma}, there exists $\epsilon \in (0, Z_t - Z_{t+})$ such that $[X_t-\epsilon,X_t+\epsilon]\times(0,Z_{t+}+\epsilon]\subset \NJ$. There is also $t^* > t$ such that $X_s \in [X_t-\epsilon,X_t+\epsilon]$ for all $s \in [t, t^*]$. Fix $s \in (t, t^*)$. By the construction of $(Z_t)$ there is $N$ such that $Z^{\delta_n}_{s} < Z_{t+} + \epsilon$ for all $n \ge N$. Since $Z^{\delta_n}_t \ge Z_t > Z_{t+} + \epsilon$, there is $u_n \in (t, s)$ such that $Z^{\delta_n}_{u_n} = Z_{t+} + \epsilon$.  We have
\begin{align*}
Z^{\delta_n}_{s}&=Z^{\delta_n}_{u_n}-\int^{s}_{u_n}Z^{\delta_n}_u \frac{r Y(Z^{\delta_n}_u)-D(X_u)}{(m(X_u; Y(Z^{\delta_n}_u))- Y(Z^{\delta_n}_u))\vee\delta_n}\ind{X_u\leq\alpha( Y(Z^{\delta_n}_u))}du\\
&\ge Z_{t+}+\epsilon-\frac{(s-u_n)c_2}{\rho \vee \delta_n}\\
&\ge Z_{t+}+\epsilon-\frac{(s-t) c_2}{\rho},
\end{align*}
where 
\begin{align*}
c_2 &= \sup_{(x',z')\in[X_t-\epsilon,X_t+\epsilon]\times(0,Z_{t+}+\epsilon]} z'\big(r Y(z') - D(x')\big) \ind{x'\leq\alpha(Y(z'))} <\infty, \\
\rho &= \inf_{(x',z')\in[X_t-\epsilon,X_t+\epsilon]\times(0,Z_{t+}+\epsilon]} m(x'; Y(z'))- Y(z') > 0.
\end{align*}
Sending $n \to \infty$ in the above inequality gives $Z_s \ge Z_{t+}+\epsilon-\frac{(s-t) c_2}{\rho}$, which implies, by the arbitrariness of $s$, $Z_{t+} \ge Z_{t+}+\epsilon$, a contradiction.

Let $\tau = \inf\{t \ge 0: X_t \le \beta(\theta_L)\}$. By Blumenthal's zero-one law, $\P_{xz}(\tau = 0) \in \{0,1\}$. If $\tau = 0$ $\P_{xz}$\as then $Z_{0+} = 0$ by Lemma \ref{lem:3.14} and $Z_t = 0$ for $t > 0$. Otherwise, Lemma \ref{tgreaterthan0lemma} implies that $(X_t, Z_t) \notin J$ for $t > 0$. If $Z_t = 0$ then trivially $Z_{t+} = 0$. If $Z_t > 0$ and $t < \tau$, then $Z_{t+} = Z_t$ by the first part of the proof as $X_t > \beta(\theta_L)$. Finally, we show that it is impossible to have $Z_t > 0$ for $t \ge \tau$. Indeed, this would mean $Z_{\tau} > 0$, $X_{\tau} = \beta(\theta_L)$ so $(X_{\tau}, Z_{\tau}) \in J$, a contradiction with the assertion of Lemma \ref{tgreaterthan0lemma}.

When $(x,z)\in \NJ$, we have $x > \beta(Y(z))$ so the process $(Z_t)$ is continuous at $0$ by construction. When $z = 0$, we have $Z_t = 0$ for all $t \ge 0$, and the continuity is trivial. It remains to study $(x,z) \in \cI \times (0,1]$ such that $x = \beta(Y(z))$. In this case, the continuity at $0$ is proved by contradiction in an analogous way as the continuity for $t > 0$. Taking $t=0$, we notice that $Z_t = z = Z(\beta^{-1}(x)) = Z(\beta^{-1}(X_t))$ from $x = \beta(Y(z))$. The remaining arguments can be repeated verbatim. 
\end{proof}


We finish this section with an easy auxiliary result which shows that the process $(Z_t)$ satisfies the dynamics \eqref{eqn:Z} in the Carath\'{e}odory sense when $(X_t, Z_t)$ is in $\NJ$. 

\begin{lemma}\label{lem:dynamics_Z}
Fix $\eta > 0$ and define a random set 
\[
A(\omega) = \{t \in (0, \infty):\ X_t(\omega) - \beta(Y(Z_t(\omega))) \ge \eta\}. 
\]
For any $\omega$ and any interval $[t_1, t_2] \subset A(\omega)$, we have
\begin{equation}\label{eqn:dynamics_Z}
Z_{t_2} = Z_{t_1} - \int_{t_1}^{t_2} Z_s \frac{r Y(Z_s)-D(X_s)}{m(X_s; Y(Z_s))- Y(Z_s)}\ind{X_s\leq\alpha( Y(Z_s))}ds.
\end{equation}
\end{lemma}
\begin{proof}
We will suppress $\omega$ in the notation for the sake of clarity. Take $[t_1, t_2] \subset A$. By the joint continuity of $(x, z) \mapsto m(x; Y(z))$, we have $\xi := \inf_{s \in [t_1, t_2]} \big(m(X_s; Y(Z_s))- Y(Z_s) \big) > 0$; $\xi = 0$ is impossible as, by continuity, the value would be attained at some $s \in [t_1, t_2]$ and hence $(X_s, Z_s)$ would belong to the stopping set, which would contradict $s \in A$. By Lemma \ref{lem:Z_cont}, the process $(Z_t)$ is continuous on $[t_1, t_2]$ as $t_1 > 0$ from the definition of $A$. By Dini's theorem, $(Z^{\delta_n}_t)$ converges uniformly to $(Z_t)$ on $[t_1, t_2]$. We can therefore fix $N$ such that $\inf_{s \in [t_1, t_2]}\big\{m\big(X_{s}; Y(Z^{\delta_N}_{s}) \big) - Y(Z^{\delta_N}_{s})\big\} > \xi/2$. Notice that $\inf_{s \in [t_1, t_2]} \big\{m\big(X_{s}; Y(Z^{\delta_n}_{s}) \big) - Y(Z^{\delta_n}_{s})\big\} > \xi/2$ for any $n \ge N$ due to Lemma \ref{decreasem} and the monotonicity of $n \mapsto Z^{\delta_n}_s$. Hence, for any $\delta_n < \xi/2$ and $n \ge N$, we have
\[
Z^{\delta_n}_{t_2} = Z^{\delta_n}_{t_1} - \int_{t_1}^{t_2} Z^{\delta_n}_s \frac{r Y(Z^{\delta_n}_s)-D(X_s)}{m(X_s; Y(Z^{\delta_n}_s))- Y(Z^{\delta_n}_s)}\ind{X_s\leq\alpha( Y(Z^{\delta_n}_s))}ds.
\]
Sending $n \to \infty$ and using the dominated convergence theorem we obtain \eqref{eqn:dynamics_Z}.
\end{proof}

Notice that the above lemma immediately implies that $(Z_t)$ is continuous at times $t > 0$ such that $(X_t, Z_t) \in \NJ$. This, however, does not guarantee that $(X_t, Z_t)$ does not jump from $\NJ$ to its boundary and does not address the continuity at the boundary. Hence, we could not have used Lemma \ref{lem:dynamics_Z} to simplify the proof of Lemma \ref{lem:Z_cont}.


\section{Best response value function, the stopping set and the equilibrium}\label{sec:value}

In this section, we study the best response optimal stopping problem 
\[
\sup_{\sigma \in \cT(\cF_t)} J(x, \sigma, \tau_i; \theta).
\] 
We refer the reader to Proposition \ref{prop:equil} to recall how it is related to the equilibrium of the game. We rewrite this problem as a Markovian optimal stopping problem with the state process $(X_t, Z_t)$. We will study the properties of the associated value function and of the stopping set. Finally, we will use them to construct an equilibrium of the game.

We start with integrating out the exit value of the opponent as it is independent from the filtration (observation) of the player. Recall that we assume that both players apply the strategy \eqref{eqn:def_hat_tau_Z} with the same process $(Z_t)$. Following similar arguments as in \cite[Lemma~2]{kwon}, we obtain
\begin{equation}\label{eqn:reformulation}
\begin{aligned}
J(x,\sigma,\tau_i;\theta) 
&= \E_{x1}\bigg[\int^\sigma_0e^{-rs}Z_sD(X_s)ds+\theta e^{-r\sigma}Z_{\sigma-}-\int_{[0,\sigma)}e^{-rs}m(X_s,\theta)dZ_s\bigg],
\end{aligned}
\end{equation}
where the subscript by the expectation indicates the initial point of the process $(X_t, Z_t)$: $(X_0, Z_{0}) = (x, 1)$.

\begin{remark}
On the right-hand side of \eqref{eqn:reformulation}, we integrate with respect to $dZ_t$. Recall that $(Z_t)$ is a \caglad process so, for each $\omega$, it induces a measure on $[0, \infty)$ and the integral is understood as a classical Lebesgue-Stieltjes integral. The reader may be more familiar with integration with respect to \cadlag functions but both are equivalent: the integral with respect to \caglad $(Z_t)_{t \ge 0}$ is the same as with respect to \cadlag $(Z_{t+})_{t \ge 0}$. 
\end{remark}

In order to maximise \eqref{eqn:reformulation} over $\sigma \in \cT(\cF_t)$, we define the corresponding value function
\begin{equation}\label{eqn:v_def}
    v(x,z;\theta)=\sup_{\sigma \in \cT(\cF_t)}\E_{xz}\bigg[\int^\sigma_0e^{-rs}Z_sD(X_s)ds+\theta e^{-r\sigma}Z_{\sigma-}-\int_{[0,\sigma)}e^{-rs}m(X_s,\theta)dZ_s\bigg].
\end{equation}
When $(x,z) \in \widebar\NJ$, the process $(Z_t)$ is continuous, so $v$ has an equivalent representation which is aligned with standard theory of optimal stopping of Markov processes:
\begin{equation}\label{eqn:v_def_NJ}
v(x,z;\theta)=\sup_{\sigma \in \cT(\cF_t)}\E_{xz}\bigg[\int^\sigma_0e^{-rs}Z_sD(X_s)ds+\theta e^{-r\sigma}Z_{\sigma}-\int_{[0,\sigma]}e^{-rs}m(X_s,\theta)dZ_s\bigg].
\end{equation}

In arguments below, we will use the following definition: a function $\varphi:\cI \times [0,1] \to \R$ is said to be $C_0$-continuous (cf. \cite[Ch.~3]{shiryaev}) if
\[
\lim_{t \to 0} \varphi(X_t, Z_t) = \varphi(x,z), \quad \P_{xz}\as
\]

\begin{theorem}\label{optimalstop}
On the set $\widebar \NJ$, the function $(x,z) \mapsto v(x,z; \theta)$ is $C_0$-continuous and the optimal stopping time is given by
\begin{equation}\label{eqn:tau_of_theta}
\tau(\theta) = \inf\{t \ge 0: v(X_t, Z_t;\theta) = Z_t \theta\}
\end{equation}
with $v(X_{\tau(\theta)}, Z_{\tau(\theta)}; \theta) = Z_{\tau(\theta)} \theta$ on $\{\tau(\theta) < \infty\}$.

On the set $J$, we have
\[
v(x,z; \theta) = \max(z \theta, v(x,Z_{0+}; \theta) + m(x;\theta)(z-Z_{0+}) \big).
\]
\end{theorem}
\begin{proof}
Fix $\theta \in [\theta_L, \theta_U]$, $(x,z) \in \widebar\NJ$ and recall \eqref{eqn:v_def_NJ}. Let $F(x,z)=\E_{xz}[\Psi]$, where
\[
\Psi = \int_{[0,\infty)}e^{-rs}Z_sD(X_s)ds-\int_{[0,\infty)}e^{-rs}m(X_s;\theta)dZ_s
\]
and
\begin{equation*}
 G(x,z)=\theta z - F(x,z).
\end{equation*}
By the strong Markov property of $(X_t, Z_t)$, we have
\begin{align*}
v(x,z;\theta)=\sup_{\sigma \in \cT(\cF_t)} \E_{xz}\big[e^{-r\sigma}G(X_\sigma, Z_\sigma)\big]+F(x,z),
\end{align*}
where we used the continuity of $(Z_t)$ $\P_{xz}\as$ since $(x,z) \in \widebar\NJ$ (see Thm.~\ref{thm:prop_Z}).

We first show important properties of $F$. Fix $(x,z) \in \widebar\NJ$. Notice that $\Psi$ is $\cF^W_\infty$-measurable, where $(\cF^W_t)$ is the natural filtration of the Brownian motion $(W_t)$. This is because $(Z_t)$ is a functional of $(W_t)$. Define $N_t = \E_{xz}[ \Psi | \cF^W_t]$. It is a square-integrable martingale as $\Psi$ is bounded. So it can be written as a stochastic integral with respect to $(W_t)$, which has two consequences. Firstly, it is continuous. Secondly, it is an $(\cF_t)$-martingale (as a stochastic integral with respect to $(\cF_t)$-Brownian motion), so for any $(\cF_t)$-stopping time $\eta$, we have $N_\eta = \E_{xz}[\Psi | \cF_\eta]$. It remains to relate $(N_t)$ with $F(X_t, Z_t)$. By the Markov property of $(X_t, Z_t)$, we have
\begin{equation}\label{eqn:G_and_N}
e^{-rt} F(X_t, Z_t) = N_t - \int_{[0,t)}e^{-rs}Z_sD(X_s)ds + \int_{[0,t)}e^{-rs}m(X_s;\theta)dZ_s.
\end{equation}
By the strong Markov property, the above identity holds for any $(\cF_t)$-stopping time $\eta$. As both sides of \eqref{eqn:G_and_N} are optional processes (cf. Theorem \ref{thm:markov_Z}) that coincide at stopping times, they are indistinguishable as stochastic processes.
Using the continuity of $(Z_t)$ (Thm.~\ref{thm:prop_Z}), we get the pointwise convergence of the integrals on the right-hand side to $0$ as $t \downarrow 0$. The continuity of $(N_t)$ then yields that $\lim_{t \downarrow 0} F(X_t, Z_t) = N_0 = F(x,z)$, which shows that $F$ is $C_0$-continuous for $(x,z) \in \widebar\NJ$. Thanks to the continuity of $(Z_t)$, the function $G$ is also $C_0$-continuous for $(x,z) \in \widebar\NJ$.

Theorem~3 in \cite[Ch.~3]{shiryaev} implies that $\hat v(x, z) = \sup_{\sigma \in \cT(\cF_t)} \E_{xz}\big[e^{-r\sigma}G(X_\sigma, Z_\sigma)\big]$ is $C_0$-continuous and the stopping time
\[
\sigma_\epsilon = \inf\{t \ge 0\colon \hat v(X_t, Z_t) \le G(X_t, Z_t) + \epsilon\}
\]
is $\epsilon$-optimal for any $\epsilon > 0$. Furthermore, $\hat v(X_{\sigma_\epsilon}, Z_{\sigma_\epsilon}) \le G(X_{\sigma_\epsilon}, Z_{\sigma_\epsilon}) + \epsilon$. To construct an optimal stopping time, we adapt the proof of the second part of \cite[Ch.~3, Thm.~3]{shiryaev} to our setting replacing the assumption of the upper semicontinuity of $G$ with a more direct argument. 

Let $\sigma_* = \lim_{\epsilon \downarrow 0} \sigma_\epsilon$ which exists and is a stopping time since the mapping $\epsilon \mapsto \sigma_\epsilon$ is decreasing. By $\epsilon$-optimality of $\sigma_\epsilon$ and Fatou's lemma, we have
\begin{equation}\label{eqn:hat-v}
\begin{aligned}
\hat v(x,z) 
&\le \limsup_{\epsilon \downarrow 0} \E_{xz}\big[e^{-r\sigma_\epsilon}G(X_{\sigma_\epsilon}, Z_{\sigma_\epsilon})\big]\\
&\le \E_{xz}\Big[ \limsup_{\epsilon \downarrow 0} e^{-r\sigma_\epsilon}G(X_{\sigma_\epsilon}, Z_{\sigma_\epsilon})\Big]\\
&= \E_{xz}\Big[ \ind{\sigma_* < \infty} \limsup_{\epsilon \downarrow 0} e^{-r\sigma_\epsilon}G(X_{\sigma_\epsilon}, Z_{\sigma_\epsilon})\Big].
\end{aligned}
\end{equation}
To evaluate the limit in the last expression, we write
\begin{align*}
&\ind{\sigma_* < \infty} |e^{-r\sigma_\epsilon}G(X_{\sigma_\epsilon}, Z_{\sigma_\epsilon}) - e^{-r\sigma_*}G(X_{\sigma_*}, Z_{\sigma_*})|\\
&=
\ind{\sigma_* < \infty} \bigg| N_{\sigma_*} - N_{\sigma_\epsilon} - \int_{\sigma_\epsilon}^{\sigma_*}e^{-rs}Z_sD(X_s)ds\\
&\hspace{55pt} +\int_{[\sigma_\epsilon,\sigma_*)}e^{-rs}m(X_s;\theta)dZ_s + \theta (e^{-r\sigma_\epsilon}Z_{\sigma_\epsilon}-e^{-r\sigma_*}Z_{\sigma_*})\bigg|,
\end{align*}
where we used \eqref{eqn:G_and_N}. The continuity of $(N_t)$ and $(Z_t)$ as well as the boundedness of functions under the integrals imply that the right-hand side converges to $0$ as $\epsilon \downarrow 0$.
This allows us to continue the estimates \eqref{eqn:hat-v} as follows
\[
\hat v(x,z) \le \E_{xz}\big[ \ind{\sigma_* < \infty} e^{-r\sigma_*}G(X_{\sigma_*}, Z_{\sigma_*})\big]
= \E_{xz}\big[ e^{-r\sigma_*}G(X_{\sigma_*}, Z_{\sigma_*})\big],
\]
where in the last equality we could drop the indicator function due to the boundedness of $G$. The above inequality means that $\sigma_*$ is an optimal stopping time.

It remains to prove that $\sigma_* = \sigma_0$. To this end, we repeat, verbatim, arguments from the proof of \cite[Ch.~3, Thm.~3]{shiryaev}:
\[
\E_{xz}\big[ e^{-r\sigma_*}\hat v(X_{\sigma_*}, Z_{\sigma_*})\big] 
\ge \E_{xz}\big[ e^{-r\sigma_*}G(X_{\sigma_*}, Z_{\sigma_*})\big]
= \hat v(x, z) \ge \E_{xz}\big[ e^{-r\sigma_*}\hat v(X_{\sigma_*}, Z_{\sigma_*})\big]
\]
with the last inequality by Thm.~1 and Lemma~1 in \cite[Ch.~3]{shiryaev}. This shows that 
\[
\E_{xz}\big[ e^{-r\sigma_*}\hat v(X_{\sigma_*}, Z_{\sigma_*})\big] = \E_{xz}\big[ e^{-r\sigma_*}G(X_{\sigma_*}, Z_{\sigma_*})\big].
\]
Combined with the inequality $\hat v \ge G$, it yields that $\hat v(X_{\sigma_*}, Z_{\sigma_*}) = G(X_{\sigma_*}, Z_{\sigma_*})$, $\P_{xz}$\as, i.e., $\sigma_* \ge \sigma_0$. The opposite inequality is true as $\sigma_\epsilon \le \sigma_0$. Hence, $\sigma_* = \sigma_0$ and so $\sigma_0$ is an optimal stopping time and, in the notation of the theorem, $v(X_{\tau(\theta)}, Z_{\tau(\theta)}; \theta) = Z_{\tau(\theta)} \theta$, $\P_{xz}$\as

Recall that $\hat v(x,z) = v(x,z; \theta) - F(x,z)$ and $G(x,z) = \theta z - F(x,z)$. Hence,
\[
\{t \ge 0\colon \hat v(X_t,Z_t) = G(X_t,Z_t) \} = \{t \ge 0\colon v(X_t,Z_t; \theta) = \theta Z_t \},
\]
so $\sigma_0 = \tau(\theta)$. 

Consider now $(x,z) \in J$. By Blumenthal's 0-1 law, 
\[
v(x,z; \theta) = \max \left( z \theta, \sup_{\sigma \in \cT(\cF_t),\ \sigma > 0} J(x,z, \sigma; \theta) \right).
\]

Similarly,
\[
v(x,Z_{0+}; \theta) = \max \Big( Z_{0+} \theta, \sup_{\sigma \in \cT(\cF_t),\ \sigma > 0} J(x,Z_{0+}, \sigma; \theta) \Big).
\]
By the continuity of $(Z_t)$, $\P_{xZ_{0+}}\as$ (cf. Thm.~\ref{thm:prop_Z}), we have $Z_{0+}\theta = \lim_{t \downarrow 0} J(x, Z_{0+}, t; \theta)$, so $v(x,Z_{0+}; \theta) = \sup_{\sigma \in \cT(\cF_t),\ \sigma > 0} J(x,Z_{0+}, \sigma; \theta)$. Noting that
\[
J(x,z, \sigma; \theta) = J(x,Z_{0+}, \sigma; \theta) + m(x;\theta)(z-Z_{0+})
\]
for $\sigma > 0$, we have
\begin{align*}
v(x,z; \theta) 
&= \max \Big( z \theta, \sup_{\sigma \in \cT(\cF_t),\ \sigma > 0} J(x,Z_{0+}, \sigma; \theta) + m(x;\theta)(z-Z_{0+}) \Big)\\
&= \max \big( z \theta, v(x, Z_{0+}; \theta) + m(x;\theta)(z-Z_{0+})\big).
\end{align*}
\end{proof}

\begin{remark}
The optimal stopping time $\tau(\theta)$ given by \eqref{eqn:tau_of_theta} is also optimal for $(x,z) \in J$. The proof of this fact is given in Lemma \ref{lem:tau_theta_J} in the next section.
\end{remark}

The structure of the stopping set for $v(x,z; \theta)$ is stated in the following theorem.
\begin{theorem}\label{thm:stopping-set}
For any $\theta \in (\theta_L, \theta_U)$, the stopping time
\[
\tau(\theta) = \inf\{t \ge 0\colon (X_t, Z_t) \in S(\theta) \},
\]
where 
\begin{equation*}
    S(\theta):=\bigg(\Big(\big(\beta(\theta),\alpha(\theta)\big]\times[0,Z(\theta)]\Big)\cup \Big((x_L,\beta(\theta)]\times(0,1]\Big)\cup \Big(\cI \times\{0\}\Big)\bigg) \cap \Big(\cI \times [0,1]\Big),
\end{equation*}
is optimal for $v(x,z;\theta)$.
\end{theorem}

Notice that the intersection with $\cI \times [0,1]$ in the definition of $S(\theta)$ only has an effect if $\alpha(\theta) = x_U$ and is there so that $S(\theta)$ is a subset of the state space of the process $(X_t, Z_t)$. We defer the proof of Theorem \ref{thm:stopping-set} until the next section. It consists in a number of technical results collected in Lemmas \ref{optimaliszero}-\ref{alphainclude}.

Although the optimal stopping time $\tau(\theta)$ is not in the form postulated heuristically as the equilibrium strategy of players, it can be shown to be equivalently written as \eqref{eqn:def_hat_tau_Z} as long as the initial value $z$ of the belief process is greater than or equal to $Z(\theta)$. We note that the strict inequality $D(x) < r\theta$ for $x < c_D(\theta)$ from Assumption \ref{ass:dm}(i) is instrumental for this result. 

\begin{lemma}\label{lem:stop-Z-transform}
For $\theta \in (\theta_L, \theta_U]$ and $(x,z)\in \I\times [Z(\theta),1]$, we have 
\begin{equation}\label{eqn:eq_sets}
     \inf\{t\geq 0\colon Z_t<Z(\theta)\}=\tau(\theta), \qquad \P_{xz}\as
\end{equation}
\end{lemma}
\begin{proof}
Fix $\theta$ and $(x,z)$ as in the statement of the lemma. Let $\tau = \inf\{t\geq 0:Z_t<Z(\theta)\}$, $\sigma=\inf\{t\geq 0\colon Z_t\le Z(\theta)\text{ and }X_t\leq\alpha(\theta)\}$, $\sigma'=\inf\{t\geq 0 \colon X_t\leq\beta(\theta)\}$ and $\sigma'' = \inf\{t \ge 0\colon Z_t = 0\}$. The equality \eqref{eqn:eq_sets} is equivalent to $\tau = \sigma \wedge \sigma' \wedge \sigma''$. Since $z > Z(\theta_L) = 0$, we have $\sigma'' \ge \sigma'$, so it is sufficient to study the equality $\tau = \sigma \wedge \sigma'$. We will show first that $\tau = \sigma$, $\P_{xz}\as$ and then prove that $\sigma' \wedge \sigma''\ge \sigma$, $\P_{xz}\as$

Consider the case $x\in(x_L,\beta(\theta)]$ and apply Lemma \ref{lem:3.14}. If $x \le \beta(\theta_L)$, we have $Z_{0+} = 0 < Z(\theta)$ and $\tau = \sigma = 0$. If $x \in (\beta(\theta_L), \beta(\theta))$, using that $Z(y)$ is a strictly increasing function and $\beta^{-1}$ is also strictly increasing (Lemma \ref{lem:beta}), we obtain $Z_{0+} \le Z(\beta^{-1}(x)) < Z(\theta)$, so again $\tau=\sigma=0$. The boundary case $x=\beta(\theta)$ is split further into two: $z > Z(\theta)$ or $z = Z(\theta)$. For $z = Z(\theta)$, we notice that $Z_t \le Z^\delta_t < Z^\delta_0 = Z(\theta)$ for $t > 0$ because $x < \alpha(\theta)$. Hence, $\tau = \sigma =0$. When $z > Z(\theta)$, we similarly have $Z_t < Z_{0+} \le Z(\theta) $ for $t > 0$, so $\tau = \sigma = 0$. 

When $x \in (\beta(\theta), \alpha(\theta)]$ and $z = Z(\theta)$ we argue similarly as above. Namely, if $x < \alpha(\theta)$, $Z_t \le Z^\delta_t < Z^\delta_0 = Z(\theta)$ for $t > 0$ (with the strict inequality by Assumption \ref{ass:dm}(i)). In the boundary case $x = \alpha(\theta)$, we employ the regularity of $\alpha(\theta)$ for the diffusion $(X_t)$, i.e., $\P_{xz}(\inf\{t \ge 0\colon X_t < \alpha(\theta) \} = 0 ) = 1$, which validates the above sequence of inequalities $Z_t < Z(\theta)$ for $t > 0$.

The rest of the proof will rest on the use of the strong Markov property of $(X_t, Z_t)$ and the fact, just proved, that $\P_{xz}(\tau = \sigma = 0) = 1$ for $(x,z) \in \cO_\theta=(x_L, \alpha(\theta)] \times \{Z(\theta)\}$. Take $x > \beta(\theta)$ and $z \in [Z(\theta), 1]$. Let $\eta = \inf\{t \ge 0\colon (X_t, Z_t) \in \cO_\theta\}$. By Lemma \ref{lem:Z_cont}, the process $(Z_t)$ is continuous for $t > 0$ $\P_{xz}$\as If $x > \beta(Y(z))$ then $Z_{0+} = z \ge Z(\theta)$, so the continuity extends to $t \ge 0$, whereas when $x \le \beta(Y(z))$ we have $Z_{0+} = Z(\beta^{-1}(x)) > Z(\theta)$ by Lemma \ref{0lemma}, so $\eta > 0$ and the discontinuity at $t=0$ does not influence the arguments below.  We obviously have $\tau \ge \eta$ and $\sigma \ge \eta$ $\P_{xz}$\as, so $\tau = \tau(\theta)=\infty$ on $\{\eta = \infty\}$. On the set $\{\eta < \infty\}$, we have $(X_\eta, Z_\eta) \in \cO_\theta$ $\P_{xz}$\as and, by the strong Markov property, $\P_{xz}(\tau = \sigma) = \E_{xz}[ \P_{X_\eta, Z_\eta}(\tau = \sigma)] \ge \E_{xz}[ \P_{X_\eta, Z_\eta}(\tau = \sigma = 0)] = 1$ with the last equality by $(X_\eta, Z_\eta) \in \cO_\theta$ and $\tau = \sigma = 0$ on $\cO_\theta$.

It remains to show that $\sigma'\geq\tau$, $\P_{xz}$\as As we have shown above that $\tau = \sigma$, $\P_{xz}$\as, we will demonstrate that $\sigma' \ge \sigma$, $\P_{xz}$\as By the continuity of $(X_t)$, we have $\P_{xz}(X_{\sigma'} \le \beta(\theta)) = 1$, so the strong Markov property of $(X_t)$ allows us to write
\[
\P_{xz}(\sigma' < \sigma) = \E_{xz}[\ind{\sigma' \le \sigma} \P_{X_{\sigma'}, Z_{\sigma'}}(\sigma > 0)] \le \E_{xz}[\P_{X_{\sigma'}, Z_{\sigma'}}(\sigma > 0)] = 0,
\]
where the last equality follows from $\P_{xz}(\sigma = 0) = 1$ for $x \le \beta(\theta)$, which we have shown at the beginning of the proof. 
\end{proof}

The main result of the paper is the equilibrium of the game.
\begin{theorem}\label{thm:equilibrium}
Define $\tau_i = \inf\{t\ge 0 \colon Z_t < Z(\theta_i)\}$. Then $(\tau_1,\tau_2)$ is a Nash equilibrium.
\end{theorem}
\begin{proof}
The proof rests on the verification of conditions of Proposition \ref{prop:equil}. Notice that $\tau_i = \hat\tau_i(\omega, \theta_i(\omega))$, where $\hat\tau_i(\omega, \theta) = \inf\{t\ge 0\colon Z_t(\omega) < Z(\theta) \}$. By the symmetry of the problem, it is sufficient to show that for every $\theta\in (\theta_L,\theta_U)$ and $x\in \cI$, the stopping time $\hat\tau_1(\cdot,\theta)$ solves the optimal stopping problem $\sup_{\sigma\in \cT(\cF_t)} J(x,\sigma,\tau_2;\theta)$; indeed, the set $\{\theta_L, \theta_U\}$ has $F$-measure zero as $F$ is continuous, so it is sufficient to show optimality of $\hat\tau_1(\cdot, \theta)$ for $\theta$ in the open interval $(\theta_L, \theta_U)$. Recall that by the definition of $v(x,z;\theta)$ we have 
\[
v(x,1;\theta)=\sup_{\sigma\in \cT(\cF_t)} J(x,\sigma,\tau_2;\theta).
\]
The assumption $Z(\theta) \le 1$ of Lemma \ref{lem:stop-Z-transform} is satisfied. By combining this finding with Theorem \ref{thm:stopping-set} it emerges that the optimal stopping time for $\sup_{\sigma\in \cT(\cF_t)} J(x,\sigma,\tau_2;\theta)$ is $\hat \tau_1(\cdot; \theta)$ as required. 
\end{proof}

\begin{remark}
The generality of assumptions of Thm.~\ref{thm:stopping-set} and Lemma \ref{lem:stop-Z-transform} implies that a version of Bayesian perfection in the spirit of \cite[Section~4.7]{kwon} can be established here.
\end{remark}

\section{Proof of Theorem \ref{thm:stopping-set}}\label{sec:proof-thm}

We will construct the stopping set separately on different regions of the state space. The methods of proof differ depending on the region. We will exploit that the optimal stopping sets for $v$ and $\Tilde{v}$ are identical, where
\begin{equation}\label{tildevexpression}
\begin{aligned}
\Tilde{v}(x,z;\theta)&= v(x,z;\theta) - z\theta\\
&= \sup_{\sigma \in \cT(\cF_t)} \E_{xz} \bigg[ \int_0^\sigma e^{-rs} Z_s \big(D(X_s) -r\theta\big) ds - \int_{[0,\sigma)} e^{-rs} \big(m(X_s; \theta) - \theta\big) dZ_s \bigg].
\end{aligned}
\end{equation}
Notice that for $(x,z) \in \widebar\NJ$ or when $\sigma > 0$, the continuity of the process $(Z_t)$ at $\sigma$ yields an equivalent formulation of the expectation in \eqref{tildevexpression} with the second integral over $[0, \sigma]$ which is aligned with the standard approach in stochastic analysis:
\[
\E_{xz} \bigg[ \int_0^\sigma e^{-rs} Z_s \big(D(X_s) -r\theta\big) ds - \int_{[0,\sigma]} e^{-rs} \big(m(X_s; \theta) - \theta\big) dZ_s \bigg].
\]
We will use both forms depending on which is more convenient for arguments.

\begin{lemma}\label{optimaliszero}
The stopping time $\sigma=0$ is optimal for $v(x,z;\theta)$ for $x\le\beta(\theta)$ and $\theta \in [\theta_L, \theta_U]$.
\end{lemma}
\begin{proof}
Reversing the integration of the exit value $\theta_i$ of the opponent, see \eqref{eqn:reformulation}, we obtain
    \begin{align*}
        &v(x,z;\theta)\\
        &=\sup_{\sigma \in \cT(\cF_t)}\E_{xz}\bigg[\int^\sigma_0e^{-rs}Z_sD(X_s)ds+\theta e^{-r\sigma}Z_{\sigma-}-\int_{[0,\sigma)} e^{-rs}m(X_s;\theta)dZ_s\bigg]\\
        &=z \sup_{\sigma \in \cT(\cF_t)}\E_{xz}\bigg[\int^\sigma_0e^{-rs}Z_s/z D(X_s)ds+\theta e^{-r\sigma}Z_{\sigma-}/z -\int_{[0,\sigma)} e^{-rs}m(X_s;\theta) d(Z_s/z)\bigg]\\
        &=z \sup_{\sigma \in \cT(\cF_t)}\E_{xz}\bigg[\int^{\sigma\wedge\hat\tau(\theta_i)}_0e^{-rs}D(X_s)ds+\theta e^{-r\sigma}\ind{\sigma\leq\hat\tau(\theta_i)}+\ind{\sigma>\hat\tau(\theta_i)}e^{-r\hat\tau(\theta_i)}m(X_{\hat\tau(\theta_i)};\theta)\bigg],
    \end{align*}
where $\hat\tau(y) = \inf\{t \ge 0\colon Z_t/z < Z(y)\}$ is an $(\cF_t)$-stopping time for each $y \in (\theta_L, \theta_U)$. Recalling that $m(x; \theta) \ge \theta$, we obtain the first inequality
    \begin{align*}
        v(x,z;\theta) &\le z \sup_{\sigma \in \cT(\cF_t)}\E_{xz}\bigg[\int^{\sigma\wedge\hat\tau(\theta_i)}_0e^{-rs}D(X_s)ds + e^{-r(\sigma \wedge \hat\tau(\theta_i))} m(X_{\sigma\wedge\hat\tau(\theta_i)};\theta)\bigg]\\
        &\le z \sup_{\sigma \in \cT(\cF_t)}\E_{xz}\bigg[\int^{\sigma\wedge\hat\tau(\theta_i)}_0e^{-rs}M(X_s)ds + e^{-r(\sigma \wedge \hat\tau(\theta_i))} m(X_{\sigma\wedge\hat\tau(\theta_i)};\theta)\bigg]\\
        &\le zm(x;\theta),
    \end{align*}
while the second inequality is because $M \ge D$ and the final inequality is by the tower property (conditioning of $\cF_\infty$) and the dynamic programming principle for the optimal stopping problem $m(x;\theta)$. Since $x \le \beta(\theta)$, we have $m(x;\theta) = \theta$, so $v(x, z; \theta) \le z \theta$. The upper bound $z\theta$ for $v$ is attained by $\sigma=0$.
\end{proof}

\begin{lemma}\label{conti}
$\Tilde{v}(x,z;\theta)>0$ for $(x,z)\in (\alpha(\theta), x_U) \times (0,1]$ for any $\theta \in [\theta_L, \theta_U]$.
\end{lemma}
\begin{proof}
Recall \eqref{eqn:a1}. Then $\Tilde{d}(x;\theta) = d(x; \theta) - \theta$ satisfies for every $x \in (\alpha(\theta), x_U)$
\[
(\cL - r)\Tilde{d}(x; \theta) + D(x) - r\theta = 0.
\]
Fix $x>\alpha(\theta)$ and $0<z\leq1$. Let $\sigma = \inf\{t \ge 0\colon X_t \le \alpha(\theta)\}$. Using the above equality and the standard localisation procedure to integrate out the martingale part from It\^o formula,  we have
\begin{align*}
&\E_{xz}\big[e^{-r\sigma} Z_\sigma \Tilde{d}(X_\sigma; \theta)\big] - z \Tilde{d}(x; \theta)\\
&=
\E_{xz}\bigg[\int_0^\sigma e^{-rs} Z_s (\cL - r) \Tilde d(X_s; \theta) ds + \int_{[0,\sigma]} e^{-rs} \Tilde d(X_s; \theta)dZ_s \bigg]\\
&=
\E_{xz}\bigg[\int_0^\sigma e^{-rs} Z_s \big(r\theta - D(X_s)\big) ds + \int_{[0,\sigma]} e^{-rs} \Tilde d(X_s; \theta)dZ_s \bigg].
\end{align*}
Recall that $m \ge d$ and $(Z_t)$ is decreasing. This yields the first inequality in 
\begin{align*}
\Tilde{v}(x,z;\theta)&\ge
\E_{xz} \bigg[ \int_0^\sigma e^{-rs} Z_s \big(D(X_s) -r\theta\big) ds - \int_{[0,\sigma]} e^{-rs} \Tilde d(X_s; \theta)dZ_s \bigg]\\
&= \E_{xz} \big[z \Tilde{d}(x; \theta) - e^{-r\sigma} Z_\sigma \Tilde{d}(X_\sigma; \theta)\big] = z \Tilde{d}(x; \theta) > 0,
\end{align*}
where in the second equality we use the continuous fit $\tilde d(\alpha(\theta); \theta) = 0$ and the fact that $e^{-r\sigma} Z_\sigma \Tilde{d}(X_\sigma; \theta) = 0$ on $\{\sigma = \infty\}$; the last inequality is because $x > \alpha(\theta)$ is in the continuation set of $\Tilde{d}(\cdot; \theta)$. 
\end{proof}

\begin{lemma}\label{lem7.3}
For any $\theta \in [\theta_L, \theta_U]$, we have $\Tilde{v}(x,z;\theta)>0$ for $(x,z)\in \I\times(Z(\theta),1]$ such that $\beta(Y(z))<x<\alpha(Y(z))$.
\end{lemma}
\begin{proof}
For any $(x,z)\in I\times(0,1]$ satisfying $\beta(Y(z))<x<\alpha(Y(z))$ and $z>Z(\theta)$, we have:
\begin{equation}\label{eqn:integralgreater0}
\begin{aligned}
&D(x)-r\theta+\lambda(x,Y(z))\big(m(x;\theta)-\theta\big)\\
&=D(x)-r\theta+\frac{rY(z)-D(x)}{m(x;Y(z))-Y(z)}\big(m(x;\theta)-\theta\big)\\
&>D(x)-r\theta+\frac{rY(z)-D(x)}{m(x;\theta)-\theta}\big(m(x;\theta)-\theta\big) = r(Y(z) - \theta),
\end{aligned}
\end{equation}
where we used that $m(x;Y(z))-Y(z)< m(x;\theta)-\theta$ for $z > Z(\theta)$ by Lemma \ref{decreasem}; indeed, $x > \beta(Y(z)) \ge \beta(\theta)$ as $\beta$ is increasing by Lemma \ref{lem:beta}.

Fix $(x,z)$ and $\theta$ as in the statement of the lemma. Let $0 < \delta < x - \beta(Y(z))$ and take a stopping time
\[
\eta=\inf\big\{t\geq 0\colon X_t<\beta(Y(Z_t)) + \delta \text{ or }X_t\geq\alpha(Y(Z_t))\big\} \wedge \inf\{t \ge 0\colon Z_t \le Z(\theta)\}.
\]
Since $(x,z) \in \NJ$, the process $(X_t, Z_t)$ is continuous and so $\P_{xz}(\eta > 0) = 1$.
By Lemma \ref{lem:dynamics_Z}, \eqref{tildevexpression} and \eqref{eqn:integralgreater0}, we obtain
\begin{align*}
\Tilde{v}(x,z;\theta)
&\geq\E_{xz}\bigg[\int^\eta_0e^{-rs}Z_s\big(D(X_s)-r\theta\big)ds-\int_{[0, \eta]}e^{-rs}(m(X_s,\theta)-\theta)dZ_s\bigg]\\
&=\E_{xz}\bigg[\int^\eta_0 e^{-rs} Z_s \Big(D(X_s)-r\theta+\lambda(X_s,Y(Z_s))\big(m(X_s,\theta)-\theta\big)\Big)ds\bigg]>0.
\end{align*}
\end{proof}

\begin{lemma}\label{continuationsetcomp}
For any $\theta \in [\theta_L, \theta_U]$, we have $\Tilde{v}(x,z;\theta)>0$ for $(x,z)\in \I\times(Z(\theta),1]$ such that $\beta(\theta)<x\le\beta(Y(z))$.
\end{lemma}
\begin{proof}
Take $(x,z)$ and $\theta$ as in the statement of the lemma. As $(x, z) \in J$, the process $(Z_t)$ has a discontinuity at $t=0$ and it is continuous afterwards, see Thm.~\ref{thm:prop_Z}.

Consider the case $x < \beta(Y(z))$ first. Let
\[
\eta = \inf \{t \ge 0\colon X_t > \alpha(Y(Z_t)) \text{ or } Z_t \le Z(\theta) \}.
\]
By Lemma \ref{0lemma}, $Z_{0+} = Z(\beta^{-1}(x)) < z$. From $\beta(\theta) < x$, we get $\theta < \beta^{-1}(x)$, so $Z(\theta) < Z(\beta^{-1}(x)) = Z_{0+}$ as both functions $\beta$ and $Z$ are increasing. Certainly $x < \alpha(\beta^{-1}(x))$ as $\alpha(\cdot) > \beta(\cdot)$ by Lemma \ref{lem:alpha_beta_ineq}. Hence, $\P_{xz}(\eta > 0) = 1$. 

For any $\epsilon > 0$, we have
\begin{align*}
\Tilde{v}(x,z; \theta) 
&\ge \E_{xz}\bigg[ \int^{\eta\wedge \epsilon}_{0} e^{-rs}Z_s\big(D(X_s)-r\theta\big)ds-\int_{[0, \eta \wedge \epsilon]} e^{-rs}\big(m(X_s;\theta)-\theta\big)dZ_s\bigg]\\
&=  \E_{xz}\bigg[ \int^{\eta\wedge \epsilon}_{0} e^{-rs}Z_s\big(D(X_s)-r\theta\big)ds-\int_{(0, \eta \wedge \epsilon]} e^{-rs}\big(m(X_s;\theta)-\theta\big)dZ_s \bigg]\\
&\quad + \big(m(x;\theta)-\theta\big)\big(z -Z(\beta^{-1}(x)) \big) .
\end{align*}
By the dominated convergence theorem, when $\epsilon \downarrow 0$, the expectation above vanishes yielding
\[
\Tilde{v}(x,z; \theta) \ge \big(m(x;\theta)-\theta\big)\big(z -Z(\beta^{-1}(x)) \big) > 0,
\]
where in the last inequality we used that $z > Z(\beta^{-1}(x))$ and $m(x; \theta) > \theta$ because $x > \beta(\theta)$.

It remains to study $x = \beta(Y(z))$. Recall that in this case $(Z_t)$ is continuous for $t\ge 0$, see Lemma \ref{lem:Z_cont}. Consider
\begin{align*}
\eta = \inf \{t \ge 0\colon &m(X_t; Y(Z_t)) - Y(Z_t) > \epsilon/2 \text{ or } m(X_t; \theta) - \theta < \epsilon/2 \\
&\text{or } X_t \ge \alpha(Y(Z_t)) \text{ or } Z_t \le Z(\theta) \},
\end{align*}
where $\epsilon = m(x;\theta)-\theta > 0$. Referring to the separation of boundaries $\alpha$ and $\beta$ (Lemma \ref{lem:alpha_beta_ineq}), we have $x = \beta(Y(z)) < \alpha(Y(z))$. Noting that $m(x; Y(z)) - Y(z) = 0$ as $x = \beta(Y(z))$ and recalling the continuity of $t \mapsto  m(X_t; Y(Z_t)) - Y(Z_t)$ and of $x \mapsto m(x; \theta) - \theta$ (Lemma \ref{lem:beta}), we conclude that $\P_{xz}(\eta > 0) = 1$. Consider a lower bound
\begin{equation}\label{eqn:ttv}
\Tilde{v}(x,z; \theta) \ge \E_{xz}\bigg[ \int^{\eta}_{0} e^{-rs}Z_s\big(D(X_s)-r\theta\big)ds-\int_{[0, \eta]} e^{-rs}\big(m(X_s;\theta)-\theta\big)dZ_s\bigg].
\end{equation}
We will justify at the end of the proof that
\[
-dZ_s \ge Z_s (rY(Z_s) - D(X_s)) 2/\epsilon\, ds.
\]
This enables a bound of the second term under the expectation
\[
\E_{xz}\bigg[-\int_{[0, \eta]} e^{-rs}\big(m(X_s;\theta)-\theta\big)dZ_s\bigg]
\ge \E_{xz}\bigg[\int_0^\eta e^{-rs}\big(m(X_s;\theta)-\theta\big)Z_s(rY(Z_s) - D(X_s)) 2/\epsilon\, ds\bigg].
\]
Since $X_s < \alpha(Y(Z_s))$ for $s < \eta$, we have $rY(Z_s) - D(X_s) > 0$ (Assumption \ref{ass:dm}). We also have $m(X_s; \theta) - \theta > \epsilon/2$, so
\[
\E_{xz}\bigg[-\int_{[0, \eta]} e^{-rs}\big(m(X_s;\theta)-\theta\big)dZ_s\bigg]
\ge
\E_{xz}\bigg[\int_0^\eta e^{-rs}Z_s(rY(Z_s) - D(X_s)) ds\bigg].
\]
We insert this bound into \eqref{eqn:ttv} and combine both integrals
\begin{align*}
&\E_{xz}\bigg[ \int^{\eta}_{0} e^{-rs}Z_s\big(D(X_s)-r\theta\big)ds-\int_{[0, \eta]} e^{-rs}\big(m(X_s;\theta)-\theta\big)dZ_s\bigg]\\
&\ge
\E_{xz}\bigg[ \int^{\eta}_{0} e^{-rs}Z_s\big(D(X_s)-r\theta\big)ds+ \int_0^\eta e^{-rs}Z_s(rY(Z_s) -D(X_s)) ds\bigg]\\
&=
\E_{xz}\bigg[ \int^{\eta}_{0} e^{-rs}Z_s r\big(Y(Z_s) - \theta\big)ds\bigg]
> 0 
\end{align*}
where we used that $Y(Z_s) > \theta$ for $s < \eta$.

It still remains to show $-dZ_s \ge Z_s (rY(Z_s) - D(X_s)) 2/\epsilon\, ds$ for $s \le \eta$. Since $(Z_t)$ is decreasing, it is sufficient to show that the Borel measure induced by the left-hand side dominates the Borel measure induced by the right-hand side for $\P_{xz}$-almost every $\omega$. Hence, it is enough to evaluate those measures $\omega$ by $\omega$ on intervals $(t_1, t_2]$ for $t_2 \le \eta(\omega)$, i.e., to show that
\[
Z_{t_1} - Z_{t_2} \ge \int_{t_1}^{t_2} Z_s \big(rY(Z_s) - D(X_s)\big) \frac2\epsilon\, ds.
\]
For any $\delta \le \epsilon/2$, we have
\begin{align*}
Z^\delta_{t_1} - Z^\delta_{t_2} 
= -(Z^\delta_{t_2} - Z^\delta_{t_1})
&= \int_{t_1}^{t_2} Z^\delta_s \lambda^\delta(X_s, Y(Z^\delta_s)) \ind{m(X_s; Y(Z^\delta_s)) - Y(Z^\delta_s) \le \epsilon/2} ds\\
&\quad + \int_{t_1}^{t_2} Z^\delta_s \lambda^\delta(X_s, Y(Z^\delta_s)) \ind{m(X_s; Y(Z^\delta_s)) - Y(Z^\delta_s) > \epsilon/2} ds\\
&\ge \int_{t_1}^{t_2} Z^\delta_s \big(rY(Z^\delta_s) - D(X_s)\big) \frac{2}{\epsilon} \ind{m(X_s; Y(Z^\delta_s)) - Y(Z^\delta_s) \le \epsilon/2} ds.
\end{align*}
Taking the limit as $\delta \downarrow 0$, we obtain
\[
Z_{t_1} - Z_{t_2}  \ge 
\int_{t_1}^{t_2} Z_s \big(rY(Z_s) - D(X_s)\big) \frac2\epsilon ds,
\]
because $m(X_s; Y(Z^\delta_s)) - Y(Z^\delta_s) \le \epsilon/2$ for all $t \in (t_1, t_2]$ with $t_2 < \eta(\omega)$. 
\end{proof}

\begin{lemma}\label{lem:tau_theta_J}
On the set $J$, the optimal stopping time for $v(x,z; \theta)$ is given by \eqref{eqn:tau_of_theta}.
\end{lemma}
\begin{proof}
Take $(x,z) \in J$, i.e., $z > 0$. Thanks to Theorem \ref{optimalstop}, we have 
\[
v(x, z;\theta) = \max\big(z\theta, v(x, Z_{0+}; \theta) + (z - Z_{0+})m(x; \theta)\big).
\]
If $Z_{0+} = 0$, we have $v(x, Z_{0+}; \theta) = 0$ and also $x \le \beta(\theta_L)$, so $m(x; \theta) = \theta$. Hence, $v(x, z; \theta) = z\theta$ and $\tau(\theta) = 0$ is optimal. 

As $(x, Z_{0+}) \in \overline{\NJ}$, Theorem \ref{optimalstop} applies to $v(x, Z_{0+}; \theta)$ and the form of the postulated optimal stopping time $\tau(\theta)$ follows; indeed, we have $(X_t, Z_{t}) \in \widebar\NJ \cup (\cI \times \{0\})$ for $t > 0$ (cf. Theorem \ref{thm:prop_Z}). It remains to study the case $\tau(\theta) = 0$. We need to show that $v(x, z; \theta) = z\theta$ on $\{\tau(\theta) = 0\}$. If $v(x, Z_{0+}; \theta) + m(x;\theta)(z-Z_{0+}) \le z\theta$, then $v(x, z; \theta) = z\theta$ and $\tau(\theta) = 0$ is clearly optimal. We will demonstrate that $\tau(\theta) > 0$ when $v(x, Z_{0+}; \theta) + m(x;\theta)(z-Z_{0+}) > z\theta$. Assume otherwise, i.e., $v(x, Z_{0+}; \theta) = Z_{0+}\theta$. We split the argument into two cases. When $x \le \beta(\theta)$, then $v(x, Z_{0+}; \theta) + m(x;\theta)(z-Z_{0+}) = v(x, Z_{0+}; \theta) + \theta(z-Z_{0+}) = Z_{0+} \theta + \theta(z-Z_{0+}) = z \theta$, a contradiction. Consider $x > \beta(\theta)$. Since $(x,z)\in J$ and $x>\beta(\theta)\ge\beta(\theta_L)$, we have $x\in(\beta(\theta_L),\beta(\theta_U))$. Theorem \ref{thm:prop_Z} then gives $Z_{0+}=Z(\beta^{-1}(x))$, so $Z_{0+}>Z(\theta)$. Lemma \ref{continuationsetcomp} implies that $0 < \tilde v(x, Z_{0+}; \theta) = v(x, Z_{0+}; \theta) - \theta Z_{0+}$, which leads again to a contradiction. This completes the proof.
\end{proof}

\begin{lemma}\label{alphainclude}
We have
\[
O(\theta) = [\beta(\theta),\alpha(\theta)] \times (0,Z(\theta)] \subseteq S_v(\theta)
\]
for any $\theta \in (\theta_L, \theta_U]$, where $S_v(\theta)$ denotes the stopping set for $v(\cdot, \cdot; \theta)$.
\end{lemma}
\begin{proof}
Fix $\theta \in (\theta_L, \theta_U]$. Recall $\Tilde{v}$ defined in \eqref{tildevexpression} and that the stopping sets for $v$ and $\Tilde v$ coincide. Similarly, $\Tilde{d}(x; \theta) = d(x; \theta) - \theta$ and their stopping sets are equal. The idea of the proof is to show that $\Tilde{v}(x,z; \theta) \le z\Tilde{d}(x; \theta)$ when $(x,z) \in O(\theta)$. Any $x \le \alpha(\theta)$ belongs to the stopping set of $\Tilde{d}(\cdot; \theta)$, so $\Tilde{d}(x; \theta) =0$. Consequently, for any $(x,z) \in O(\theta)$, $\Tilde{v}(x, z; \theta) = 0$ and $(x,z)$ is in the stopping set of $\Tilde{v}(\cdot, \cdot; \theta)$, which, by equality of stopping sets for $\Tilde{v}$ and $v$ means that $(x,z) \in S_v(\theta)$.

Define
\[
R(\theta) := \big\{(x,z) \in \cI \times [0, Z(\theta)]:\ z = 0 \text{ or }  x \ge \beta(Y(z))  \big\}.
\]
The set $R(\theta)$ contains $O(\theta)$. Indeed, from $z \le Z(\theta)$, we have $Y(z) \le \theta$ and $\beta(Y(z)) \le \beta(\theta)$. So, for $(x,z) \in O(\theta)$, we have $x \ge \beta(\theta) \ge \beta(Y(z))$ and $(x,z) \in R(\theta)$. We also note that $R(\theta) \subseteq \widebar\NJ \cup (\cI \times \{0\})$ and it is invariant for the process $(X_t, Z_t)$: if $(x,z) \in R(\theta)$, then $(X_t, Z_t) \in R(\theta)$ for all $t \ge 0$, $\P_{xz}$\as To see it, take $(x,z) \in R(\theta)$. Then $Z_t \le z \le Z(\theta)$ by the monotonicity of $(Z_t)$. Theorem \ref{thm:prop_Z} asserts that $(X_t, Z_t) \in \widebar\NJ \cup (\cI \times \{0\})$ for $t \ge 0$ $\P_{xz}$\as We also have that the process $(Z_t)$ is continuous $\P_{xz}$\as This continuity will be needed below to apply general theory of optimal stopping of \cadlag Markov processes.

We will show that $\phi(x,z)=\Tilde{v}(x,z;\theta)$ is the smallest non-negative function satisfying
\begin{equation}\label{eqn:phi_ineq}
\E_{xz}\big[e^{-rt}\phi(X_t,Z_t) +\int^{t}_0e^{-rs}Z_s(D(X_s)-r\theta)ds-\int_{[0,t]}e^{-rs}(m(X_s,\theta)-\theta)dZ_s \big]\leq\phi(x,z)
\end{equation}
for $(x,z)\in R(\theta)$, while $\phi(x,z) = z\Tilde{d}(x; \theta)$ is a non-negative function satisfying \eqref{eqn:phi_ineq} for $(x,z) \in R(\theta)$. 

Recall the notation in the proof of Theorem \ref{optimalstop}; there, all functions depend implicitly on $\theta$. We recall that $R(\theta)$ is the invariant set of the process $(X_t, Z_t)$. By \cite[Ch.~3, Thm.~1]{shiryaev}, the function $\hat v$ is the smallest excessive majorant of $G$, hence, for any $t \ge 0$, we have $\E_{xz}[e^{-rt}\hat v(X_t, Z_t)] \le \hat v(x, z)$ and $\hat v(x,z) \ge G(x,z)$. Since $\hat v(x,z) = v(x, z; \theta) - F(x,z) = \Tilde{v}(x, z; \theta) + \theta z - F(x,z) = \tilde{v}(x, z; \theta) + G (x,z)$, the inequality $\hat v \ge G$ is equivalent to $\Tilde{v} \ge 0$. The excessive property of $\hat v$ translates into
\begin{equation}\label{eqn:exc1}
\E_{xz} \big[ e^{-rt} \Tilde{v}(X_t, Z_t; \theta) + e^{-rt} G(X_t, Z_t) - G(x,z) \big] \le \Tilde{v}(x,z; \theta).
\end{equation}
Using the definition of $G$
\[
G(x,z) = \theta z - \E_{xz} \Big[-\int_0^\infty e^{-rs}Z_sD(X_s)ds + \int_{[0,\infty)}e^{-rs}m(X_s;\theta)dZ_s\Big],
\]
the inequality \eqref{eqn:exc1} takes the form
\begin{align*}
&\Tilde{v}(x,z; \theta)\\
&\ge
\E_{xz} \Big[ e^{-rt} \Tilde{v}(X_t, Z_t; \theta) + e^{-rt} \theta Z_t - \theta z +\int_0^t e^{-rs}Z_sD(X_s)ds - \int_{[0,t]}e^{-rs}m(X_s;\theta)dZ_s \Big]\\
&=
\E_{xz} \Big[ e^{-rt} \Tilde{v}(X_t, Z_t; \theta) + \int_0^t e^{-rs}Z_s(D(X_s)-r\theta)ds - \int_{[0,t]}e^{-rs}\big(m(X_s;\theta)-\theta\big)dZ_s \Big],
\end{align*}
which shows that $\tilde{v}$ is the smallest non-negative function satisfying \eqref{eqn:phi_ineq} for $(x,z) \in R(\theta)$.

From \eqref{eqn:a1} and $\Tilde d(x; \theta) = d(x; \theta) - \theta$, we obtain
\begin{equation}\label{eqn:tldd}
(\cL - r) \Tilde d (x; \theta) + D(x) - r\theta \le 0.
\end{equation}
Notice that $\Tilde d$ is continuously differentiable (Proposition \ref{prop:dm})  with the second derivative potentially discontinuous at the stopping boundary $\alpha(\theta)$. However, it is locally bounded there as for $x > \alpha(\theta)$, the second derivative is determined by $(\cL-r) \Tilde d(x; \theta) = 0$ and for $x < \alpha(\theta)$ we have $\Tilde{d}(x; \theta) = 0$, so the second derivative is zero. This means that we can apply It\^o formula (\cite[Ch.~3, Problem~7.3]{karatzas2012brownian}) and obtain
\begin{align*}
d\big(e^{-rt} \Tilde{d}(X_t; \theta)\big)
&=
e^{-rt} (\cL - r) \Tilde{d}(X_t; \theta) dt + dM_t,
\end{align*}
where $(M_t)$ is a local martingale. Further, using a product rule we obtain
\begin{align*}
&\E_{xz} \big[ e^{-r(t \wedge \tau_n)} Z_{t \wedge \tau_n} \Tilde{d}(X_{t \wedge \tau_n}; \theta) - z \Tilde{d}(x; \theta) \big]\\
&=
\E_{xz} \Big[ \int_0^{t \wedge \tau_n} e^{-rs} Z_s (\cL - r) \Tilde{d}(X_s; \theta) ds + \int_{[0, t \wedge \tau_n]} e^{-rs} \Tilde{d}(X_s; \theta) dZ_s \Big]
\end{align*}
where $(\tau_n)$ is the localising sequence for $(M_t)$. We split the first term as
\begin{align*}
&\int_0^{t \wedge \tau_n} \ind{X_s > \alpha(Y(Z_s))} e^{-rs} Z_s (\cL - r) \Tilde{d}(X_s; \theta) ds\\
&+
\int_0^{t \wedge \tau_n} \ind{X_s \le \alpha(Y(Z_s))} e^{-rs} Z_s (\cL - r) \Tilde{d}(X_s; \theta) ds\\
&\le \int_0^{t \wedge \tau_n} \ind{X_s > \alpha(Y(Z_s))} e^{-rs} Z_s \big(r\theta - D(X_s)\big) ds,
\end{align*}
where we used that $\tilde d(x; \theta) = 0$ for $x \le \alpha(Y(Z_s)) \le \alpha(Y(z)) \le \alpha(\theta)$ in the case of the second term and \eqref{eqn:tldd} for the first term. Hence, we have
\begin{align*}
&\E_{xz} \big[ e^{-r(t \wedge \tau_n)} Z_{t \wedge \tau_n} \Tilde{d}(X_{t \wedge \tau_n}; \theta) - z \Tilde{d}(x; \theta) \big]\\
&\le
\E_{xz} \Big[ \int_0^{t \wedge \tau_n} \ind{X_s > \alpha(Y(Z_s))} e^{-rs} Z_s \big(r\theta - D(X_s)\big) ds+ \int_{[0, t \wedge \tau_n]} e^{-rs} \Tilde{d}(X_s; \theta) dZ_s \Big].
\end{align*}
Thanks to the boundedness of $\Tilde{d}$, $D$ and $(Z_t)$, we apply the dominated convergence theorem in order to pass to the limit $n \to \infty$. Rearranging, this yields
\[
z \Tilde{d}(x; \theta) \ge 
\E_{xz} \Big[ e^{-rt} Z_{t} \Tilde{d}(X_{t}; \theta) 
+ \int_0^{t} \ind{X_s > \alpha(Y(Z_s))} e^{-rs} Z_s \big(D(X_s) - r\theta\big) ds - \int_{[0, t]} e^{-rs} \Tilde{d}(X_s; \theta) dZ_s\Big].
\]
Consider the last term and split it as follows
\[
\int_{[0, t]} \ind{X_s > \alpha(Y(Z_s))} e^{-rs} \Tilde{d}(X_s; \theta)  dZ_s
+
\int_{[0, t]} \ind{X_s \le \alpha(Y(Z_s))} e^{-rs} \Tilde{d}(X_s; \theta)  dZ_s.
\]
On the set $\{X_s > \alpha(Y(Z_s))\}$, the process $(Z_t)$ is constant so the first integral on the right-hand side is zero and we can replace $\Tilde{d}(X_s; \theta)$ with $m(X_s; \theta) - \theta$. On the other hand, when $X_s \le \alpha(Y(Z_s)) \le \alpha(Y(z)) \le \alpha(\theta)$, we have $\Tilde{d}(X_s; \theta) = 0$, so the second integral is zero too. This leads to a final expression
\begin{equation}\label{eqn:z-tilde-d}
\begin{aligned}
z \Tilde{d}(x; \theta) &\ge 
\E_{xz} \Big[ e^{-rt} Z_{t} \Tilde{d}(X_{t}; \theta)
+ \int_0^{t}  \ind{X_s > \alpha(Y(Z_s))} e^{-rs} Z_s \big(D(X_s)- r\theta\big)  ds\\
&\hspace{33pt} - \int_{[0, t]}  \ind{X_s > \alpha(Y(Z_s))} e^{-rs} \big(m(X_s; \theta) - \theta\big)  dZ_s\Big].
\end{aligned}
\end{equation}
To show that $\phi(x,z) = z\Tilde{d}(x; \theta)$ satisfies \eqref{eqn:phi_ineq} for $(x,z) \in R(\theta)$, it remains to prove that
\begin{equation}\label{eqn:integ}
\int_0^{t}  e^{-rs} Z_s \big(D(X_s)- r\theta\big) \ind{X_s \le \alpha(Y(Z_s))} ds
- \int_{[0, t]}  e^{-rs} \big(m(X_s; \theta) - \theta\big) \ind{X_s \le \alpha(Y(Z_s))} dZ_s \le 0.
\end{equation}
The left-hand side is trivially zero $\P_{xz}$\as for $z = 0$. Take $(x,z) \in R(\theta)$. Recall that $(X_t, Z_t) \in R(\theta)$ $\P_{xz}\as$, so $X_t \ge \beta(Y(Z_t))$ as long as $Z_t > 0$; when $Z_t$ reaches $0$ it is absorbed and the remaining integration yields $0$. Consider a sequence $\delta_n \downarrow 0$ with $\delta_0 = x_U - x_L$; then $\beta(Y(Z_s))+\delta_0 \ge \alpha(Y(Z_s))$.  We have
\begin{equation}\label{eqn:integs}
\begin{aligned}
&\int_0^{t}  e^{-rs} Z_s \big(D(X_s)- r\theta\big) \ind{X_s \le \alpha(Y(Z_s))} ds
- \int_{[0, t]}  e^{-rs} \big(m(X_s; \theta) - \theta\big) \ind{X_s \le \alpha(Y(Z_s))} dZ_s\\
&=
\sum_{i=1}^\infty \bigg(\int_0^{t}  e^{-rs} Z_s \big(D(X_s)- r\theta\big) \ind{(X_s,Z_s) \in A_i} ds
- \int_{[0, t]}  e^{-rs} \big(m(X_s; \theta) - \theta\big) \ind{(X_s,Z_s) \in A_i} dZ_s \bigg)\\
&\quad+ \int_0^{t}  e^{-rs} Z_s \big(D(X_s)- r\theta\big) \ind{\beta(Y(Z_s)) = X_s} ds
- \int_{[0, t]}  e^{-rs} \big(m(X_s; \theta) - \theta\big) \ind{\beta(Y(Z_s)) = X_s} dZ_s,
\end{aligned}
\end{equation}
where
\[
A_i = \{(x,z) \in \I \times [0,1]\colon \beta(Y(z)) + \delta_i < x \le \big(\beta(Y(z)) + \delta_{i-1}\big) \wedge \alpha(Y(z))\}.
\]
By \cite[Lemma~3]{kwon}, we have $\int_0^{t}  e^{-rs} Z_s \big(D(X_s)- r\theta\big) \ind{\beta(Y(Z_s)) = X_s} ds = 0$. Observe that  when $\beta(Y(Z_s)) = X_s$, we also have $X_s \le \beta(\theta)$, so $m(X_s;\theta)-\theta=0$ and $\int_{[0, t]}  e^{-rs} \big(m(X_s; \theta) - \theta\big) \ind{\beta(Y(Z_s)) = X_s} dZ_s = 0$. 

Turning attention to the elements of the sum in \eqref{eqn:integs}, $(Z_t)$ follows the dynamics \eqref{eqn:dynamics_Z} on $\{(X_s, Z_s) \in A_i\}$ (see Lemma \ref{lem:dynamics_Z}). Hence
\begin{align*}
&-\int_{[0, t]}  e^{-rs} \big(m(X_s; \theta) - \theta\big) \ind{(X_s,Z_s) \in A_i} dZ_s\\
&=
\int_{[0, t]}  e^{-rs} \big(m(X_s; \theta) - \theta\big) \ind{(X_s,Z_s) \in A_i} Z_s \lambda(X_s, Y(Z_s))ds\\
&\le
\int_{[0, t]}  e^{-rs} Z_s \big(r \theta - D(X_s)\big) \ind{(X_s,Z_s) \in A_i} ds
\end{align*}
where the inequality follows from the following bound: for $(x', z') \in A_i$, we either have $x' \le \beta(\theta)$ when
\[
\lambda(x',Y(z'))(m(x';\theta)-\theta) = 0 \le r\theta - D(x'),
\]
where the final inequality follows from Proposition \ref{prop:dm} and $\beta(\theta) \le \alpha(\theta)$;
or $x' > \beta(\theta)$ when
\begin{equation*}
\lambda(x',Y(z'))(m(x';\theta)-\theta)\le \lambda(x',\theta)(m(x';\theta)-\theta)=r\theta - D(x'),
\end{equation*}
by the monotonicity of $y\mapsto\lambda(x,y)$ (by Lemma \ref{lambdadeltamontone}) and the definition \eqref{eqn:lambda} of $\lambda$.
Hence, each element of the sum in \eqref{eqn:integs} is bounded from above by $0$, which completes the proof of \eqref{eqn:integ}.
\end{proof}

\begin{proof}[Proof of Theorem \ref{thm:stopping-set}]
Fix $\theta \in (\theta_L, \theta_U)$. Denote by $S_v(\theta)=\{(x,z) \in \cI \times [0,1]\colon v(x,z;\theta)=z\theta\}$ the stopping set for $v$. By Theorem \ref{optimalstop} and Lemma \ref{lem:tau_theta_J}, its hitting time \eqref{eqn:tau_of_theta} is optimal. Lemmas \ref{optimaliszero} and \ref{alphainclude} as well as an obvious inclusion $\cI \times \{0\} \subset S_v(\theta)$ imply that $S(\theta) \subseteq S_v(\theta)$. It remains to show that the complement of $S(\theta)$ has an empty intersection with $S_v(\theta)$, i.e., is contained in the continuation set $C_v(\theta)$. We have
\begin{align*}
S^c(\theta) &= \Big( (\alpha(\theta), x_U) \times (0, Z(\theta)]\Big) \cup \Big( (\beta(\theta), x_U) \times (Z(\theta), 1] \Big)\\
&= \Big( (\alpha(\theta), x_U) \times (0,1] \Big) \cup \Big( (\beta(\theta), \alpha(\theta)] \times (Z(\theta), 1] \Big).
\end{align*}
It is shown in Lemma \ref{conti} that $(\alpha(\theta), x_U) \times (0,1] \subset C_v(\theta)$. A combination of Lemma \ref{lem7.3} and 
Lemma \ref{continuationsetcomp} demonstrates that $(\beta(\theta), \alpha(\theta)] \times (Z(\theta), 1] \subset C_v(\theta)$. 
This concludes the proof.
\end{proof}

\section{Numerical Illustration}\label{sec:example}
\begin{figure}[tb]
    \centering
    \includegraphics[width=0.6\textwidth]{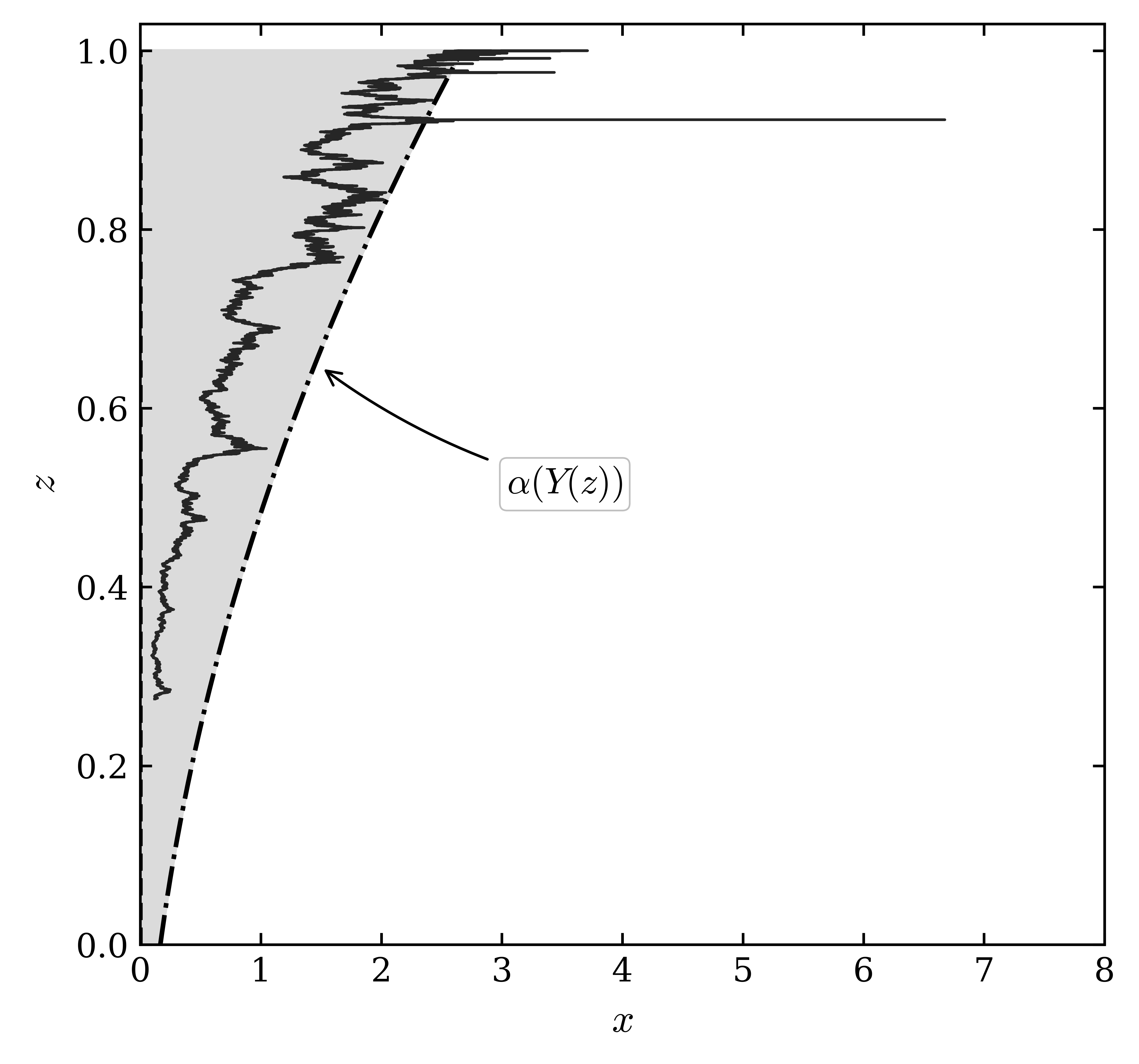}
    \caption{A plot of $(X_t, Z_t)$ when $\theta_L = 0.5$ and $\theta_U = 2$.}
    \label{fig:agreecvae}
\end{figure}

We conclude the paper with a presentation of a specific example which we used for Figure \ref{fig:example_plot}. The underlying process is a geometric Brownian motion
\[
dX_t= -0.5 X_t dt+ X_t dW_t
\]
with the state space $\I = (0, \infty)$. The interest rate is $r=1$. The duopoly payoff takes the form
\[
D(x) = \begin{cases}
x^{1/2}, & x \in (0, x_M),\\
x_M^{1/2}, & x \ge x_M,
\end{cases}
\]
while the monopoly payoff is $M(x) = D(x) + 2$. We study two cases:
\begin{itemize}
\item $\theta_i$ uniformly distributed over $[0.5, 5]$, when the monopolist may exit the market for larger values of $\theta$ (see Figure \ref{fig:example_plot}); 
\item $\theta_i$ uniformly distributed over $[0.5, 2]$, when the monopolist never exits the market reducing our problem to the setting studied in \cite{kwon} (see Figure \ref{fig:agreecvae}).
\end{itemize}
It can be shown using results in \cite[Appendix~B]{kwon} that Assumption \ref{ass:dm} is satisfied for sufficiently large $x_M$. We verified numerically that $x_M = 1000$ is large enough for our model.

When $\theta_i$ is uniformly distributed on $[0.5, 2]$, the monopolist never exits the market. The present value of the future profit flow always dominates the exit value. The boundary $\beta(\theta) = x_L = 0$ for all $\theta \in [0.5, 2]$, see Figure \ref{fig:agreecvae}. This is exactly the setting studied in \cite{kwon}. To enable comparisons between problems with different ranges of $\theta$, the trajectory of the process $(X_t)$ is the same as in Figure \ref{fig:example_plot}. Notice that the process $(Z_t)$ decreases faster when $\theta_U$ -- the upper end of the range of exit values -- increases. This has two reasons. Firstly, the duopoly boundary is shifted to the right as the mapping $\theta \mapsto \alpha(\theta)$ is increasing and $z \in [0,1]$ is mapped into $[\theta_L, \theta_U]$. Secondly, the appearance of the monopoly boundary for a larger value of $\theta_U$ puts an additional pressure increasing the rate at which the process $(Z_t)$ decreases. Indeed, the rate explodes when the process $(X_t, Z_t)$ approaches the monopoly boundary $(\beta(Y(z)), z)_{z \in [0,1]}$ while in the framework of \cite{kwon} and Figure \ref{fig:agreecvae} the rate is bounded. 

\begin{figure}[tb]
    \centering
    \includegraphics[width=0.6\textwidth]{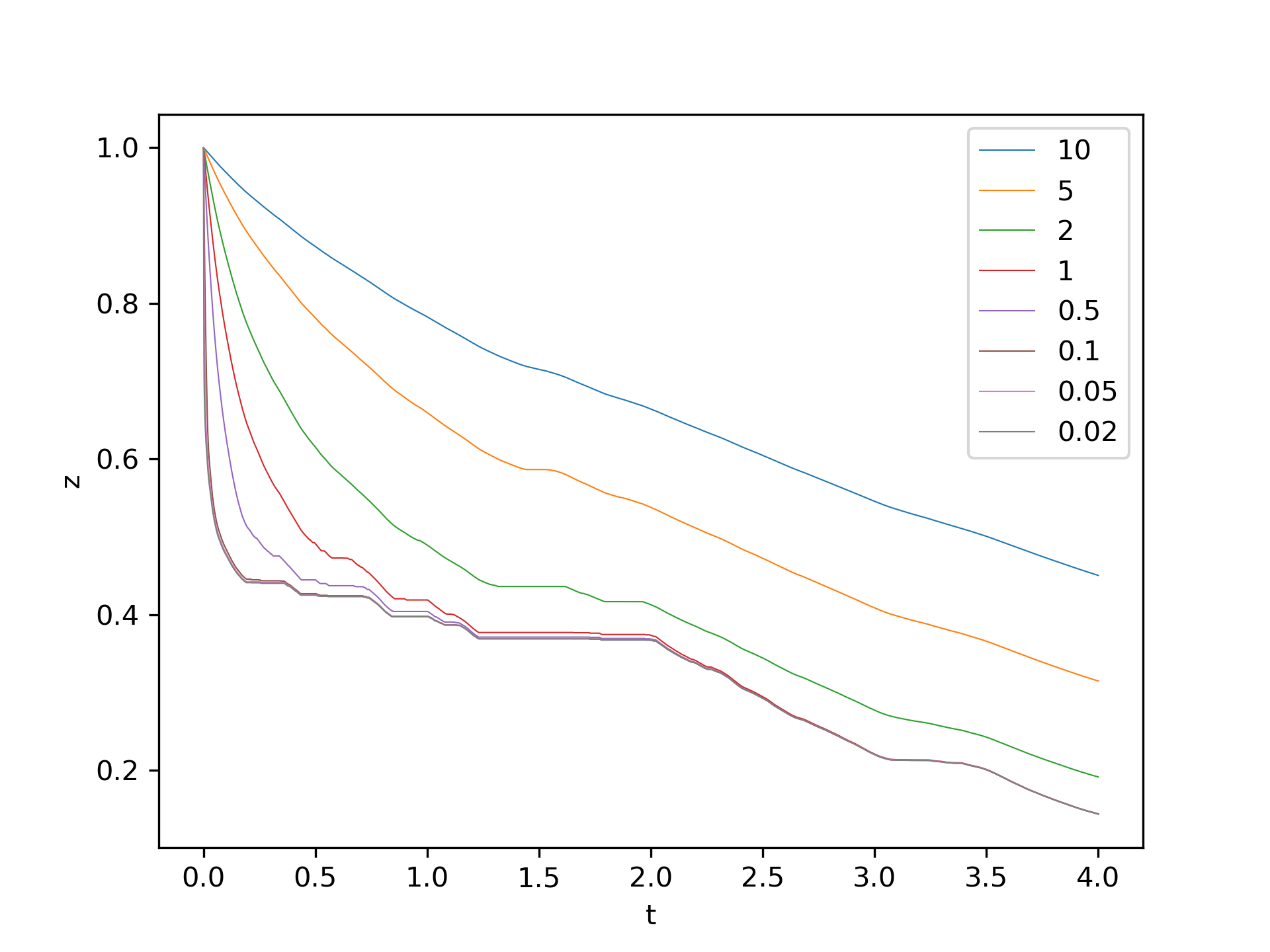}
    \caption{Evolution of $(Z^\delta_t)$ with different values of $\delta$ and a common trajectory of $(X_t)$.}
    \label{fig:differentdelta}
\end{figure}

The process $(Z_t)$ was constructed in Section \ref{sec:belief} as a limit of approximations $(Z^\delta_t)$ as $\delta \downarrow 0$. Figure \ref{fig:differentdelta} shows how these approximations depend on $\delta$; the graph is based on a single trajectory of the process $(X_t)$. Note that the limit process $(Z_t)$ exhibits a jump at time $0$.

\bibliographystyle{siam}
\bibliography{exitgame}

\end{document}